\documentclass[aos,preprint]{imsart}

\usepackage{graphicx} 

\usepackage{amsfonts,amssymb,float}

\DeclareSymbolFont{sfletters}{OML}{cmbrm}{m}{it}

\DeclareMathSymbol{\salpha}{\mathord}{sfletters}{"0B}

\usepackage{mathrsfs}
\RequirePackage{natbib}
\RequirePackage[colorlinks,citecolor=blue,urlcolor=blue]{hyperref}

\usepackage{xr}
\usepackage{xcite}

\externalcitedocument{boot_LLM_aos_galley_SUPP_18july2019}
\usepackage{mathtools}
\usepackage{subfigure} 
\RequirePackage{natbib}
\usepackage{comment}
\usepackage{enumerate}
\usepackage{eufrak}

\usepackage{array,varwidth}

\RequirePackage[OT1]{fontenc}
\RequirePackage{amsthm,amsmath}

\startlocaldefs
\numberwithin{equation}{section}
\theoremstyle{plain}
\newtheorem{theorem}{Theorem}[section]
\newtheorem{lemma}{Lemma}[section]
\newtheorem{proposition}{Proposition}[section]
\newtheorem{assumption}{Assumption}[section]

\newcommand{\M}{\mathcal{M}}

\newcommand{\R}{\mathbb{R}}

\renewcommand{\hat}[1]{\widehat{#1}}

\newcommand{\I}{\textsc{I}}
\newcommand{\II}{\textsc{II}}
\newcommand{\III}{\textsc{III}}
\newcommand{\bI}{\textup{\texttt{I}}}
\newcommand{\bII}{\textup{\texttt{II}}}
\newcommand{\bIII}{\textup{\texttt{III}}}

\newcommand{\dK}{d_{\textup{K}}}

\newcommand{\J}{\mathcal{J}}

\newcommand{\G}{\mathbb{G}}

\renewcommand{\P}{\mathbb{P}}
\newcommand{\E}{\mathbb{E}}

\renewcommand{\L}{\mathcal{L}}

\newcommand{\e}{\epsilon}

\newcommand{\tr}{\operatorname{tr}}

\newcommand{\var}{\operatorname{var}}

\newcommand{\ttop}{^{\top}}

\newcommand{\op}{_{\textup{op}}}

\newcommand{\tsum}{\textstyle\sum}
\newcommand{\ts}{\textstyle}

\newcommand{\diag}{\textup{diag}}

\newcommand{\hilbert}{\mathcal{H}}

\newcommand{\sci}{\hat{\mathcal{I}}}

\newcommand{\siglevel}{\varrho}
\newcommand{\covarop}{\mathcal{C}}

\endlocaldefs

\def\ci{\cite}
\def\cp{\citep}
\def\references{\bibliography{decay_boot_bib_galley,fda_galley}}

\begin{document}

\begin{frontmatter}
\title{Bootstrapping Max Statistics in High Dimensions: Near-Parametric Rates Under Weak Variance Decay and Application to Functional and Multinomial Data}
\runtitle{Bootstrapping Max Statistics}

\begin{aug}
\author{\fnms{Miles E.} \snm{Lopes}\thanksref{t1}\ead[label=e1]{melopes@ucdavis.edu}},
\author{\fnms{Zhenhua} \snm{Lin}\thanksref{t2}\ead[label=e2]{zhnlin@ucdavis.edu}}
\and
\author{\fnms{Hans-Georg} \snm{M\"uller}\thanksref{t3}
\ead[label=e3]{hgmueller@ucdavis.edu}}

\thankstext{t1}{Supported in part by NSF grant DMS 1613218}
\thankstext{t2}{Supported in part by NIH grant  5UG3OD023313-03}
\thankstext{t3}{Supported in part by NSF grant DMS 1712864 and NIH grant 5UG3OD023313-03}
\runauthor{M.~E. Lopes et al.}

\affiliation{University of California, Davis}

\end{aug}



\begin{abstract}

In recent years, bootstrap methods have drawn attention for their ability to
approximate the laws of ``max statistics'' in high-dimensional problems.
A leading example of such a statistic is the coordinate-wise maximum of a sample average of $n$ random vectors in $\R^p$. Existing results for this statistic show that the bootstrap can work when $n\ll p$, and rates of approximation (in Kolmogorov distance) have been obtained with only logarithmic dependence in $p$.
Nevertheless, one of the challenging aspects of this setting is that established rates tend to scale like $n^{-1/6}$ as a function of $n$.

The main purpose of this paper is to demonstrate that improvement in rate is possible when extra model structure is available. Specifically, we show that if the coordinate-wise variances of the observations exhibit decay, then a nearly $n^{-1/2}$ rate can be achieved, \emph{independent of $p$}. Furthermore, a surprising aspect of this dimension-free rate is that it holds even when the decay is \emph{very weak}. Lastly, we provide examples showing how these ideas can be applied to inference problems dealing with functional and multinomial data.

\end{abstract}



\end{frontmatter}

\section{Introduction}
One of the current challenges in theoretical statistics is to understand when bootstrap methods work in high-dimensional problems. In this direction, there has been a surge of recent interest in connection with ``max statistics'' such as
$$T=\max_{1\leq j\leq p} S_{n,j},$$
where $S_{n,j}$ is the $j$th coordinate of the sum \mbox{$S_n=\frac{1}{\sqrt n}\sum_{i=1}^n (X_i-\E[X_i])$,} involving i.i.d.~vectors $X_1,\dots,X_n$ in $\R^p$.

This type of statistic has been a focal point in the literature for at least two reasons. First, it is an example of a statistic for which bootstrap methods can succeed in high dimensions under mild assumptions, which was established 
in several pathbreaking works~\cp{ArlotI:2010,ArlotII:2010,CCK:2013,CCK:AOP}. Second, the statistic $T$ is closely linked to several fundamental topics, such as suprema of empirical processes, nonparametric confidence regions, and multiple testing problems. Likewise, many applications of bootstrap methods for max statistics have ensued at a brisk pace in recent years 
\cp[see, e.g.,][]{CCKAnti:2014,Wasserman:2014,ChenRidges:2015,Yao:2017,ChengJASA:2017,Dezeure:2017,ChenUstat:2018,fan:18,Belloni:2018}.

One of the favorable aspects of bootstrap approximation results for the distribution $\mathcal{L}(T)$ is that rates have been established with only logarithmic dependence in $p$. For instance, the results in~\cite{CCK:AOP} imply that under certain conditions, the Kolmogorov distance $\dK$ between $\mathcal{L}(T)$ and its bootstrap counterpart $\mathcal{L}(T^*|X)$ satisfies the bound
\begin{equation}\label{eqn:introrate}
\dK\Big( \L(T) \, , \L(T^*|X)\Big) \ \leq  \frac{c \log(p)^{b}}{n^{1/6}}
\end{equation}
with high probability, where $c,b>0$ are constants not depending on $n$ or $p$, and $X$ denotes the matrix whose rows are $X_1,\dots,X_n$. (In the following, symbols such as $c$ will be often re-used to designate a positive constant not depending on $n$ or $p$, possibly with a different value at each occurrence.) Additional refinements of this result can be found in the same work, with regard to the choice of metric, or choice of bootstrap method. Also, recent progress in sharpening the exponent $b$ has been made by~\cite{Deng:2017}. However, this mild dependence on $p$ is offset by the $n^{-1/6}$ dependence on $n$, which differs from the $n^{-1/2}$ rate in the multivariate Berry-Esseen theorem when $p\ll n$.

Currently, the general problem of determining the best possible rates for Gaussian and bootstrap approximations is largely open in the high-dimensional setting. In particular, if we let $\tilde T$ denote the counterpart of $T$ that arises from replacing $X_1,\dots,X_n$ with independent Gaussian vectors \mbox{$\tilde X_1,\dots,\tilde X_n$} satisfying $\text{cov}(X_i)=\text{cov}(\tilde X_i)$, then a conjecture of~\cite{CCK:AOP} indicates that a bound of the form $d_{\text{K}}(\L(T),\L(\tilde T))\leq cn^{-1/6}\log(p)^b$ is optimal under certain conditions. A related conjecture in the setting of high-dimensional U-statistics may also be found in~\cite{ChenUstat:2018}. (Further discussion of related work on Gaussian approximation is given in Appendix~\ref{app:Gaussian}.)
Nevertheless, the finite-sample performance of bootstrap methods for max statistics is often more encouraging than what might be expected from the $n^{-1/6}$ dependence on $n$~\cite[see, e.g.][]{ChengJASA:2017, fan:18,Belloni:2018}. This suggests that improved rates are possible in at least some situations. 

The purpose of this paper is to quantify an instance of such improvement when additional model structure is available. Specifically, we consider the case when the coordinates of $X_1,\dots,X_n$ have decaying variances. If we let $\sigma_j^2=\var(X_{1,j})$ for each $1\leq j\leq p$, and write $\sigma_{(1)}\geq \dots\geq \sigma_{(p)}$, then this condition may be formalized as
\begin{equation}\label{eqn:varcondintro}
\sigma_{(j)} \leq c \, j^{-\alpha} \ \ \ \ \text{for all} \ \ \ j\in\{1,\dots,p\},
\end{equation}
where $\alpha>0$ is a parameter not depending on $n$ or $p$. (A complete set of assumptions, including a weaker version of~\eqref{eqn:varcondintro}, is given in Section~\ref{sec:prelim}.) This type of condition arises in many contexts, and in Section~\ref{sec:prelim} we discuss examples related to principal component analysis, count data, and Fourier coefficients of functional data. Furthermore, this condition can be assessed in practice, due to the fact that the parameters $\sigma_1,\dots,\sigma_p$ can be accurately estimated, even in high dimensions (cf.~Lemma~\ref{lem:corbasic}).

Within the setting of decaying variances, our main results show that a nearly parametric rate can be achieved for both Gaussian and bootstrap approximation of $\L(T)$. More precisely, this means that for any fixed \mbox{$\delta\in(0,1/2)$}, the bound $\dK( \L(T) , \L(\tilde T))\leq cn^{-1/2+\delta}$ holds, and similarly, the event
\begin{equation}\label{eqn:improved}
\dK\Big( \L(T) \, , \L(T^*|X)\Big)  \ \leq  \ c\, n^{-1/2+\delta}
\end{equation}
holds with high probability.
Here, it is worth emphasizing a few basic aspects of these bounds. First, they are non-asymptotic and \emph{do not depend on $p$}. Second, the parameter $\alpha$ is allowed to be \emph{arbitrarily small}, and in this sense, the decay condition~\eqref{eqn:varcondintro} is very weak. Third, the result for $T^*$ holds when it is constructed using the standard multiplier bootstrap procedure~\citep{CCK:2013}.

With regard to the existing literature, it is important to clarify that our near-parametric rate does not conflict with the conjectured optimality of the rate $n^{-1/6}$ for Gaussian approximation. The reason is that the $n^{-1/6}$ rate has been established in settings where the values $\sigma_1,\dots,\sigma_p$ are restricted from becoming too small. A basic version of such a requirement is that
\begin{equation}\label{eqn:lower}
\min_{1\leq j\leq p} \sigma_j \geq c.
\end{equation}
Hence, the conditions~\eqref{eqn:varcondintro} and~\eqref{eqn:lower} are complementary. Also, it is interesting to observe that the two conditions ``intersect'' in the limit $\alpha\to 0^+$, suggesting there is a phase transition in rates at the ``boundary'' corresponding to $\alpha=0$.

Another important consideration that is related to the conditions~\eqref{eqn:varcondintro} and~\eqref{eqn:lower} is the use of standardized variables. Namely, it is of special interest to approximate the distribution of the statistic
$$T'=\max_{1\leq j\leq p} S_{n,j}/\sigma_j,$$
which is equivalent to approximating $\L(T)$ when each $X_{i,j}$ is standardized to have variance 1. Given that standardization eliminates variance decay, it might seem that the rate $n^{-1/2+\delta}$ has no bearing on approximating $\L(T')$. However, it is still possible to take advantage of variance decay, by using a basic notion that we refer to as ``partial standardization''.

The idea of partial standardization is to slightly modify $T'$ by using a fractional power of each $\sigma_j$. Specifically, if we let $\tau_n\in [0,1]$ be a free parameter, then we can consider the partially standardized statistic
\begin{equation}\label{eqn:mdef}
M=\max_{1\leq j\leq p} S_{n,j}/\sigma_j^{\tau_n},
\end{equation}
which interpolates between $T$ and $T'$ as $\tau_n$ ranges over $[0,1]$.
This statistic has the following significant property: If $X_1,\dots,X_n$ satisfy the variance decay condition~\eqref{eqn:varcondintro}, and if $\tau_n$ is chosen to be slightly less than 1, then our main results show that the rate $n^{-1/2+\delta}$ holds for bootstrap approximations of $\L(M)$. In fact, this effect occurs even when $\tau_n\to 1$ as $n\to\infty$. Further details can be found in Section~\ref{sec:main}. Also note that our main results are formulated entirely in terms of $M$, which covers the statistic $T$ as a special case.

In practice, simultaneous confidence intervals derived from approximations to $\L(M)$ are just as easy to use as those based on $\L(T')$. Although there is a slight difference between the quantiles of $M$ and $T'$ when $\tau_n<1$, the important point is that the quantiles of $\L(M)$ may be preferred, since faster rates of bootstrap approximation are available. (See also Figure~\ref{fig:tau} in Section~\ref{sec:expt}.)
In this way, the statistic $M$ offers a simple way to blend the utility of standardized variables with the beneficial effects of variance decay.

\paragraph{Outline} 
The remainder of the paper is organized as follows. In Section~\ref{sec:prelim}, we outline the problem setting, with a complete statement of the theoretical assumptions, as well as some motivating facts and examples. Our main results are given in Section 3, which consist of a Gaussian approximation result for $\L(M)$ (Theorem~\ref{THM:G}), and a corresponding bootstrap approximation result (Theorem~\ref{THM:BOOT}). To provide a numerical illustration of our results, we discuss  a problem in functional data analysis in Section~\ref{sec:expt},  where the variance decay condition naturally arises. Specifically, we show how bootstrap approximations to $\L(M)$ can be used to derive simultaneous confidence intervals for the Fourier coefficients of a mean function. A second application to high-dimensional multinomial models is described in Section~\ref{sec:multinomial}, which offers both a theoretical bootstrap approximation result, as well as some numerical results.
  Lastly, our conclusions are summarized in Section~\ref{sec:conc}. All proofs are given in the appendices, found in the supplementary material.

\paragraph{Notation} 
 The standard basis vectors in $\R^p$ are denoted $e_1,\dots,e_p$, and the identity matrix of size $p\times p$ is denoted  $\mathbf{I}_p$. For any symmetric matrix $A\in\R^{p\times p}$, the ordered eigenvalues are denoted $\lambda(A)=(\lambda_1(A),\dots,\lambda_p(A))$, where $\lambda_{\max}(A)=\lambda_1(A)\geq \cdots \geq \lambda_p(A)=\lambda_{\min}(A)$. The operator norm of a matrix $A$, denoted $\|A\|_{\text{op}}$, is the same as its largest singular value.
If $v\in\R^p$ is a fixed vector, and $r>0$, we write $\|v\|_r=(\sum_{j=1}^p |v_j|^r)^{1/r}$. In addition, the weak-$\ell_r$ (quasi) norm is given by $\|v\|_{w\ell_r}=\max_{1\leq j\leq p} j^{1/r}|v|_{(j)},$ where \mbox{$|v|_{(1)}\geq \cdots \geq |v|_{(p)}$} are the sorted absolute entries of $v$. Likewise, the notation $v_{(1)}\geq\cdots\geq v_{(p)}$ refers to the sorted entries.
In  a slight abuse of notation, we write $\|\xi\|_r=\E[|\xi|^r]^{1/r}$ to refer to the $L^r$ norm of a scalar random variable $\xi$, with $r\geq 1$. The $\psi_1$-Orlicz norm is $\|\xi\|_{\psi_1}=\inf\{t>0\, | \, \E[\exp(|\xi|/t)]\leq 2\}$.  If $\{a_n\}$ and $\{b_n\}$ are sequences of non-negative real numbers, then the relation $a_n\lesssim b_n$ means that there is a constant $c>0$ not depending on $n$, and an integer $n_0\geq 1$, such that $a_n\leq c b_n$ for all $n\geq n_0$.  Also, we write $a_n\asymp b_n$ if $a_n\lesssim b_n$ and $b_n\lesssim a_n$.
Lastly, define the abbreviations $a_n\vee b_n=\max\{a_n,b_n\}$ and $a_n\wedge b_n=\min\{a_n,b_n\}$.

\section{Setting and preliminaries}\label{sec:prelim}

We consider a sequence of models indexed by $n$, with all parameters depending on $n$, except for those that are stated to be fixed. In particular, the dimension $p=p(n)$ is regarded as a function of $n$, and hence, if a constant does not depend on $n$, then it does not depend on $p$ either.

\noindent \begin{assumption}[Data-generating model]\label{A:model}
~\\[-0.3cm]
\begin{enumerate}[(i).]
\item There is a vector $\mu=\mu(n)\in\R^p$ and positive semi-definite matrix $\Sigma=\Sigma(n)\in\R^{p\times p}$, such that the observations $X_1,\dots,X_n\in\R^p$ are generated as $X_i=\mu+\Sigma^{1/2}Z_i$ for each $1\leq i\leq n$, where the random vectors $Z_1,\dots,Z_n\in\R^p$ are i.i.d.\\[-0.2cm]
\item The random vector $Z_1$ satisfies $\E[Z_1]=0$ \!and $\E[Z_1Z_1\ttop]=\mathbf{I}_p$, as well as $\sup_{\|u\|_2=1}\|Z_1\ttop u\|_{\psi_1}\leq c_0$, for some constant $c_0>0$ that does not depend on $n$.
\end{enumerate}
\end{assumption}
\paragraph{Remarks} Note that no constraints are placed on the ratio $p/n$.
Also, the sub-exponential tail condition in part \emph{(ii)} is similar to other tail conditions that have been used in previous works on bootstrap methods for max statistics \cp[][]{CCK:2013,Deng:2017}.

 To state our next assumption, it is necessary to develop some notation. For any $d\in\{1,\dots,p\}$, let $\J(d)$  denote a set of indices corresponding to the $d$ largest values among $\sigma_1,\dots,\sigma_p$, i.e., 
 $\{\sigma_{(1)},\dots,\sigma_{(d)}\}= \{\sigma_j| \ j\in\J(d)\}$. In addition, let $R(d)\in\R^{d\times d}$ denote the correlation 
 matrix of the random variables $\{X_{1,j}\,|\, j\in\J(d)\}$.
Lastly, let $a\in(0,\ts\frac{1}{2})$ be a constant fixed with respect to $n$, and define the integers $\ell_n$ and $k_n$ according to
\begin{align}\
\ell_n&=\big\lceil (1\vee \log(n)^3)\wedge p\big\rceil\\[0.2cm]
k_n&=\big\lceil \big(\ell_n\vee n^{\frac{1}{\log(n)^{a}}}\big)\wedge p\big\rceil.\label{eqn:kndef}
\end{align}
Note that both $\ell_n$ and $k_n$ grow slower than any fractional power of $n$, and always satisfy $1\leq  \ell_n\leq k_n\leq p$.

\begin{assumption}[Structural assumptions]\label{A:cor}
~\\[-0.3cm]
\begin{enumerate}[(i).]
 \item  The parameters $\sigma_1,\dots,\sigma_p$ are positive, and there are positive constants $\alpha$, $c$, and $c_{\circ}\in(0,1)$, not depending on $n$, such that
\begin{equation}\label{eqn:varcond1}
\ \sigma_{(j)} \leq c \, j^{-\alpha} \ \ \ \ \text{for all} \ \ \ j\in\{k_n,\dots,p\},
\end{equation}
\vspace{-0.5cm}
\begin{equation}\label{eqn:varcond2}
\ \sigma_{(j)} \geq \ts c_{\circ} \, j^{-\alpha}\  \ \ \text{for all} \ \ \ j\in\{1,\dots,k_n\}.
\end{equation}

\item There is a constant $\e_0\in(0,1)$, not depending on $n$, such that 
\begin{equation}\label{eqn:cormax}
\max_{i\neq j}R_{i,j}(\ell_n)\leq 1-\e_0.
\end{equation}
Also, the matrix $R^+(\ell_n)$ with $(i,j)$ entry given by $\max\{R_{i,j}(\ell_n),0\}$ is positive semi-definite, and there is
 a constant $C>0$ not depending on $n$ such that
\begin{equation}\label{eqn:corbound}
\sum_{1\leq i<j\leq \ell_n}\!\!R_{i,j}^+(\ell_n) \ \leq \ C \ell_n.
\end{equation}
\end{enumerate}
\end{assumption}

\paragraph{Remarks}
Since $\ell_n,k_n\ll n$, it is possible to accurately estimate the parameters $\sigma_{(1)},\dots,\sigma_{(k_n)}$, as well as the matrix $R(\ell_n)$, even when $p$ is large (cf.~Lemmas~\ref{lem:cor} and~\ref{lem:corbasic}). In this sense, it is possible to empirically assess the conditions above.
When considering the size of the decay parameter $\alpha$, note that if $\Sigma$ is viewed as a covariance operator acting on a Hilbert space, then the condition $\alpha>1/2$ essentially corresponds to the case of a trace-class operator --- a property that  is typically assumed in functional data analysis \cp{Hsing2015}.  From this perspective, the condition $\alpha>0$ is very weak, and allows the trace of $\Sigma$ to diverge as $p\to\infty$. 

With regard to the conditions on the correlation matrix $R(\ell_n)$, it is important to keep in mind that they only apply to a small set of variables of size $\mathcal{O}(\log(n)^3)$ --- and the dependence among the variables outside of $\J(\ell_n)$ is \emph{completely unrestricted}.
The interpretation of~\eqref{eqn:corbound} is that it prevents excessive dependence among the coordinates with the largest variances. Meanwhile, the condition that $R^+(\ell_n)$ is positive semi-definite is more technical in nature, and is only used in order to apply a specialized version of Slepian's lemma (Lemma~\ref{lem:slepian}). Nevertheless, this condition always holds in the important case where $R(\ell_n)$ is non-negative. Perturbation arguments may also be used to obtain other examples where some entries of $R(\ell_n)$ are negative.

\subsection{Examples of correlation matrices} 
Some correlation matrices satisfying Assumption~\eqref{A:cor}(ii) are given below.
\small
\begin{align*}
\ \ \ \ \ \ & \text{$\bullet$ \emph{Autoregressive:}} && R_{i,j}=\rho_{0}^{|i-j|} \, , \quad\quad\quad\quad\quad\quad\quad \, \ \text{for any } \rho_0\in(0,1). &&& \\[0.7cm]
\ \ \ \ \ \ &\text{$\bullet$ \emph{Algebraic decay:}} &&  R_{i,j}=1\{i=j\}+\ts\frac{1\{i\neq j\}}{4|i-j|^{\gamma}} \, , \, \ \, \, \ \ \ \text{ for any} \ \gamma \geq 2.  &&&  \\[0.7cm]
\ \ \ \ \ \ &\text{$\bullet$ \emph{Banded:}} && R_{i,j}=\Big(1-\ts\frac{|i-j|}{c_0}\Big)_+ \ ,\quad\quad\quad \,\ \  \text{ for any } c_0>0.\\[0.7cm]
\ \ \ \ \ \ &\text{$\bullet$ \emph{Multinomial:}} && R_{i,j}=1\{i=j\}-\sqrt{\ts\frac{\pi_i\pi_j}{(1-\pi_i)(1-\pi_j)}}1\{i\neq j\}\,\\[0.1cm]
 & \ && \text{where $(\pi_1,\dots,\pi_p)$ is a probability vector.}
\end{align*}
\normalsize
By combining these types of correlation matrices with choices of $(\sigma_1,\dots,\sigma_p)$ that satisfy~\eqref{eqn:varcond1} and~\eqref{eqn:varcond2}, it is straightforward to construct examples of $\Sigma$ that satisfy all aspects of  Assumption~\ref{A:cor}.

\subsection{Examples of variance decay}\label{sec:decayexamples}
To provide additional context for the decay condition~\eqref{eqn:varcond1}, we describe some general situations where it occurs.

\begin{itemize}
\item \emph{Principal component analysis (PCA).} The broad applicability of PCA rests on the fact that many types of data have an underlying covariance matrix with weakly sparse eigenvalues. Roughly speaking, this means that most of the eigenvalues of $\Sigma$ are small in comparison to the top few. Similar to the condition~\eqref{eqn:varcond1}, this situation can be modeled with the decay condition
\begin{equation}\label{eqn:eigdecay}
\lambda_j(\Sigma) \leq c j^{-\gamma},
\end{equation}
for some parameter $\gamma>0$~\cite[e.g.][]{Bunea:Xiao:2015}. Whenever this holds, it can be shown that the variance decay condition \emph{must} hold for some associated parameter $\alpha>0$, and this is done in Proposition~\ref{PROP:DECAY} below. So, in a qualitative sense, this indicates that if a dataset is amenable to PCA, then it is also likely to fall within the scope of our setting.\\

Another way to see the relationship between PCA and variance decay is through the measure of ``effective rank'', defined as
\begin{equation}
{\tt{r}}(\Sigma)=\ts\frac{\tr(\Sigma)}{\|\Sigma\|_{\text{op}}}.
\end{equation}
This quantity has played a key role in a substantial amount of recent work on PCA, because it offers a useful way to describe covariance matrices with an ``intermediate'' degree of complexity, which may be neither very low-dimensional, nor very high-dimensional. We refer to~\cite{Vershynin:2012},~\cite{Lounici:2014},~\cite{Bunea:Xiao:2015}, \cite{Reiss:2016},~\cite{Koltchinskii:Bernoulli:2017,Koltchinskii:2017}, \cite{Nickl:2017}, \cite{Naumov:2017}, and~\cite{Jung:2018}, among others.
Many of these works have focused on regimes where
\begin{equation}\label{eqn:effrank1}
{\tt{r}}(\Sigma)=o(n),
\end{equation}
which conforms naturally with variance decay. Indeed, within a basic setup where $n\asymp p$ and $\|\Sigma\|\op\asymp 1$, the condition~\eqref{eqn:effrank1} holds under $\sigma_{(j)}\leq cj^{-\alpha}$ for any $\alpha>0$.\\[0.2cm]

\item \emph{Count data.} Consider a multinomial model based on $p$ cells and $n$ trials, parameterized by a vector of cell proportions $\boldsymbol\pi=(\pi_1,\dots,\pi_p)$. If the $i$th trial is represented as a vector $X_i\in\R^p$ in the set of standard basis vectors $\{e_1,\dots,e_p\}$, then the marginal distributions of $X_i$ are binomial with $\sigma_j^2=\pi_j(1-\pi_j)$. In particular, it follows that \emph{all} multinomial models satisfy the variance decay condition~\eqref{eqn:varcond1}, because if we let $\boldsymbol \sigma=(\sigma_1,\dots,\sigma_p)$, then the weak-$\ell_2$ norm of $\boldsymbol \sigma$ must satisfy $\|\boldsymbol\sigma\|_{w\ell_2}\leq \|\boldsymbol\sigma\|_2\leq 1$, which implies
\begin{equation}
\sigma_{(j)} \ \leq \ j^{-1/2}
\end{equation}
for all $j\in\{1,\dots,p\}$. In order to study the consequences of this further, we offer some detailed examples in Section~\ref{sec:multinomial}.
More generally, the variance decay condition also arises for other forms of count data. For instance, in the case of a high-dimensional distribution with sparse Poisson marginals, the relation $\var(X_{i,j})=\E[X_{i,j}]$ shows that sparsity in the mean vector can lead to variance decay.\\[-0.2cm]

\item \emph{Fourier coefficients of functional data.} Let $Y_1,\dots,Y_n$ be an i.i.d.~sample of functional data, taking values in a separable Hilbert space $\mathcal{H}$. In addition, suppose that the covariance operator $\mathcal{C}=\text{cov}(Y_1)$ is trace-class, which implies an eigenvalue decay condition of the form~\eqref{eqn:eigdecay}.
Lastly, for each $i\in\{1,\dots,n\}$, let $X_i\in\R^p$ denote the first $p$ generalized Fourier coefficients of $Y_i$ with respect to some fixed orthonormal basis $\{\psi_j\}$ for $\mathcal{H}$. That is,
$X_i=(\langle Y_i,\psi_1\rangle,\dots,\langle Y_i,\psi_p\rangle).$

\ \ \ \ Under the above conditions, it can be shown that no matter which basis $\{\psi_j\}$ is chosen, the vectors $X_1,\dots, X_n$ always satisfy the variance decay condition. (This follows from Proposition~\ref{PROP:DECAY} below.)
 In Section~\ref{sec:expt}, we explore some consequences of this condition as it relates to simultaneous confidence intervals for the Fourier coefficients of the mean function $\E[Y_1]$.
\end{itemize}

To conclude this section, we state a proposition that was used in the examples above. This basic result shows that decay among the eigenvalues $\lambda_1(\Sigma),\dots,\lambda_p(\Sigma)$ requires at least some decay among $\sigma_1,\dots,\sigma_p$. 
\begin{proposition}\label{PROP:DECAY}
Fix two numbers $s\geq 1$, and $r\in(0,s)$. Then, there is a constant $c_{r,s}>0$ depending only on $r$ and $s$, such that for any symmetric matrix $A\in\R^{p\times p}$, we have
$$\|\textup{diag}(A) \|_{w\ell_s}\leq c_{r,s} \| \lambda(A)\|_{w\ell_r}.$$
In particular, if $A=\Sigma$, and if there is a constant $c_0>0$ such that the inequality
$$\lambda_j(\Sigma)\leq c_0\, j^{-1/r}$$
holds for all $1\leq j\leq p$,
then the inequality
$$\sigma_{(j)}^2 \, \leq \, c_0c_{r,s} \, j^{-1/s}$$
holds for all $1\leq j\leq p$.
\end{proposition}
\noindent The proof is given in Appendix~\ref{app:intro}, and follows essentially from the Schur-Horn majorization theorem, as well as inequalities relating $\|\cdot\|_r$ and $\|\cdot\|_{w\ell_r}$.

\section{Main results}\label{sec:main}

In this section, we present our main results on Gaussian approximation and bootstrap approximation.

\subsection{Gaussian approximation}

Let $\tilde S_n\sim N(0,\Sigma)$ and define 
the Gaussian counterpart of the partially standardized statistic $M$  (\ref{eqn:mdef}) according to 
\begin{equation}\label{eqn:Mtildedef}
 \tilde{M}=\max_{1\leq j\leq p} \tilde{S}_{n,j}/\sigma_j^{\tau_n}.
\end{equation}
Our first theorem shows that in the presence of variance decay, the distribution $\mathcal{L}(\tilde M)$ can approximate $\mathcal{L}(M)$ at a nearly parametric rate in Kolmogorov distance. Recall that for any random variables $U$ and $V$, this distance is given by $\dK(\L(U),\L(V))=\sup_{t\in\R}|\P(U\leq t)-\P(V\leq t)|$.

\begin{theorem}[Gaussian approximation]\label{THM:G}
Fix any number $\delta\in(0,1/2)$, and suppose that Assumptions~\ref{A:model} and~\ref{A:cor} hold. In addition, suppose that $\tau_n\in [0,1)$ with $(1-\tau_n)\sqrt{\log(n)}\gtrsim 1$. Then, 
\begin{equation}\label{eqn:thmG}
d_{\textup{K}}\big(\mathcal{L}(M) \, , \, \mathcal{L}(\tilde{M})\big) \ \lesssim  \  n^{-\frac 12+\delta}.
\end{equation}
\end{theorem}
\paragraph{Remarks} As a basic observation, note that the result handles the ordinary max statistic $T$ as a special case with $\tau_n=0$. In addition, it is worth emphasizing that the rate does not depend on the dimension $p$, or the variance decay parameter $\alpha$, provided that it is positive. In this sense, the result shows that even a small amount of structure can have a substantial impact on Gaussian approximation (in relation to existing $n^{-1/6}$ rates that hold when $\alpha=0$). Lastly, the reason for imposing the lower bound on $1-\tau_n$ is that if $\tau_n$  quickly  approaches 1 as $n\to\infty$, then the variances $\var(S_{n,j}/\sigma_j^{\tau_n})$ will also quickly approach 1, thus eliminating the beneficial effect of variance decay.

\subsection{Multiplier bootstrap approximation}\label{sec:MB}

In order to define the multiplier bootstrap counterpart of $\tilde M$, first define the sample covariance matrix 
\begin{equation}\label{eqn:sighat}
\hat\Sigma_n=\frac{1}{n}\sum_{i=1}^n(X_i-\bar X)(X_i-\bar X)\ttop,
\end{equation}
where $\bar X=\frac 1n \sum_{i=1}^n X_i$. Next, let $S_n^{\star}\sim N(0,\hat\Sigma_n)$,
and define the associated max statistic as
\begin{equation}\label{eqn:Mstardef}
M^{\star}=\max_{1\leq j\leq p}S_{n,j}^{\star}/\hat\sigma_j^{\tau_n},
\end{equation}
where $(\hat\sigma_1^2,\dots,\hat\sigma_p^2)=\text{diag}(\hat\Sigma_n)$. 
In the exceptional case when $\hat\sigma_j=0$ for some $j$, the expression $S_{n,j}^{\star}/\hat\sigma_j$ is understood to be 0. This convention is natural, because the event $S_{n,j}^{\star}=0$ holds with probability 1, conditionally on $\hat\sigma_j=0$.

\paragraph{Remarks} 
The above description of $M^{\star}$ differs from some previous works insofar as we have suppressed the role of ``multiplier variables'', and have defined $S_n^{\star}$ as a sample from $N(0,\hat\Sigma_n)$. From a mathematical standpoint, this is equivalent to the multiplier formulation~\citep{CCK:2013}, where $S_n^{\star}= \frac{1}{\sqrt n}\sum_{i=1}^n \xi_i^{\star}(X_i-\bar X)$ and $\xi_1^{\star},\dots,\xi_n^{\star}$ are i.i.d.~$N(0,1)$ random variables, generated independently of $X$.

\begin{theorem}[Bootstrap approximation]\label{THM:BOOT} Fix any number $\delta\in(0,1/2)$, and suppose the conditions of Theorem~\ref{THM:G} hold. Then, there is a constant $c>0$ not depending on $n$, such that the event
\begin{equation}
d_{\textup{K}}\big(\L( \tilde M)\, ,\, \L(M^{\star}|X)\big) \ \leq c\, n^{-\frac 12+\delta}
\end{equation}
occurs with probability at least $1-\frac cn$.
\end{theorem}

\paragraph{Remarks} At a high level, the proofs of Theorems~\ref{THM:G} and~\ref{THM:BOOT} are based on the following observation: When the variance decay condition holds, there is a relatively small subset of $\{1,\dots,p\}$ that is likely to contain the  maximizing index for $M$. In other words, if $\hat\j\in\{1,\dots,p\}$ denotes a random index satisfying $M=S_{n,\hat \j}/\sigma_{\hat \j}^{\tau_n}$, then the ``effective range'' of \, $\hat \j$ \, is fairly small. Although this situation is quite intuitive when the decay parameter $\alpha$ is large, what is more surprising is that the effect persists even for small values of $\alpha$. 

Once the maximizing index $\hat\j$ has been localized to a small set, it becomes possible to use tools that are specialized to the regime where $p\ll n$. For example, Bentkus' multivariate Berry-Esseen theorem~\citep{Bentkus:2003} (cf.~Lemma~\ref{lem:bentkus}) is helpful in this regard. Another technical aspect of the proofs worth mentioning is that they make essential use of the sharp constants in Rosenthal's inequality, as established in~\citep{Zinn:1985} (Lemma~\ref{lem:rosenthal}).

\section{Numerical illustration with functional data}\label{sec:expt}

Due to advances in technology and data collection, functional data have become ubiquitous
in the past two decades, and statistical methods for their analysis have received growing interest. General references and surveys may be found in~\ci{Ramsay2005,Ferraty2006,horv:12,Hsing2015,Wang2016}. 

The purpose of this section is to present an illustration of how the partially standardized statistic $M$ and the bootstrap can be employed to do inference on functional
data.
More specifically, we consider a one-sample test for a mean function, which proceeds by constructing simultaneous confidence intervals (SCI) for its Fourier coefficients. With regard to our theoretical results, this is a natural problem for illustration, because the Fourier coefficients of functional data typically satisfy the variance decay condition~\eqref{eqn:varcondintro}, as explained in the third example of Section~\ref{sec:decayexamples}. Additional background and recent results on mean testing for functional data may be found in \ci{Benko:2009,degras2011simultaneous,Cao:2012,horvath2013,zhen:14,Choi2018,zhan:18}, as well as the references therein.

\subsection{Tests for the mean function}
To set the stage, let $\hilbert$ be a separable Hilbert space of functions, and
let $Y\in \hilbert$ be a random function with mean $\E[Y]=\mu$.
Given a sample $Y_{1},\ldots,Y_{n}$ of i.i.d.~realizations of $Y$,
a basic goal  is to test 
\begin{equation}\label{eqn:test}
H_0: \mu=\mu^{\circ} \text{ \ \ \ \ \ versus \ \ \ \ \ } H_1: \mu\neq \mu^{\circ},
\end{equation}
where $\mu^{\circ}$ is a fixed function in $\hilbert$. 

This testing problem can be naturally formulated in terms of  SCI, as follows.
Let $\{\psi_j\}$ denote any fixed orthonormal basis for $\hilbert$. Also, let $\{u_j\}$ and $\{u_{j}^{\circ}\}$ respectively denote the generalized Fourier coefficients of $\mu$ and $\mu^{\circ}$ with respect to $\{\psi_j\}$, so that
$$\mu=\tsum_{j=1}^{\infty}u_{j}\psi_{j} \text{\  \ \ \ \ and \ \ \ \ \ } \mu^{\circ}=\tsum_{j=1}^{\infty}u_{j}^{\circ}\,\psi_{j}.$$
 Then, the null hypothesis is equivalent to
$u_{j}=u_{j}^{\circ}$ for all $j\geq 1$.
 To test this condition, one can construct a confidence interval $\sci_{j}$ for each $u_{j}$, and reject the null if $u_{j}^{\circ}\not\in \sci_{j}$ for at least one $j\geq 1$. In practice, due to the 
infinite dimensionality of  $\hilbert$, one will choose a sufficiently large integer
$p$, and reject the null if \smash{$u_{j}^{\circ}\not\in \sci_{j}$} for at least one $j\in\{1,\dots,p\}$.

Recently, a similar  general strategy was pursued by  \cite{Choi2018}, hereafter CR, who developed a test for the problem~\eqref{eqn:test} based on a hyper-rectangular confidence region for $(u_1,\dots,u_p)$  --- which is equivalent to constructing SCI.  In the CR approach, the basis is taken to be the eigenfunctions $\{\psi_{\mathcal{C},j}\}$ of the covariance operator \smash{$\mathcal{C}=\text{cov}(Y)$}, and $p$  
is chosen as the number of eigenfunctions $\psi_{\mathcal{C},1},\dots,\psi_{\mathcal{C},p}$ required to explain a certain fraction (say 99\%)  of variance in the data. However, since $\mathcal{C}$ is unknown, the eigenfunctions must be estimated from the available data. 

When $p$ is large, estimating the eigenfunctions $\psi_{\mathcal{C},1},\dots,\psi_{\mathcal{C},p}$ is a well-known challenge in functional data analysis. For instance, a large choice of $p$ may be needed to explain 99\% of the variance if the sample paths of $Y_1,\dots, Y_n$ are not sufficiently smooth. Another example occurs when $H_1$ holds but $\mu$ and $\mu^{\circ}$ are not well separated,
 which may require a large choice of $p$ in order to distinguish $(u_1,\dots,u_p)$ and $(u_1^{\circ},\dots,u_p^{\circ})$. 
In light of these considerations, we will pursue an alternative approach to constructing SCI that does not require estimation of eigenfunctions. 

\subsection{Applying the bootstrap}\label{sec:applyboot} Let $\{\psi_j\}$ be any pre-specified orthonormal basis for $\hilbert$. For instance, when $\hilbert=L^2[0,1]$, a natural option is the standard Fourier basis. For a sample $Y_1,\dots,Y_n\in\hilbert$ as considered before, define random vectors \smash{$X_1,\dots,X_n$} in $\R^p$ according to
$$X_i=(\langle Y_i,\psi_{1}\rangle,\ldots, \langle Y_i,\psi_{p}\rangle),$$
and note that $\E[X_1]=(u_{1},\ldots,u_{p})$.
For simplicity, we retain the previous notations associated with $X_1,\dots,X_n$, so that $S_{n,j}=n^{-1/2}\sum_{i=1}^{n}(X_{i,j}-u_{j})$, and likewise for other quantities. In addition, for any $\tau_n\in[0,1]$, let
\[
L=\min_{1\leq j\leq p}S_{n,j}/\sigma_{j}^{\tau_{n}} \text{ \ \ \ and \ \ \  }M=\max_{1\leq j\leq p}S_{n,j}/\sigma_{j}^{\tau_{n}}.
\]
For a given significance level $\siglevel \in (0,1)$, the $\siglevel$-quantiles of $L$ and $M$ are denoted $q_L(\siglevel)$ and $q_M(\siglevel)$. Thus, the following event occurs with probability at least $1-\siglevel$,
\begin{equation}\label{eq:sci}
\bigcap_{j=1}^{p}\bigg\{ \ts\frac{q_L(\siglevel/2)\sigma_{j}^{\tau_{n}}}{\sqrt{n}}\leq\bar{X}_{j}-u_{j}\leq \ts\frac{q_M(1-\siglevel/2) \sigma_{j}^{\tau_{n}}}{\sqrt{n}}\bigg\}, 
\end{equation}
which leads to theoretical SCI for $(u_1,\dots,u_p)$. 

We now apply the bootstrap from Section~\ref{sec:MB} to estimate $q_L(\siglevel/2)$ and $q_M(1-\siglevel/2)$. Specifically, we generate $B\geq 1$ independent samples of $M^{\star}$ as in~\eqref{eqn:Mstardef}, and then  define $\hat q_M(1-\siglevel/2)$ to be the empirical $(1-\siglevel/2)$-quantile of the $B$ samples (and similarly for $\hat q_L(\siglevel/2)$), leading to the bootstrap SCI 
\begin{equation}\label{eqn:SCIdef}
\hat{\mathcal{I}}_j=\bigg[\bar X_j - \ts\frac{\hat q_M(1-\siglevel/2) \hat \sigma_{j}^{\tau_{n}}}{\sqrt{n}}  \ , \ \bar X_j-\ts\frac{\hat q_L(\siglevel/2) \hat \sigma_{j}^{\tau_{n}}}{\sqrt{n}}\bigg]
\end{equation}
for each $j\in\{1,\dots,p\}$.

It remains to select the value of $\tau_n$, for which we adopt the following simple rule. For each choice of $\tau_n$ in a set of possible candidates, say $\mathcal{T}=\{0,0.1,\ldots,0.9,1\}$, we construct the associated intervals $\hat{\mathcal{I}}_1,\dots,\hat{\mathcal{I}}_p$ as in~\eqref{eqn:SCIdef}, and then select the value $\tau_n\in\mathcal{T}$ for which the average width $\frac 1p \sum_{j=1}^p |\,\hat{\mathcal{I}}_j|$ is the smallest, where  $|[a,b]|=b-a$.

\begin{figure}[H]
\begin{centering}
\includegraphics[width=0.5\textwidth,height=\textheight,keepaspectratio]{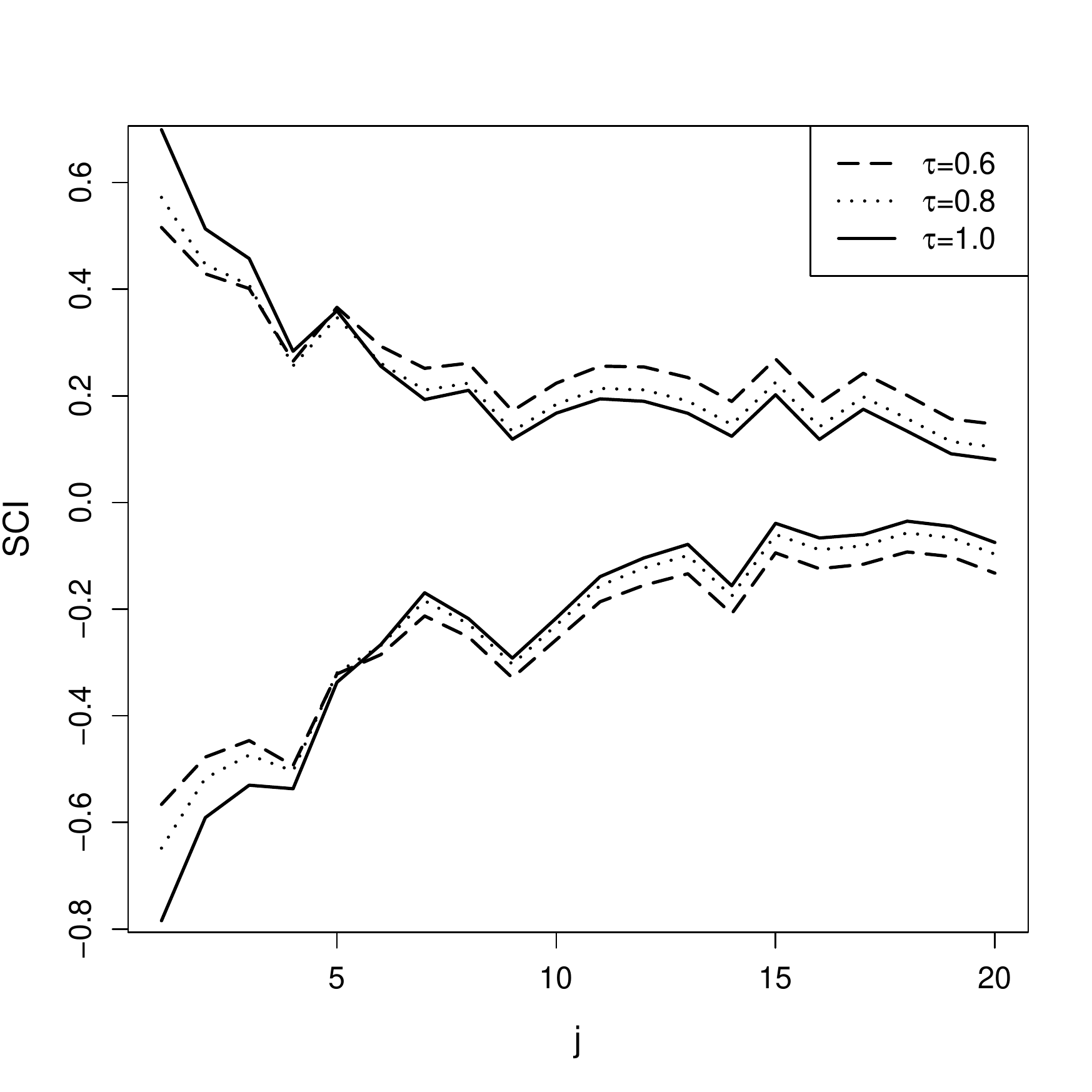}
\par\end{centering}
\caption{Illustration of the impact of $\tau_n$ on the shape of simultaneous confidence intervals (SCI). The curves represent upper and lower endpoints of the respective SCI, where the Fourier coefficients are indexed by $j$. Overall, the plot shows that the SCI change very gradually as a function of $\tau_n$, and that there is a trade-off in the widths of the  intervals. Namely, as $\tau_n$ decreases, the intervals for the leading coefficients (small $j$) become tighter, while the intervals for the subsequent coefficients (large $j$) become wider. }
\label{fig:tau}
\end{figure}

 In Figure \ref{fig:tau}, we illustrate the influence of $\tau_n$ on the shape of the SCI. There are two main points to notice: (1) The intervals change very gradually as a function of $\tau_n$, which shows that partial standardization is at most a mild adjustment of ordinary standardization. (2) The choice of $\tau_n$ involves a tradeoff, which controls the ``allocation of power'' among the $p$ intervals. When $\tau_n$ is close to 1, the intervals are wider for the leading coefficients (small $j$), and narrower for the subsequent coefficients (large $j$). However, as $\tau_n$ decreases from 1, the widths of the intervals gradually become more uniform, and the intervals for the leading coefficients become narrower. Hence, if the vectors $(u_1,\dots,u_p)$ and $(u_1^{\circ},\dots,u_p^{\circ})$ differ in the leading coefficients, then choosing a smaller value of $\tau_n$ may lead to a gain in power. One last interesting point to mention is that in the simulations reported below, the selection rule of ``minimizing the average width'' typically selected values of $\tau_n$ around 0.8, and hence strictly less than 1.

\subsection{Simulation settings}\label{sec:simstudy}
\begin{figure}[]
\begin{centering}
\includegraphics[width=\textwidth,height=\textheight,keepaspectratio]{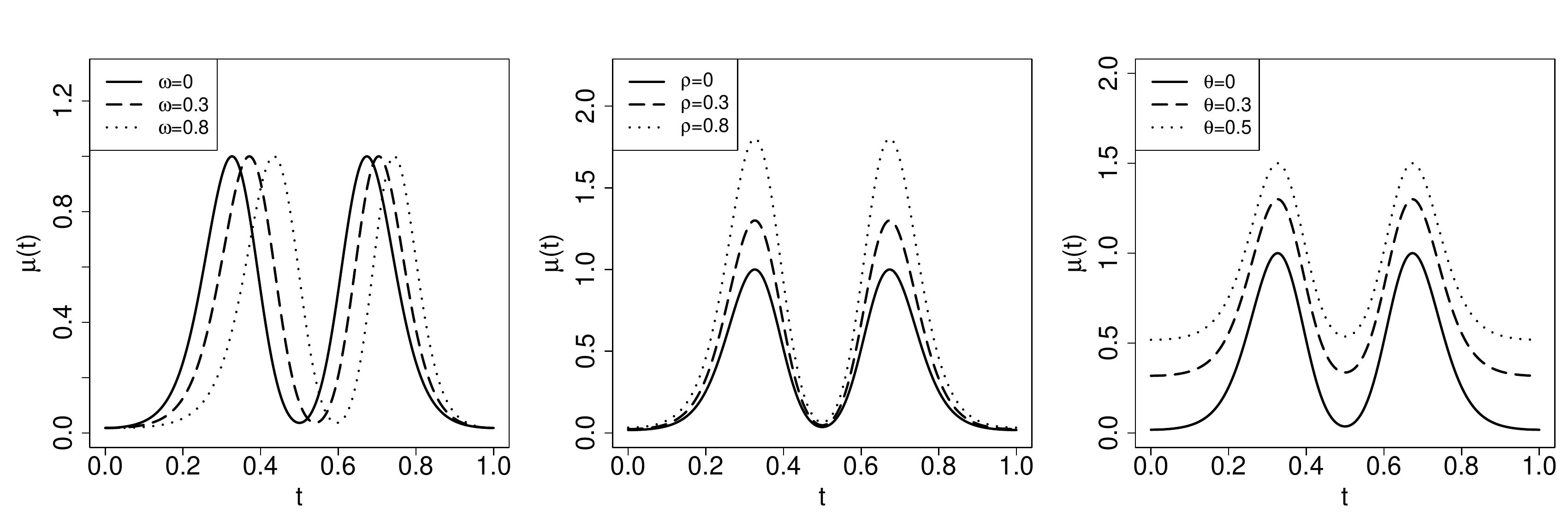}
\par\end{centering}
\caption{Left: Mean functions for varying shape parameters $\omega$ with
$\rho=\theta=0$. Middle: Mean functions for varying scale parameters
$\rho$ with $\omega=\theta=0$. Right: Mean functions with different
shift parameters $\theta$ with $\omega=\rho=0$.}
\label{fig:mean-func-family}
\end{figure}

To study the numerical performance of the SCI described above,
we generated  i.i.d. samples from a Gaussian process
on $[0,1]$, with population mean function 
$$\mu_{\omega,\rho,\theta}(t)=(1+\rho)\cdot\left(\exp[-\{g_{\omega}(t)+2\}^{2}]+\exp[-\{g_{\omega}(t)-2\}^{2}]\right)+\theta$$
indexed by parameters $(\omega,\rho,\theta)$, where $g_{\omega}(t):=8h_{\omega}(t)-4$,
and $h_{\omega}(t)$ denotes the Beta distribution function
with shape parameters $(2+\omega,2)$. This family of functions was
considered in \citet{Chen2012}. To interpret the parameters, note that $\omega$ determines
the shape of the mean function (see Figure \ref{fig:mean-func-family}), whereas
$\rho$ and $\theta$ are scale and shift parameters. In terms of these parameters, the null hypothesis corresponds to $\mu=\mu^{\circ}:=\mu_{0,0,0}$.

\begin{figure}
\begin{centering}
\includegraphics[width=0.45\textwidth,height=\textheight,keepaspectratio]{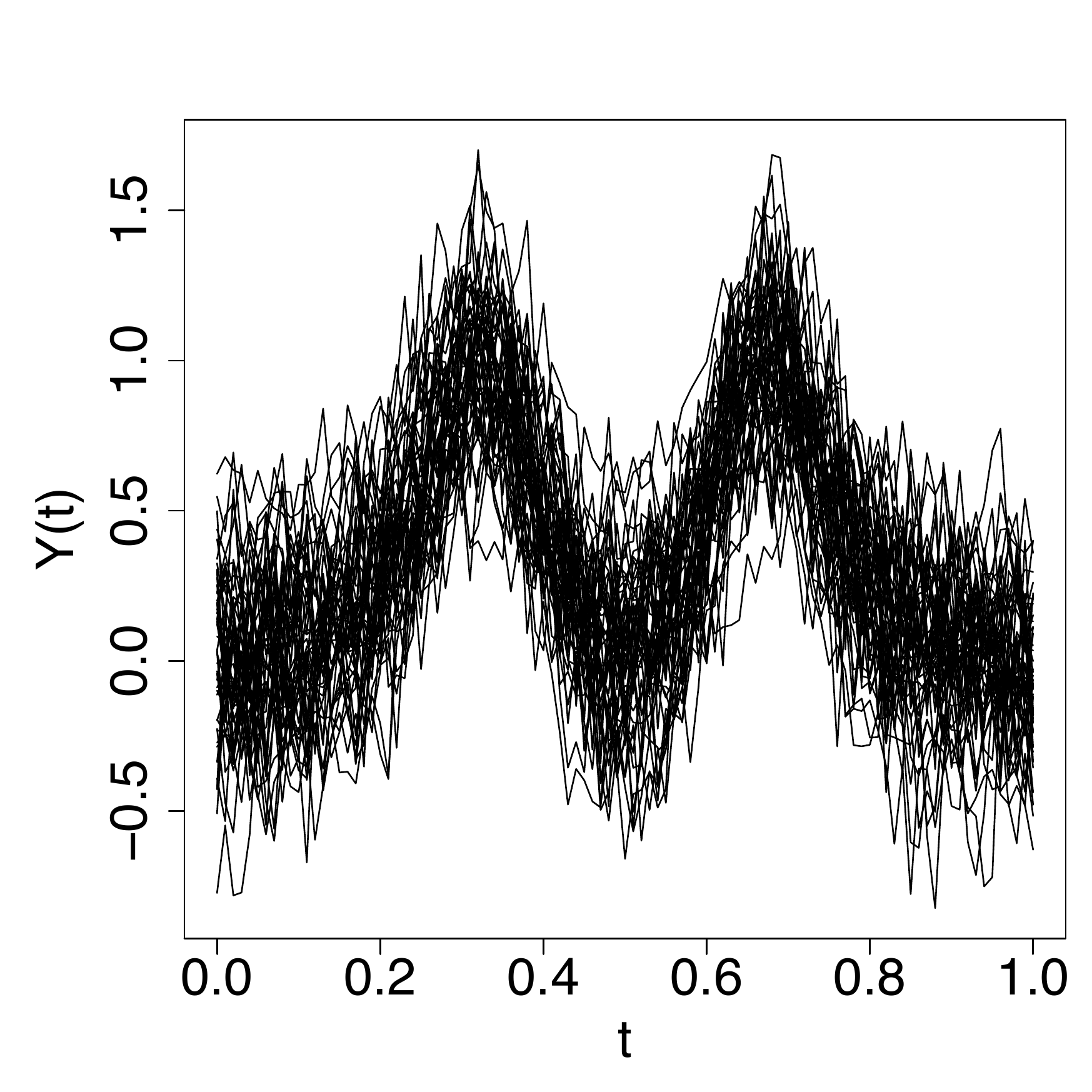} \ \ \ \ \ \
\includegraphics[width=0.45\textwidth,height=\textheight,keepaspectratio]{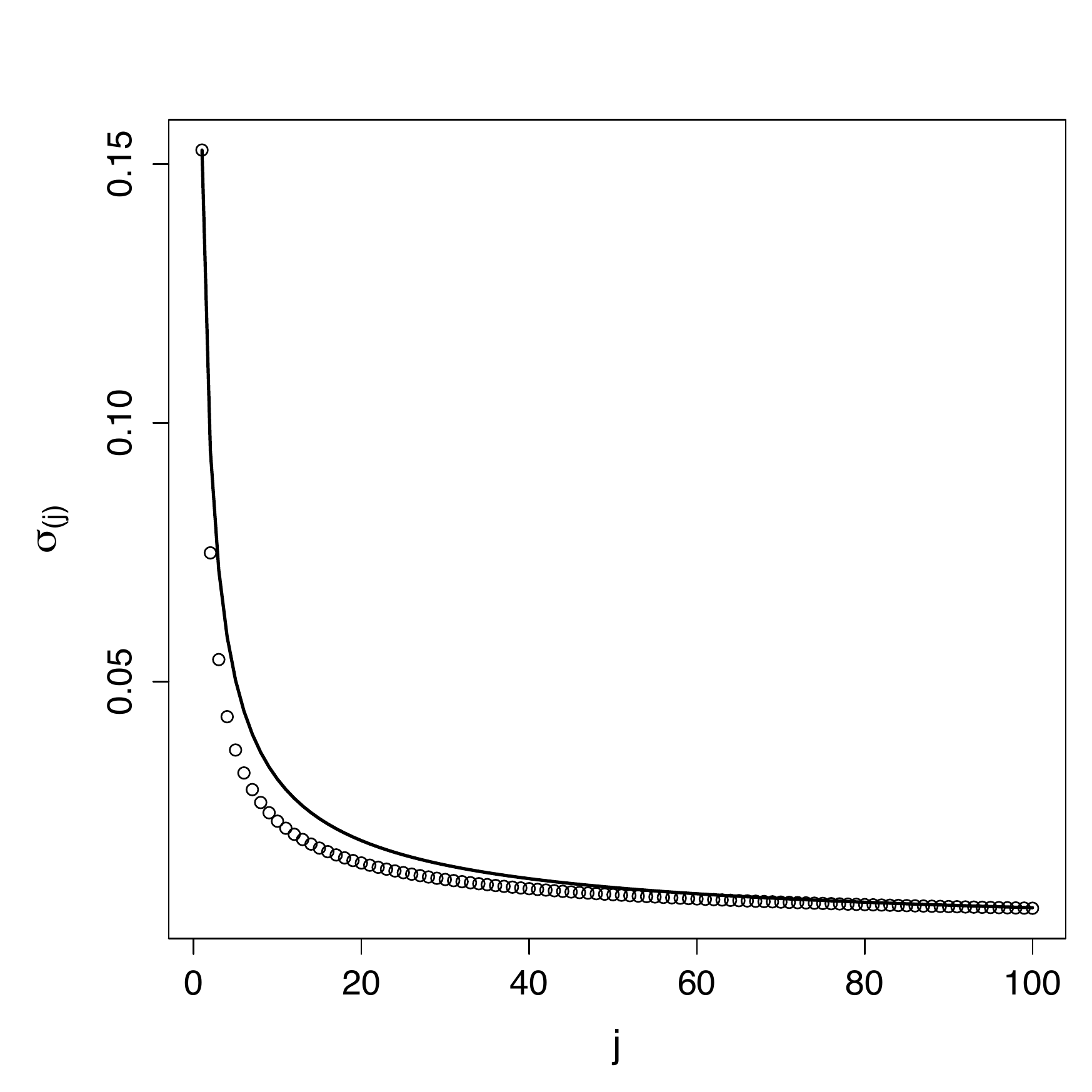}
\par\end{centering}
\caption{Left: A sample of the functional data $Y_1,\dots,Y_n$ in the simulation study. Right: The ordered values $\sigma_{(j)}=\sqrt{\var(X_{1,j})}$ are represented by dots, which are approximated by the decay profile $0.15j^{-0.69}$ (solid line).}
\label{fig:raw-data}
\end{figure}

 The population covariance function was taken to be the Mat\'ern function 
$$\covarop(s,t)=\ts\frac{(\sqrt{2\nu}|t-s|)^{\nu}}{16\Gamma(\nu)2^{\nu-1}}K_{\nu}(\sqrt{2\nu}|t-s|),$$
which was previously considered in CR, with $K_\nu$ being a modified Bessel function of the second kind. 
We set $\nu=0.1$, which results in relatively rough sample paths, as illustrated in the left panel of Figure \ref{fig:raw-data}. Also, the significant presence of variance decay is shown in the right panel.

When implementing the bootstrap in Section~\ref{sec:applyboot}, we  used the first $p=100$ functions from the standard Fourier basis on [0,1]. (In principle, an even larger value $p$ could have been selected, but we chose $p=100$ to limit computation time.) For comparison purposes, we also implemented the `$R_{zs}$' version of the method proposed in CR, using the  accompanying R package~{\tt{fregion}}~\citep{fregion} under default settings, which typically utilized estimates of  the first $p\approx 50$ eigenfunctions of $\mathcal{C}$.

\paragraph{Results on type I error} The nominal significance level was set to $5\%$ in all simulations. To assess the actual type I error, we carried out 5,000 simulations under the null hypothesis, for both  $n=50$ and $n=200$. When $n=50$, the type I error was 6.7\% for the bootstrap method, and 1.6\% for CR. When $n=200$, the results were 5.7\% for the bootstrap method, and  2.6\% for CR. So, in these cases, the bootstrap respects the nominal significance level relatively well. In addition, our numerical results support the idea that partial standardization can be beneficial, because in the fully standardized case where $\tau_n=1$, we observed less accurate type I error rates of 7.0\% for $n=50$, and 6.4\% for $n=200$.

\begin{figure}[H]
\begin{centering}
\includegraphics[width=\textwidth,height=\textheight,keepaspectratio]{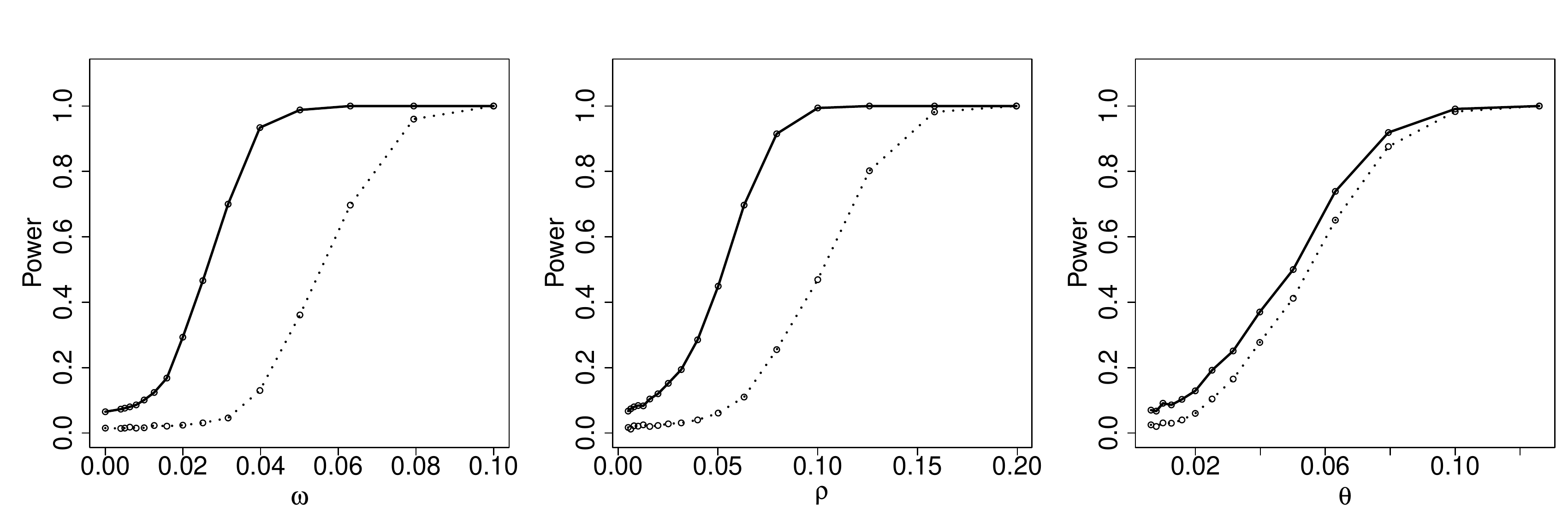}
\par\end{centering}
\caption{Empirical power for the partially standardized bootstrap method (solid) and the CR method (dotted) Left: Empirical power for varying shape parameters $\omega$ while $\rho=\theta=0$. Middle: Empirical power for varying scale parameters
$\rho$ while $\omega=\theta=0$. Right: Empirical power for varying   
shift parameters $\theta$ while $\omega=\rho=0$.}
\label{fig:power}
\end{figure}

\paragraph{Results on power} To  consider power,  we varied each of the parameters $\omega$, $\rho$ and $\theta$, one at a time, while keeping  the other two at their baseline value of zero. In each parameter setting, we carried out 1,000 simulations with sample size $n=50$. The results are summarized in Figure \ref{fig:power}, showing that the bootstrap achieves relative gains in power --- especially with respect to the shape ($\omega$) and scale ($\rho$) parameters. In particular, it seems that using a large number of basis functions can help to catch small differences in these parameters (see also Figure \ref{fig:mean-func-family}).

\section{Examples with multinomial data}\label{sec:multinomial} 
When multinomial models are used in practice, it is not uncommon for the number of cells $p$ to be quite large. Indeed, the challenges of this situation have been a topic of sustained interest, and many inferential questions remain unresolved~\citep[e.g.][]{Hoeffding:1965,Holst:1972,Fienberg:1973,cres:84,zelt:87,Paninski2008,Chafai:2009,Balakrishnan:2019}. A recent survey is~\citep{Balakrishnan:2018}. As one illustration of how our approach can be applied to such models, this section will look at the task of constructing SCI for the cell proportions. Although this type of problem has been studied from a variety of perspectives over the years~\citep[e.g.][]{Quesenberry:1964, Goodman:1965,Fitzpatrick:1987,Glaz:1995,Wang:2008,Chafai:2009}, relatively few theoretical results directly address the high-dimensional setting --- and in this respect, our example offers some progress. Lastly, it is notable that multinomial data are of a markedly different character than the functional data considered in Section~\ref{sec:expt}, which demonstrates how our approach has a broad scope of potential applications.

\subsection{Theoretical example}\label{sec:M1}
Recall from Section~\ref{sec:decayexamples} that we regard the observations in the multinomial model as lying in the set of standard basis vectors $\{e_1,\dots,e_p\}\subset\R^p$. In this context, we also write $\hat\pi_j=\bar{X}_j$ to indicate that the $j$th coordinate of the sample mean is an estimate of the $j$th cell proportion $\pi_j$. In addition, it is important to clarify that a variance decay condition of the form~\eqref{eqn:varcondintro} is automatically satisfied in this model (as explained in Section~\ref{sec:decayexamples}), and so it is not necessary to include this as a separate assumption. Below, we retain the definition of $k_n$ in~\eqref{eqn:kndef}.
\begin{assumption}[Multinomial model]\label{A:mult}
~\\[-0.3cm]
\begin{enumerate}[(i)]
\item The observations $X_1,\dots,X_n\in\R^p$ are i.i.d.,~with $\P(X_1=e_j)=\pi_j$ for each $j\in\{1,\dots,p\}$, where $\boldsymbol\pi=(\pi_1,\dots,\pi_p)$ is a probability vector that may vary with $n$.\\[-0.2cm]
\item There are constants $\alpha>0$ and $\e_0\in(0,1)$, with neither depending on $n$, such that 
\begin{equation}
\sigma_{(j)} \,  \geq \, \e_0 \, j^{-\alpha} \ \ \ \text{ for all } \ \ \ j\in \{1,\dots,(k_n+1)\wedge p\}.
\end{equation}
\end{enumerate}
\end{assumption}
\paragraph{Remarks}
A concrete set of examples satisfying the conditions of Assumption~\ref{A:mult} is given by probability vectors of the form $\pi_{(j)}\propto j^{-\eta}$, with $\eta>1$. Furthermore, the condition $\eta>1$ is mild, since the inequality $\pi_{(j)}\leq j^{-1}$ is satisfied by every probability vector.

\paragraph{Applying the bootstrap} In the high-dimensional setting, the multinomial model differs in an essential way from the model in Section~\ref{sec:prelim}, because there will often be many empty cells (indices) $j\in\{1,\dots,p\}$ for which $\hat\sigma_j=0$. For the indices where this occurs, the usual confidence intervals of the form~\eqref{eqn:SCIdef} have zero width, and thus cannot be used. More generally, if the number of observations in cell $j$ is small, then it is inherently difficult to construct a good confidence interval around $\pi_j$. Consequently, we will restrict our previous SCI~\eqref{eqn:SCIdef} by focusing on a set of cells that contain a sufficient number of observations.
For theoretical purposes, such a set may be defined as
\begin{equation}\label{eq:Jn}
\hat{\mathcal{J}}_n=\Big\{j\in \{1,\dots,p\} \ \Big| \ \hat\pi_j\geq \sqrt{\ts\frac{\log(n)}{n}}\Big\}.
\end{equation}
Accordingly, the max statistic and its bootstrapped version are defined by taking maxima over the indices  in $\hat{\mathcal{J}}_n$, and we denote them as
$$\M = \max_{j\in\hat\J_n} S_{n,j}/\sigma_j^{\tau_n}$$
and 
$$\M^{\star} = \max_{j\in\hat\J_n} S_{n,j}^{\star}/\hat\sigma_j^{\tau_n},$$
where we arbitrarily take $\mathcal{M}$ and $\mathcal{M}^{\star}$ to be zero in the exceptional case when $\hat{\mathcal{J}}$ is empty.

Although the presence of the random index set $\hat\J_n$ complicates the distributions of $\mathcal{M}$ and $\mathcal{M}^{\star}$, it is a virtue of the bootstrap that this source of randomness is automatically accounted for in the resulting inference.
In addition, the following result shows that the bootstrap continues to achieve a near-parametric rate of approximation. 
\begin{theorem}\label{THM:M} Fix any $\delta\in(0,1/2)$, and suppose that Assumption~\ref{A:mult} holds. In addition, suppose that $\tau_n\in[0,1)$ with $(1-\tau_n)\sqrt{\log(n)}\gtrsim 1$. Then, there is a constant $c>0$ not depending on $n$ such that the event
\begin{equation}
 d_{\textup{K}}\big(\L(\M), \L(\M^{\star}|X)\big) \ \leq \ c \, n^{-1/2+\delta},
\end{equation}
occurs with probability at least $1-\frac cn$.
\end{theorem}

\paragraph{Remarks} The proof of this result shares much of the same structure as the proofs of Theorems~\ref{THM:G} and~\ref{THM:BOOT}, but there are a few differences. First, the use of the random index set $\hat\J_n$ in the definition of $\mathcal{M}$ and $\mathcal{M}^{\star}$ entails some extra technical considerations, which are handled with the help of Kiefer's inequality (Lemma~\ref{lem:Kiefer}). Second, we develop a lower bound for $\lambda_{\min}(\Sigma(k_n))$, where $\Sigma(k_n)$ is the covariance matrix of the variables indexed by $\J(k_n)$~(see Lemma~\ref{lem:lambdamin}). This bound may be of independent interest for problems involving multinomial distributions, and does not seem to be well known; see also~\citep{Benasseni:2012} for other related eigenvalue bounds.

\subsection{Numerical example}\label{sec:M2}
We illustrate the bootstrap procedure in the case of the model $\pi_j\propto j^{-1}$, which was considered in a recent numerical study of \cite{Balakrishnan:2018}. 
Taking $p=1000$ and \mbox{$n\in\{500,1000\}$,} we applied the bootstrap method to construct 95\% SCI for the proportions $\pi_j$ corresponding to the cells with at least 5 observations. 
The cutoff value of $5$ is based on a guideline that is commonly recommended in textbooks,~e.g.,~\cite[p.19]{Agresti:2003},~\cite[p.519]{Rice:2007}.
Lastly, the parameter $\tau_n$ was chosen in the same way as described in Section \ref{sec:applyboot}.

Based on  5000 Monte Carlo runs, the average coverage probability was found to be 93.7\% for $n=500$, and 94.4\% for $n=1000$, demonstrating satisfactory performance. Regarding the parameter $\tau_n$, the selection rule typically produced values close to 0.8, for both $n=500$ and $n=1000$. As a point of comparison, it is also interesting to mention the coverage probabilities that occurred when $\tau_n$ was set to 1 (which eliminates all variance decay). In this case, the coverage probabilities became less accurate, with values of $92.7\%$ for $n=500$, and $93.1\%$ for $n=1000$. Hence, this shows that taking advantage of variance decay can enhance coverage probability.

\section{Conclusions}\label{sec:conc}

The main conclusion to draw from our work is that a modest amount of variance decay in a high-dimensional model can substantially improve rates of bootstrap approximation for max statistics --- which helps to reconcile some of the empirical and theoretical results in the literature. In particular, there are three aspects of this type of model structure that are worth emphasizing. First, the variance decay condition~\eqref{eqn:varcondintro} is very weak, in the sense that the parameter $\alpha>0$ is allowed to be arbitrarily small. Second, the condition is approximately checkable in practice, since the parameters $\sigma_1,\dots,\sigma_p$ can be accurately estimated when $n\ll p$. Third, this type of structure arises naturally in a variety of contexts.

Beyond our main theoretical focus on rates of bootstrap approximation, we have also shown that the technique of partial standardization leads to favorable numerical results. Specifically, this was illustrated with examples involving both functional and multinomial data, where variance decay is an inherent property that can be leveraged. 
Finally, we note that these applications are by no means exhaustive, and the adaptation of the proposed approach to other types of data may provide further opportunities for future work.

\references

\newpage
\begin{center}
\Large{Supplementary Material}
\end{center}

\paragraph{Organization of appendices} 
 In Appendix~\ref{app:intro} we prove Proposition~\ref{PROP:DECAY}, and in Appendices~\ref{app:thmg} and~\ref{app:thmboot} we prove Theorems~\ref{THM:G} and~\ref{THM:BOOT} respectively. These proofs rely on numerous technical lemmas, which are stated and proved in Appendix~\ref{app:technical}. Next, the proof of Theorem~\ref{THM:M} for the multinomial model is given in Appendix~\ref{app:mult}, and the associated technical lemmas are given in Appendix~\ref{app:lemmamult}. Lastly, in Appendix~\ref{app:background} we provide statements of background results, and in Appendix~\ref{app:Gaussian} we provide a discussion of related work in the Gaussian approximation literature.

\paragraph{General remarks and notation} Based on the formulation of Theorems~\ref{THM:G},~\ref{THM:BOOT}, and~\ref{THM:M}, it is sufficient to show that these results hold for all large values of $n$, and it will simplify some of the proofs to make use of this reduction. For this reason, it is  understood going forward that $n$ is sufficiently large for any given expression to make sense. Another convention is that all proofs in 
Appendices~\ref{app:thmg},~\ref{app:thmboot},~\ref{app:technical},~\ref{app:mult}, and~\ref{app:lemmamult} will implicitly assume that $p>k_n$ (unless otherwise stated), because once the proofs are given for this case, it will follow that the low-dimensional case where $p\leq k_n$ can be handled as a direct consequence (which is explained on page~\pageref{lowdimcase}).

To fix some notation that will be used throughout the appendices, let $d\in\{1,\dots,p\}$, and define a generalized version of $M$ as
$$M_d=\max_{j\in\J(d)}S_{n,j}/\sigma_j^{\tau_n}.$$
In particular, the statistic $M$ defined in equation~\eqref{eqn:mdef} is the same as $M_p$.
Similarly, the Gaussian and bootstrap versions of $M_d$ are defined as
$$\tilde M_d=\max_{j\in\J(d)}\tilde S_{n,j}/\sigma_j^{\tau_n},$$
and
$$M_d^{\star} =\max_{j\in\J(d)}S_{n,j}^{\star}/\hat\sigma_j^{\tau_n}.$$
 In addition, define the parameter
\begin{equation}
\beta_n=\alpha(1-\tau_n).
\end{equation}
Lastly, we will often use the fact that if a random variable $\xi$ satisfies the bound $\|\xi\|_{\psi_1}\leq c$ for some constant $c$ not depending on $n$, then there is another constant $C>0$ not depending on $n$, such that $\|\xi\|_r\leq C\, r$ for all $r\geq 1$~\cite[Proposition 2.7.1]{VershyninHDP}.

\appendix

\section{Proof of Proposition~\ref{PROP:DECAY}} \label{app:intro}
\proof 
It is a standard fact that for any $s\geq 1$, the $\ell_s$ norm dominates its $w\ell_s$ counterpart, and so
$\|\text{diag}(A)\|_{w\ell_s}\leq \|\text{diag}(A)\|_s$.
Next, since $A$ is symmetric, the Schur-Horn Theorem implies that the vector $\text{diag}(A)$ is majorized by $\lambda(A)$~\cite[p.300]{Olkin:2011}. Furthermore, when $s\geq 1$, the function $\|\cdot \|_s$ is Schur-convex on $\R^p$, which means that if  $u\in\R^p$ is majorized by $v\in\R^p$, then $\|u\|_s\leq \|v\|_s$~\cite[p.138]{Olkin:2011}. Hence,
$$\|\text{diag}(A)\|_{w\ell_s}\leq \| \lambda(A)\|_s.$$
Finally, if $r\in(0,s)$, then for any $v\in\R^p$, the inequality 
$$\|v\|_s \leq \big(\zeta(s/r)\big)^{1/s}\,\|v\|_{w\ell_r}$$ 
holds, where $\zeta(x):=\sum_{j=1}^{\infty}j^{-x}$ for $x>1$. This bound may be derived as in~\citep[p.257]{Johnstone:2017}, 
$$\|v\|_s^s \ = \ \sum_{j=1}^p |v|_{(j)}^s \ \leq \ \sum_{j=1}^p \big(\|v\|_{w\ell_r} j^{-1/r}\big)^s \ \leq \ \zeta(s/r)\cdot \|v\|_{w\ell_r}^s,
$$
which completes the proof.
\qed

\section{Proof of Theorem~\ref{THM:G}}\label{app:thmg}

\proof Consider the inequality
\begin{equation}
d_{\textup{K}}(\mathcal{L}(M_p),\mathcal{L}(\tilde{M}_p)) \, \leq \, \I_n+\II_n+\III_n,
\end{equation}
where we define
\begin{align}
\I_n &= d_{\text{K}}\Big(\mathcal{L}(M_p) \, , \, \mathcal{L}(M_{k_n})\Big)\\[0.2cm]
\II_n &=d_{\text{K}}\Big( \mathcal{L}(M_{k_n})\, , \, \mathcal{L}(\tilde M_{k_n})\Big)\\[0.2cm]
\III_n &=d_{\text{K}}\Big( \mathcal{L}(\tilde M_{k_n})\, ,\, \mathcal{L}(\tilde M_p)\Big).\label{eqn:IIIdef}
\end{align}
Below, we show that the term $\II_n$ is at most of order $n^{-\frac{1}{2}+\delta}$ in Proposition~\ref{prop:bentkus}. Later on, we establish a corresponding result for $\I_n$ and $\III_n$ in Proposition~\ref{prop:IandIII}. Taken together, these results complete the proof of Theorem~\ref{THM:G}.\qed

\begin{proposition}\label{prop:bentkus}
Fix any number $\delta\in(0,1/2)$, and suppose the conditions of Theorem~\ref{THM:G} hold. Then, 
\begin{equation}\label{eqn:IInbound}
\II_n \, \lesssim \,  n^{-\frac 12+\delta}.
\end{equation}
\end{proposition}
\proof 
Let $\Pi_{k_n}\in\R^{k_n\times p}$ denote the projection onto the coordinates indexed by $\J(k_n)$. This means that if we write $\J(k_n)=\{j_1,\dots,j_{k_n}\}$ so that $(\sigma_{j_1},\dots,\sigma_{j_{k_n}})=(\sigma_{(1)},\dots,\sigma_{(k_n)})$, then the $l$th row of $\Pi_{k_n}$ is the standard basis vector $e_{j_l}\in\R^p$.
Next, define the diagonal matrix $D_{k_n}=\diag(\sigma_{(1)},\dots,\sigma_{(k_n)})$. It follows that
$$M_{k_n}=\max_{1\leq j\leq k_n} e_j\ttop D_{k_n}^{-\tau_n}\Pi_{k_n} S_n.$$
Define the matrix $\mathfrak{C}\ttop =D_{k_n}^{-\tau_n} \Pi_{k_n}\Sigma^{1/2}$, which is of size $k_n\times p$. Also, let $r$ denote the rank of $\mathfrak{C}$, and note that $r\leq k_n$, since the matrix $\Sigma$ need not be invertible. Next, consider a decomposition
$$\mathfrak{C}=QR,$$
where the columns of $Q\in\R^{p\times r}$ are an orthonormal basis for the image of $\mathfrak{C}$, and $R\in\R^{r\times k_n}$. Hence, if we define the random vector
\begin{equation}\label{eqn:Zrep}
\breve Z=\ts\frac{1}{\sqrt n}\sum_{i=1}^n Q\ttop Z_i,
\end{equation}
then we have
$$D_{k_n}^{-\tau_n}\Pi_{k_n} S_n =R\ttop \breve Z.$$
It is simple to check that for any fixed $t\in\R$, there exists a Borel convex set $\mathcal{A}_t\subset\R^{r}$ such that $\P(M_{k_n}\leq t)=\P(\breve Z\in\mathcal{A}_t)$. By the same reasoning, we also have $\P(\tilde M_{k_n}\leq t)=\gamma_r(\mathcal{A}_t)$, where $\gamma_{r}$ is the standard Gaussian distribution on $\R^{r}$. Therefore, the quantity $\II_n$  satisfies the bound
 \begin{equation}
 \begin{split}
 \II_n
 & \ \leq  \ \sup_{\mathcal{A}\in\mathscr{A}}\, \Big| \P\big(\breve Z\in \mathcal{A} \big)- \gamma_{r}(\mathcal{A})\Big|,
 \end{split}
\end{equation}
where $\mathscr{A}$ denotes the collection of all Borel convex subsets of $\R^{r}$.

We now apply Theorem 1.1 of~\cite{Bentkus:2003} (Lemma~\ref{lem:bentkus}), to handle the supremum above. First observe that
the definition of $\breve Z$ in~\eqref{eqn:Zrep} satisfies the conditions of that result, since the terms $Q\ttop Z_1,\dots, Q\ttop Z_n$ are i.i.d.~with zero mean and identity covariance matrix. Therefore,
$$\II_n \ \lesssim \ r^{1/4}\cdot  \E\big[ \|Q\ttop Z_1\|_2^3] \cdot n^{-1/2}.$$

It remains to bound the middle factor on the right side. By Lyapunov's inequality, 
\begin{equation}\label{eqn:fourthcalc}
\begin{split}
\E\big[\|Q\ttop Z_1\|_2^3\big] &\leq \Big(\E\Big[\big(Z_1\ttop QQ\ttop Z_1\big)^2\Big]\Big)^{3/4}.
\end{split}
\end{equation} 
Next, if $v_1,\dots,v_r$ denote the orthonormal columns of $Q$, then we have
$$QQ\ttop = \tsum_{j=1}^r v_jv_j\ttop.$$
Hence, if we put $\zeta_j=Z_1\ttop v_j$, then
\begin{equation}
\begin{split}
 \E\Big[(Z_1\ttop QQ\ttop Z_1)^2\Big] \ &= \ \Big\|\tsum_{j=1}^r \zeta_j^2\Big\|_2^2\\[0.2cm]
 & \ \leq \ \Big(\tsum_{j=1}^r \|\zeta_j^2\|_2\Big)^2\\[0.2cm]
 & \ \lesssim \ k_n^2,
 \end{split}
\end{equation}
where we have used the fact that $r\leq k_n$ and $\|\zeta_j^2\|_2=\|Z_1\ttop v_j\|_4^2\lesssim 1$, based on Assumption~\ref{A:model}.
Combining the last few steps gives $\E[\|Q\ttop Z_1\|_2^3]\lesssim k_n^{6/4}$, and hence
$$ \II_n  \ \lesssim \ k_n^{7/4}n^{-1/2}\ \lesssim \ n^{-\frac 12+\delta},$$
as needed.
\qed 

~\\

\begin{proposition}\label{prop:IandIII}
Fix any number $\delta\in(0,1/2)$, and suppose the conditions of Theorem~\ref{THM:G} hold. Then, 
\begin{equation}\label{IandIII}
\I_n \, \lesssim \, n^{-\frac 12+\delta} \ \ \ \text{ and } \ \ \ \III_n \, \lesssim \, n^{-\frac 12+\delta}.
\end{equation}
\end{proposition}
\proof We only prove the bound for $\I_n$, since the same argument applies to $\III_n$.
It is simple to check that for any fixed real number $t$,
$$\Big|\P\Big(\max_{1\leq j\leq p} S_{n,j}/\sigma_j^{\tau_n} \leq t\Big)-\P\Big(\max_{j\in\J(k_n)} S_{n,j}/\sigma_j^{\tau_n} \leq t\Big)\Big| = \ \P\Big(A(t)\cap B(t)\Big),$$
where we define the events
\begin{equation}\label{inclusion}
\small
A(t)=\Big\{\max_{j\in\J(k_n)} S_{n,j}/\sigma_j^{\tau_n} \leq t\Big\} \ \ \ \text{ and } \ \  \ \ B(t)=\Big\{\max_{j\in\J(k_n)^c} S_{n,j}/\sigma_j^{\tau_n}> t\Big\},
\end{equation}
and $\J(k_n)^c$ denotes the complement of $\J(k_n)$ in $\{1,\dots,p\}$.
Also, for any pair of real numbers $t_{1,n}$ and $t_{2,n}$ satisfying $t_{1,n}\leq t_{2,n}$, it is straightforward to check that the following inclusion holds for all $t\in\R$,
\begin{equation}\label{inclusion}
A(t)\cap B(t) \ \subset \ A(t_{2,n})\cup B(t_{1,n}).
\end{equation}
Applying a union bound, and then taking the supremum over $t\in\R$, we obtain
$$\I_n\,\leq \, \P(A(t_{2,n}))\,+\, \P(B(t_{1,n})).$$

The remainder of the proof consists in selecting $t_{1,n}$ and $t_{2,n}$ so that $t_{1,n}\leq t_{2,n}$ and that the probabilities $\P(A(t_{2,n}))$ and $\P(B(t_{1,n}))$ are sufficiently small. Below, Lemma~\ref{lem:bounds} shows that if $t_{1,n}$ and $t_{2,n}$ are chosen as 
\begin{align}
t_{1,n}&=c\cdot k_n^{-\beta_n}\cdot \log(n)  \label{eqn:t1}\\[0.3cm]
t_{2,n}&=c_{\circ}\cdot \ell_n^{-\beta_n}\cdot\sqrt{\log(\ell_n)},\label{eqn:t2}
\end{align}
for a certain constant $c>0$, and $c_{\circ}$ as in~\eqref{eqn:varcond2}, then  $\P(A(t_{2,n}))$ and $\P(B(t_{1,n}))$ are at most of order $n^{-\frac 12+\delta}$. Furthermore, the inequality $t_{1,n}\leq t_{2,n}$ holds for all large $n$, due to the definitions of $\ell_n$, $k_n$, and $\beta_n$, as well as the condition $(1-\tau_n)\sqrt{\log(n)}\gtrsim 1$.
\qed
\begin{lemma}\label{lem:bounds}
Fix any number $\delta\in(0,1/2)$, and suppose the conditions of Theorem~\ref{THM:G} hold. Then, there are positive constants $c$ and $c_{\circ}$, not depending on $n$, that can be selected in the definitions of $t_{1,n}$~\eqref{eqn:t1} and $t_{2,n}$\,\eqref{eqn:t2}, so that
\begin{equation}\label{aboundlem}
\P(A(t_{2,n})) \, \lesssim \, n^{-\frac 12+\delta},\tag{a}
\end{equation}
and
\begin{equation}\label{bboundlem}
\P(B(t_{1,n})) \, \lesssim \,  n^{-1}.\tag{b}
\end{equation}
\end{lemma}

\paragraph{\textsc{Proof of Lemma~\ref{lem:bounds} part} \eqref{aboundlem}} Due to Proposition~\ref{prop:bentkus} and the fact that $\J(\ell_n)\subset \J(k_n)$, we have
\begin{equation}
\begin{split}
\P(A(t_{2,n})) &\ \leq \ \P\Big(\max_{j\in\J(k_n)} \tilde S_{n,j}/\sigma_j^{\tau_n}\leq t_{2,n}\Big)+\II_n\\[0.2cm]
& \ \leq \  \ \P\Big(\max_{j\in\J(\ell_n)} \tilde S_{n,j}/\sigma_j^{\tau_n}\leq t_{2,n}\Big)+  c\, n^{-\frac 12+\delta}.
\end{split}
\end{equation}
To bound the probability in the last line, we will make use of the Gaussianity of $\tilde S_n$ to apply certain results based on Slepian's lemma, as contained in Lemmas~\ref{lem:slepian2} and \ref{lem:slepian} below.  As a preparatory step, consider some generic random variables $\{Y_j\}$ and positive scalars $\{a_j\}$ indexed by $\J(\ell_n)$, as well as a constant $b$ such that $\max_{j\in\J(\ell_n)}a_j\leq b$. Then,
\begin{equation}
\begin{split}
\P\Big(\max_{j\in\J(\ell_n)} Y_j\leq t_{2,n}\Big) & \leq \ \P\Big(\max_{j\in\J(\ell_n)}a_jY_j \leq b\, t_{2,n}\Big),
\end{split}
\end{equation}
which can be seen by expressing the left side in terms of $\cap_j\{a_jY_j\leq a_j t_{2,n}\}$, and noting that this set is contained in $\cap_j\{a_jY_j\leq b t_{2,n}\}$. Due to Assumption~\ref{A:cor} with $c_{\circ}\in(0,1)$, we have the inequality $\sigma_j^{\tau_n-1}\leq \ell_n^{\beta_n}/c_{\circ}$ for all $j\in\J(\ell_n)$, and so we may apply the previous observation with  $a_j=\sigma_j^{\tau_n-1}$, and $b=\ell_n^{\beta_n}/c_{\circ}$. Furthermore, the definition of $t_{2,n}$ gives $b\, t_{2,n}=\sqrt{\log(\ell_n)}$, and so we if we let $Y_j=\tilde S_{n,j}/\sigma_j^{\tau_n}$, it follows that
$$\P\Big(\max_{j\in\J(\ell_n)} \tilde S_{n,j}/\sigma_j^{\tau_n}\leq t_{2,n}\Big) \ \leq \  \P\Big(\max_{j\in\J(\ell_n)} \tilde S_{n,j}/\sigma_j\leq \sqrt{\log(\ell_n)}\Big).$$
The proof is completed by applying the next result (Lemma~\ref{lem:slepian2}) in conjunction with the conditions of Assumption~\ref{A:cor}. (Take $m=\ell_n$ in the statement of Lemma~\ref{lem:slepian2}.)\qed

\paragraph{Remark} The lemma below may be of independent interest, and so we have stated it in a way that can be understood independently of the context of our main assumptions. Also, the constants $1/2$ and $1/3$ in the exponent of the bound~\eqref{eqn:slepian2result} can be improved slightly, but we have left the result in this form for simplicity.

\begin{lemma}\label{lem:slepian2} 
For each integer $m\geq 1$, let $\mathsf{R}=\mathsf{R}(m)$ be a correlation matrix in $\R^{m\times m}$, and let $\mathsf{R}^+=\mathsf{R}^+(m)$ denote the matrix with $(i,j)$ entry given by $\max\{\mathsf{R}_{i,j},0\}$. Suppose the matrix $\mathsf{R}^+$ is positive semi-definite for all $m$, and that there are constants $\e_0\in(0,1)$ and $c>0$, not depending on $m$, such that the inequalities
\begin{align}
\ \max_{i\neq j}\mathsf{R}_{i,j}  & \  \leq \, 1-\e_0\\[0.2cm]
\ts\sum_{i\neq j} \mathsf{R}_{i,j}^+ & \  \leq  \ c\, m\label{eqn:rpluscondn}
\end{align}
hold for all $m$.
Lastly, let $(\zeta_1,\dots,\zeta_m)$ be a Gaussian vector drawn from $N(0,\mathsf{R})$. Then, there is a constant $C>0$, not depending on $m$, such that the inequality
\begin{equation}\label{eqn:slepian2result}
\P\Big(\max_{1\leq j\leq m} \zeta_j\leq \sqrt{\log(m)}\Big) \ \leq \ C\exp\big(-\ts\frac{1}{2}m^{1/3}\big)
\end{equation}
holds for all $m$.
\end{lemma}
\proof It is enough to show that the result holds for all large $m$, because if $m\leq m_0$, then the result is clearly true when $C= \exp(\frac{1}{2}m_0^{1/3})$. To begin the argument,  we may introduce a Gaussian vector $(\xi_1,\dots,\xi_m)\sim N(0,\mathsf{R}^+)$, since the matrix $\mathsf{R}^+$ is positive semi-definite. In turn, the version of Slepian's Lemma given in Lemma~\ref{lem:slepian} leads to
\begin{equation}\label{eqn:splitfactor}
\begin{split}
\P\Big(\max_{1\leq j\leq m}\zeta_j \leq \sqrt{\log(m)}\Big) 
& \ \leq \ \P\Big(\max_{1\leq j\leq \ell_n}\xi_j \leq  \sqrt{\log(m)}\Big)\\[0.3cm]
 & \ \leq \ K_m\cdot \Phi^{m}\big( \sqrt{\log(m)}\big),
\end{split}
\end{equation}
where we put
$$K_m=  \exp\Bigg\{\sum_{1\leq i<j\leq m} \log\Big(\ts\frac{1}{1-\frac{2}{\pi}\arcsin(\mathsf{R}_{i,j}^+)}\Big)\exp\Big(-\ts\frac{\log(m)}{1+\mathsf{R}_{i,j}^+}\Big)\Bigg\}.$$
Next, we apply the assumption  $\max_{i\neq j}\mathsf{R}_{i,j}^+\leq 1-\e_0$. Since the functions $x\mapsto \log(1/(1-x))$ and  $x\mapsto\frac{2}{\pi}\arcsin(x)$ have bounded derivatives on any closed subinterval of $[0,1)$,  it follows that
 $$\log\Big(\ts\frac{1}{1-\frac{2}{\pi}\arcsin(\mathsf{R}_{i,j}^+)}\Big) \ \leq c\, \mathsf{R}_{i,j}^+,$$
 for some constant $c>0$ not depending on $m$. Therefore, by possibly increasing $c$, the condition~\eqref{eqn:rpluscondn} gives
\begin{equation}
\begin{split}
K_m
& \ \leq \ \exp\Big\{ cm \cdot \exp\Big(-\ts\frac{\log(m)}{1+(1-\e_0)}\Big)\Big\}\\[0.3cm]
& \ = \ \exp\Big\{cm^{1-\frac{1}{2-\e_0}}\Big\}.
\end{split}
\end{equation}
To bound the earlier factor involving $\Phi^m(\sqrt{\log(m)})$, let $\eta_0\in(0,1)$ be a small constant to be optimized below, and note that the following inequality holds for all sufficiently large $s>0$,
\begin{equation}\label{numericalineq}
\begin{split}
\Phi\Big(\sqrt{(2-\eta_0\e_0)\log(s)}\Big) \ \leq \ 1-\ts\frac{1}{s},
\end{split}
 \end{equation}
which may be found in~\cite[p.337]{Massart:2013}. Taking $s=m^{\kappa_0}$ with $\kappa_0:=\frac{1}{2-\eta_0\e_0}$ shows that for all large $m$,
\begin{equation}
\begin{split}
\Phi^{m}\big(\sqrt{\log(m)}\big) & \ \leq  \ \Big(1-\ts\frac{1}{m^{\kappa_0}}\Big)^{m}\\[0.3cm]
& \ \leq \ \exp\big(-m^{1-\kappa_0}\big).
\end{split}
\end{equation}
We now collect the last several steps. If we observe that $\kappa_0<\frac{1}{2-\e_0}$, then the following inequalities hold for all large $m$,
\begin{equation}
\begin{split}
  K_m\cdot \Phi^{m}\big( \sqrt{\log(m)}\big) & \ \leq \ \exp\Big\{cm^{1-\frac{1}{2-\e_0}}-m^{1-\kappa_0}\Big\}\\[0.3cm]
  & \ \leq \ \exp\Big\{-(1-\eta_0)m^{1-\kappa_0}\Big\}.
  \end{split}
\end{equation}
So, by possibly further decreasing $\eta_0$, we have $(1-\kappa_0)>1/3$, as well as $(1-\eta_0)>1/2$. This leads to the stated result.
\qed

~\\

\paragraph{\textsc{Proof of Lemma~\ref{lem:bounds} part}~\eqref{bboundlem}} Define the random variable
$$V=\max_{j\in\J(k_n)^c} S_{n,j}/\sigma_j^{\tau_n},$$
and let
$$q=\max\big\{\ts\frac{2}{\beta_n}, \log(n), 3\big\}.$$
Clearly, for any $t>0$, we have the tail bound
\begin{equation}\label{qchebyshev}
\P\big( V\geq t)\leq \frac{\|V\|_q^q}{t^q},
\end{equation}
and furthermore
 \begin{equation}
 \begin{split}
  \|V\|_q^q &= \E\bigg[\Big|\max_{j\in\J(k_n)^c} S_{n,j}/\sigma_j^{\tau_n}\Big|^q\bigg]\\[0.3cm]
  &\leq \sum_{j\in\J(k_n)^c} \sigma_j^{q(1-\tau_n)}\,\E\big[|\ts\frac{1}{\sigma_j}S_{n,j}|^q\big].
  \end{split}
 \end{equation}
By Lemma~\ref{lem:Snjnorm}, we have $\|\frac{1}{\sigma_j}S_{n,j}\|_{q}\leq cq$, and so 
  \begin{equation}
  \begin{split}
    \|V\|_q^q \ &\leq \  (c q)^q \sum_{j\in\J(k_n)^c} \sigma_j^{q(1-\tau_n)} \\[0.3cm]
  &\lesssim \ (c q)^q \sum_{j=k_n+1}^p j^{-q \beta_n }\\[0.3cm]
  &\leq \ (cq)^q \int_{k_n}^p x^{-q\beta_n} dx\\[0.3cm]
  &\leq \ \ts\frac{(c q)^q }{q \beta_n-1}\ k_n^{-q \beta_n  +1}, 
 \end{split}
 \end{equation}
 where we recall $\beta_n=\alpha(1-\tau_n)$, and note that $ q\beta_n\geq 2$, which holds by the definition of $q$.
Hence, if we put $C_n:=\ts\frac{c}{(q\beta_n -1)^{1/q}} \cdot k_n^{1/q}$, then
$$\|V\|_q \leq \ C_n \cdot q\cdot k_n^{-\beta_n}.$$
Furthermore, it is simple to check that $C_n\lesssim 1$, and that the assumption $(1-\tau_n)\sqrt{\log(n)}\gtrsim 1$ implies
$q\lesssim \log(n).$
Therefore, from the inequality~\eqref{qchebyshev} with $t=e\|V\|_q$, as well as the definition of $q$, we obtain
$$\P\Bigg(V\geq c\cdot \log(n)\cdot k_n^{-\beta_n}\Bigg) \ \leq \ e^{-q} \ \leq \ \ts\frac 1n,$$
for some constant $c>0$ not depending on $n$, as needed.\qed

~\\

\section{Proof of Theorem~\ref{THM:BOOT}}\label{app:thmboot} 
\proof Consider the inequality
\begin{equation}
d_{\textup{K}}(\mathcal{L}(\tilde M_p),\mathcal{L}(M_p^{\star}| X)) \ \leq \ \I'_n \ + \ \II'_n(X) \ + \ \III'_n(X),
\end{equation}
where we define
\begin{align}
\I'_n & \ = \  d_{\text{K}}\Big(\mathcal{L}(\tilde M_p) \, , \, \mathcal{L}(\tilde M_{k_n})\Big)\\[0.2cm]
\II'_n(X) &\ = \ d_{\text{K}}\Big(\mathcal{L}(\tilde M_{k_n}) \, , \, \mathcal{L}( M^{\star}_{k_n}|X\big)\Big)\\[0.2cm]
\III'_n(X) & \ = \ d_{\text{K}}\Big(\mathcal{L}(M^{\star}_{k_n}| X\big) \, , \, \mathcal{L}( M_p^{\star} |X\big)\Big).
\end{align}
Note that $\I_n'$ is deterministic, whereas $\II_n'(X)$ and $\III_n'(X)$ are random variables depending on $X$.
The remainder of the proof consists in showing that each of these terms are at most of order $n^{-\frac 12+\delta}$, with probability at least $1-\frac cn$. The terms $\II_n'(X)$ and $\III_n'(X)$ are handled in Sections~\ref{sec:II'} and~\ref{sec:III'} respectively. The first term $\I_n'$ requires no further work, due to Proposition~\ref{prop:IandIII} (since $\I_n'$ is equal to $\III_n$, defined in equation~\eqref{eqn:IIIdef}).  \qed

\subsection{Handling the term $\III_n'(X)$}\label{sec:III'}

The proof of Proposition~\ref{prop:IandIII} can be partially re-used to show that for any fixed realization of $X$, and any real numbers $t_{1,n}'\leq t_{2,n}'$, the following bound holds
\begin{equation}
\III_n'(X) \ \leq \ \P\big(A'(t_{2,n}')\big| X\big) \ + \ \P\big(B'(t_{1,n}')\big| X\big),
\end{equation}
where we define the following events for any $t\in\R$,
\begin{equation}\label{inclusionboot}
\small
A'(t)=\Big\{\max_{j\in\J(k_n)} S^{\star}_{n,j}/\hat \sigma_j^{\tau_n}\leq t\Big\} \ \ \ \text{ and } \ \  \ \ B'(t)=\Big\{\max_{j\in\J(k_n)^c} S^{\star}_{n,j}/\hat \sigma_j^{\tau_n}> t\Big\}.
\end{equation} 
Below, Lemma~\ref{lem:Bboundboot} ensures that $t_{1,n}'$ and $t_{2,n}'$ can be chosen so that 
the random variables $\P(B'(t_{1,n}')\big|X)$ and $\P(A'(t_{2,n}')\big|X)$ are at most $c n^{-\frac 12+\delta}$, with probability at least $1-\frac cn$. Also, it is straightforward to check that under Assumption~\ref{A:cor}, the choices of $t_{1,n}'$ and  $t_{2,n}'$ given in Lemma~\ref{lem:Bboundboot} satisfy $t_{1,n}'\leq t_{2,n}'$ when $n$ is sufficiently large.

\begin{lemma}\label{lem:Bboundboot}
Fix any number $\delta\in(0,1/2)$, and suppose the conditions of Theorem~\ref{THM:G} hold. Then, there are positive constants $c_1$, $c_2$, and $c$, not depending on $n$,  for which the following statement is true: \\[-0.3cm]

If $t_{1,n}'$ and $t_{2,n}'$ are chosen as
\begin{align}
 t_{1,n}' &= c_1\cdot k_n^{-\beta_n}\cdot\log(n)^{3/2}   \ \ \ \ \text{ and }
 \label{eqn:t1prime}\\[0.2cm]
 t_{2,n}' &=c_2\cdot \ell_n^{-\beta_n}\cdot\sqrt{\log(\ell_n)},
\label{eqn:t2prime}
\end{align}
then the events
\begin{equation}\label{abound}\tag{a}
\P(A'(t'_{2,n})\big|X) \leq c\, n^{-\frac 12+\delta}
\end{equation}
and
\begin{equation}\label{bbound}\tag{b}
\P(B'(t'_{1,n})\big| X) \leq n^{-1}
\end{equation}
each hold with probability at least $1-\frac cn$.
\end{lemma}

\paragraph{\textsc{Proof of Lemma~\ref{lem:Bboundboot} part}~\eqref{abound}}Using the definition of $\II_n'(X)$, followed by  $\J(\ell_n)\subset \J(k_n)$, we have
\begin{equation}
\begin{split}
\P(A'(t'_{2,n})|X) 
& \ \leq \ \P\Big(\max_{j\in \J(\ell_n)} \tilde S_{n,j}/\sigma_j^{\tau_n}\leq t'_{2,n}\Big) +\II_n'(X).
\end{split}
\end{equation}
Taking $t_{2,n}'=t_{2,n}$ as in~\eqref{eqn:t2}, the proof of Lemma~\ref{lem:bounds} part \eqref{aboundlem} shows that the first term is $\mathcal{O}(n^{-1/2})$. With regard to the second term, Proposition~\ref{prop:IIprime} in the next subsection shows that  there is a constant $c>0$ not depending on $n$ such that the event
\begin{equation}
\II_n'(X) \leq c\, n^{-\frac 12+\delta}
\end{equation}
holds with probability at least $1-\frac cn$. This completes the proof.\qed


~\\

\paragraph{\textsc{Proof of Lemma~\ref{lem:Bboundboot} part}~\eqref{bbound}} Define the random variable
\begin{equation}\label{eqn:Vstar}
V^{\star}:=\max_{j\in\J(k_n)^c} S_{n,j}^{\star}/\hat\sigma_j^{\tau_n},
\end{equation}
and as in the proof of Lemma~\ref{lem:bounds}\eqref{bboundlem}, let 
$q =\max\big\{\ts\frac{2}{\beta_n},\log(n),3\big\}.$
The idea of the proof is to construct a function $b(\cdot)$ such that the following bound holds for every realization of $X$,
$$\Big( \E\big[|V^{\star}|^q\big | X\big]\Big)^{1/q} \leq  b(X),$$
and then Chebyshev's inequality gives the following inequality for any number $b_n$ satisfying $b(X)\leq b_n$,
$$\P\Big(V^{\star} \geq e b_n\, \Big | X\Big)  \ \leq \ e^{-q} \ \leq \ \ts\frac 1n.$$
In turn, we will derive an expression for $b_n$ such that the event $\{b(X)\leq b_n\}$ holds with high probability.
This will lead to the statement of the lemma, because it will turn out that $t_{1,n}'\asymp b_n$.\\

To construct the function $b(\cdot)$, observe that the initial portion of the proof of Lemma~\ref{lem:bounds}\eqref{bboundlem} shows that for any realization of $X$,
\begin{equation}\label{eqn:Vstarbound}
 \E\big[|V^{\star}|^q\big | X\big]\ \leq   \sum_{j\in\J(k_n)^c} \hat\sigma_j^{q(1-\tau_n)}\,\E\big[|\ts\frac{1}{\hat \sigma_j}S_{n,j}^{\star}|^q|X\big].
\end{equation}
Next, Lemma~\ref{lem:Snjnorm} ensures that for every $j\in\{1,\dots,p\}$, the event
$$ \E\big[|\ts\frac{1}{\hat \sigma_j}S_{n,j}^{\star}|^q|X\big] \ \leq (c\,q)^q, $$
holds with probability 1.
Consequently, if we let $s=q(1-\tau_n)$ and consider the random variable
\begin{equation}
\hat{\mathfrak{s}}:=\Bigg(\sum_{j\in\J(k_n)^c}\hat\sigma_j^{s}\Bigg)^{\frac{1}{s}},
\end{equation}
as well as
$$b(X):=c\cdot q\cdot  \hat{\mathfrak{s}}^{(1-\tau_n)},$$
then we obtain the bound
\begin{equation}
\Big( \E\big[|V^{\star}|^q\big | X\big]\Big)^{1/q} \ \leq b(X),
\end{equation}
with probability 1.
To proceed, Lemma~\ref{lem:mathfrak} implies 
\begin{equation}
\P\bigg(b(X) \geq  q\cdot \ts\frac{(c\sqrt q)^{1-\tau_n}}{(q\beta_n-1)^{1/q}}\cdot k_n^{-\beta_n+1/q}\bigg) \ \leq e^{-q}\leq \frac 1n,
\end{equation}
for some constant $c>0$ not depending on $n$.
By weakening this tail bound slightly, it can be simplified to
\begin{equation}
\P\bigg(b(X) \geq C_n'\cdot  q^{3/2} \cdot k_n^{-\beta_n}\bigg) \ \leq \ts\frac 1n,
\end{equation}
where  $C_n':=\ts\frac{ c \,k_n^{1/q}}{(q\beta_n-1)^{1/q}}$, and we recall $\beta_n=\alpha(1-\tau_n)$.
To simplify further, it can be checked that $C_n'\lesssim 1$, and that the assumption $(1-\tau_n)\sqrt{\log(n)}\gtrsim 1$ gives $q\lesssim \log(n)$. It follows that there is a constant $c$ not depending on $n$ such that if
$$b_n:=c\cdot \log(n)^{3/2} \cdot  k_n^{-\beta_n},$$
then
\begin{equation}
\P(b(X)\geq  b_n)  \ \leq \ts\frac 1n,
\end{equation}
which completes the proof.\qed

\subsection{Handling the term $\II_n'(X)$}\label{sec:II'}
\begin{proposition}\label{prop:IIprime}
Fix any number $\delta\in(0,1/2)$, and suppose the conditions of Theorem~\ref{THM:G} hold. Then, there is a constant $c>0$ not depending on $n$ such that the event
\begin{equation}
\II_n'(X) \leq c\, n^{-\frac 12+\delta}
\end{equation}
holds with probability at least $1-\frac cn$.
\end{proposition}

\proof Define the random variable
\begin{equation}
\breve M_{k_n}^{\star}:=\max_{j\in \J(k_n)} S_{n,j}^{\star}/\sigma_j^{\tau_n},
\end{equation}
which differs from $M_{k_n}^{\star}$, since $\sigma_j^{\tau_n}$ is used in place of $\hat\sigma_j^{\tau_n}$.
 Consider the triangle inequality
 \begin{equation}\label{eqn:II'}
 \small
\II_n'(X) \ \leq \ \dK\Big(\L(\tilde M_{k_n})\, , \L(\breve M_{k_n}^{\star}|X)\Big) \ + \ \dK \Big(\L(\breve M_{k_n}^{\star}|X) \, , \, \L(M_{k_n}^{\star}|X)\Big).
\end{equation}
The two terms on the right will bounded separately. 

To address the first term on the right side of~\eqref{eqn:II'}, we will apply Lemma~\ref{lem:hellinger}, for which a substantial amount of notation is needed. Recall the matrix $\mathfrak{C}=\Sigma^{1/2}\Pi_{k_n}\ttop D_{k_n}^{-\tau_n}$ of size $p\times k_n$, where the projection matrix $\Pi_{k_n}\in\R^{k_n\times p}$ is defined in the proof of Proposition~\ref{prop:bentkus}.  Note that $\tilde M_{k_n}$ is the coordinate-wise maximum of a Gaussian vector drawn from $N(0,\mathfrak{S})$, with $\mathfrak{S}=\mathfrak{C}\ttop\mathfrak{C}$. To address $\breve{M}_{k_n}^{\star}$, let
$$W_n=\ts\frac 1n \sum_{i=1}^n  (Z_i-\bar Z)(Z_i-\bar Z)\ttop$$ 
where $\bar Z=\frac 1n \sum_{i=1}^n Z_i$, and observe that $\breve{M}_{k_n}^{\star}$ is the coordinate-wise maximum of Gaussian vector drawn from $N(0,\breve{\mathfrak{S}})$, with $\breve{\mathfrak{S}}=\mathfrak{C}\ttop W_n\mathfrak{C}$. Next, consider the s.v.d.,
$$\mathfrak{C}=U\Lambda V\ttop,$$
where if $r$ denotes the rank of $\mathfrak{C}$, then we may take $U\in\R^{p\times r}$ to have orthonormal columns, $\Lambda\in\R^{r\times r}$ to be invertible, and $V\ttop\in\R^{r\times k_n}$ to have orthonormal rows. In order to apply Lemma~\ref{lem:hellinger} for a given realization of $\breve{\mathfrak{S}}$, it is necessary that the columns of $\mathfrak{S}$ and $\breve{\mathfrak{S}}$ span the same subspace of $\R^{k_n}$. This occurs with probability at least $1-\frac cn$, because the matrix $\breve{\mathfrak{S}}$ is equal to $V\Lambda (U\ttop W_n U)\Lambda V\ttop$, and the matrix $(U\ttop W_n U)$ is invertible with probability at least $1-\frac cn$ (due to Lemma~\ref{lem:white}). Another ingredient for applying Lemma~\ref{lem:hellinger} is the following algebraic relation, which is a direct consequence of the definitions just introduced,
\begin{equation}
\Big(V\ttop \mathfrak{S} V\Big)^{-1/2}\Big( V\ttop \breve{\mathfrak{S}}V\Big) \Big(V\ttop\mathfrak{S}V\Big)^{-1/2} \ = \ U\ttop W_n U.
\end{equation}
Building on this relation, Lemma~\ref{lem:hellinger} shows that if the event
 \begin{equation}\label{eqn:tempevent2}
 \|U\ttop W_n U-\mathbf{I}_{r}\|_{\text{op}} \ \leq \ \e,
 \end{equation}
 holds for some number $\e>0$, then the event
\begin{equation}
 \dK\Big(\L(\tilde M_{k_n})\, , \L(\breve M_{k_n}^{\star}|X)\Big) \leq c\cdot  k_n^{1/2} \cdot \e
\end{equation}
also holds, where $c>0$ is a constant not depending on $n$ or $\e$. Thus, it remains to specify $\e$ in the event~\eqref{eqn:tempevent2}. For this purpose, Lemma~\ref{lem:white} shows that if $\e=c \cdot n^{-1/2}\cdot k_n\cdot \log(n)$, then the event~\eqref{eqn:tempevent2} holds with probability at least $1-\frac cn$. 
So, given that 
$$n^{-1/2}\cdot k_n^{3/2}\cdot \log(n) \ \lesssim \ n^{-\frac 12+\delta},$$ 
the first term in the bound~\eqref{eqn:II'} requires no further consideration.\\

To deal with the second term in~\eqref{eqn:II'}, we proceed by considering the general inequality
\begin{equation}
\dK(\L(\xi),\L(\zeta)) \ \leq \ \sup_{t\in\R}\P\big(|\zeta-t|\leq r\big)  \ + \  \P(|\xi-\zeta|>r),
\end{equation}
which holds for any random variables $\xi$ and $\zeta$, and any real number $r>0$ (cf.~\citet[Lemma 2.1]{CCK:SPA}). Specifically, we will let $\L(\breve M_{k_n}^{\star}|X)$ play the role of $\L(\xi)$, and let $\L(M_{k_n}^{\star}|X)$ play the role of $\L(\zeta)$. In other words, we need to establish an anti-concentration inequality for $\mathcal{L}(M_{k_n}^{\star}|X)$, as well as a coupling inequality for $ M_{k_n}^{\star}$ and $\breve M_{k_n}^{\star}$, conditionally on $X$.

To establish the coupling inequality, if we put
\begin{equation}
r_n= c \cdot n^{-1/2}\cdot \log(n)^{5/2},
\end{equation}
for a suitable constant $c$ not depending on $n$, then Lemma~\ref{lem:coupling} shows that the event
\begin{equation}\label{eqn:maincouple}
\P\Big(\big|\breve M_{k_n}^{\star} - M_{k_n}^{\star}\big| > r_n\, \Big| X\Big)  \ \leq  \ \ts\frac cn
\end{equation}
holds with probability at least $1-\frac cn$.

Lastly, the anti-concentration inequality can be derived from Nazarov's inequality (Lemma~\ref{lem:Nazarov}), since $M_{k_n}^{\star}$ is obtained from a Gaussian vector, conditionally on $X$. For this purpose, let 
\begin{equation}
\hat{\underline{\sigma}}_{k_n}=\min_{j\in\mathcal{J}(k_n)}\hat\sigma_j.
\end{equation}
In turn, Nazarov's inequality implies that the event
\begin{equation}
\sup_{t\in\R}\, \P\Big(|M_{k_n}^{\star} -t|\leq r_n \, \Big| X\Big) \ \leq \ c\cdot \frac{r_n}{\hat{\underline{\sigma}}_{k_n}^{1-\tau_n}}\cdot \sqrt{\log(k_n)},
\end{equation}
holds with probability 1.
Meanwhile, Lemma~\ref{lem:cor} and Assumption~\ref{A:cor} imply that the event
\begin{equation}\label{eqn:minsigevent}
\frac{1}{\hat{\underline{\sigma}}_{k_n}^{1-\tau_n}} \ \leq c\, k_n^{\beta_n}
\end{equation}
holds with probability at least $1-\frac cn$. 
Combining the last few steps, we conclude that the following bound holds with probability at least $1-\frac cn$,
\begin{equation}
\small
\begin{split}
\sup_{t\in\R}\, \P\Big(|M_{k_n}^{\star} -t|\leq r_n \, \Big| X\Big)  \  & \leq \ts c\cdot n^{-1/2}\cdot k_n^{\beta_n}\cdot \log(n)^{5/2}\cdot\sqrt{\log(k_n)}\\[0.2cm]
& \ \leq  c\,n^{-\frac 12 +\delta},
\end{split}
\end{equation}
as needed.
\qed

\section{Technical lemmas for Theorems~\ref{THM:G} and~\ref{THM:BOOT}}\label{app:technical}
\begin{lemma}\label{lem:sighatrnorm}
Fix any number $\delta\in(0,1/2)$, and suppose the conditions of Theorem~\ref{THM:G} hold. Also, let $q=\max\{\frac{2}{\beta_n}, \log(n),3\}$. Then, there is a constant $c>0$ not depending on $n$, such that for any $j\in\{1,\dots,p\}$, we have
\begin{equation}
\|\hat\sigma_{j}\|_q \ \leq \ c \cdot \sigma_j \cdot \sqrt{q}.
\end{equation}
\end{lemma}
\proof Define the vector $u:=\frac{1}{\sigma_j}\Sigma^{1/2}e_j\in\R^p$, which satisfies $\|u\|_2=1$. Observe that
\begin{equation}\label{eqn:sigsteps}
\begin{split}
\ts\frac{1}{\sigma_j}\|\hat\sigma_{j}\|_q &= \ \Bigg\|\Big(\ts\frac 1n \sum_{i=1}^n (Z_i\ttop u)^2-(\bar Z\ttop u)^2\Big)^{1/2}\Bigg\|_q\\[0.2cm]
&\leq \ \bigg\|\Big(\ts\frac 1n \sum_{i=1}^n ( Z_i\ttop u)^2\Big)^{1/2}\bigg\|_q\\[0.2cm]%
&= \ \Big\|\ts\frac 1n \sum_{i=1}^n (Z_i\ttop u)^2\Big\|_{q/2}^{1/2}\ .
\end{split}
\end{equation}
Since the random variables $(Z_1\ttop u)^2,\dots, (Z_n\ttop u)^2$ are independent~and non-negative, part \emph{(i)} of Rosenthal's inequality in Lemma~\ref{lem:rosenthal} implies the $L^{q/2}$ norm in the last line satisfies
\begin{equation}
\small
 \Big\|\ts\frac 1n \sum_{i=1}^n (Z_i\ttop u)^2\Big\|_{q/2} \ \leq \ c \cdot q \cdot \max\bigg\{\big\|(Z_1\ttop u)^2\big\|_1\, , \, n^{-1+2/q} \big\| (Z_1\ttop u)^2\big\|_{q/2}\bigg\},
\end{equation}
for an absolute constant $c>0$.
For the first term inside the maximum, observe that since $\|u\|_2=1$ and $Z_1$ is isotropic, we have $\|(Z_1\ttop u)^2\|_1=1$. To handle the second term inside the maximum, Assumption~\ref{A:model} implies $\|(Z_1\ttop u)^2\|_{q/2}\lesssim q^2$.
Combining the last few steps, and noticing the square root on the $L^{q/2}$ norm in the last line of~\eqref{eqn:sigsteps}, we obtain 
\begin{equation}
\begin{split}
\ts\frac{1}{\sigma_j}\|\hat\sigma_{j}\|_q & \ \lesssim \sqrt q \cdot \max\Big\{ 1\, , \, n^{-1/2+1/q}q\Big\},
\end{split}
\end{equation}
and this implies the statement of the lemma.\qed

\begin{lemma}\label{lem:mathfrak}
Fix any number $\delta\in(0,1/2)$, and suppose the conditions of Theorem~\ref{THM:G} hold. Also, let $q=\max\{\frac{2}{\beta_n},\log(n),3\}$, and $s=q(1-\tau_n)$, and consider the random variables $\hat{\mathfrak{s}}$ and $\hat{\mathfrak{t}}$ defined by
$$\hat{\mathfrak{s}}=\Bigg(\sum_{j\in\J(k_n)^c} \hat\sigma_{j}^s\Bigg)^{1/s} \text{ \ \ \ \  and \ \ \ \ } \hat{\mathfrak{t}}=\Bigg(\sum_{j\in\J(k_n)} \hat\sigma_{j}^s\Bigg)^{1/s}.$$
Then,
there is a constant $c>0$ not depending on $n$ such that 
\begin{equation}
\P\Bigg(\hat{\mathfrak{s}}\geq \ts\frac{c\sqrt q}{(q\beta_n-1)^{1/s}}\cdot k_n^{-\alpha+1/s}\Bigg) \ \leq \  e^{-q},
\end{equation}
and
\begin{equation}\label{eqn:boundmathfrakt}
\P\Big(\, \hat{\mathfrak{t}}\, \geq \ts\frac{c \sqrt q}{(q\beta_n-1)^{1/s}}\Big) \ \leq \  e^{-q}.
\end{equation}
\end{lemma}
\proof 
In light of the Chebyshev inequality $\P\big(\hat{\mathfrak{s}}\geq e\|\hat{\mathfrak{s}}\|_q\big) \leq e^{-q}$, it suffices to bound $\|\hat{\mathfrak{s}}\|_q$ (and similarly for $\hat{\mathfrak{t}}$). We proceed by direct calculation,
\begin{equation*}
\footnotesize
\begin{split}
\|\hat{\mathfrak{s}}\|_q \ & = \ \Bigg\|\sum_{j\in\J(k_n)^c}\hat\sigma_{j}^s\Bigg\|_{q/s}^{1/s}\\[0.2cm]
&\leq \ \Bigg(\sum_{j\in\J(k_n)^c} \big\|\hat\sigma_{j}^s\big\|_{q/s}\Bigg)^{1/s} \ \ \ \ \ \ \ \text{(triangle inequality for $\|\cdot\|_{q/s}$, with $q/s\geq 1$)}\\[0.2cm]
&=\ \ \Bigg(\sum_{j\in\J(k_n)^c} \big\|\hat\sigma_{j}\big\|_q^s\Bigg)^{1/s}\\[0.2cm]
&\lesssim \ \sqrt{q} \cdot\Bigg( \sum_{j\in\J(k_n)^c} \sigma_j^{s}\Bigg)^{1/s} \ \ \ \ \ \ \ \ \ \text{(Lemma~\ref{lem:sighatrnorm})}\\[0.2cm]
&\lesssim \ \sqrt{q} \cdot \Bigg(\int_{k_n}^p x^{-s\alpha}dx\Bigg)^{1/s}  \\[0.2cm]
&\lesssim \ \sqrt{q}\cdot  \frac{k_n^{-\alpha+1/s}}{(s\alpha-1)^{1/s}},
\end{split}
\end{equation*}
\normalsize
and in the last step we have used the fact that $s\alpha=q\beta_n>1$, which holds since $q$ is defined to satisfy $q\beta_n>1$. The calculation for $\hat{\mathfrak{t}}$ is essentially the same, except that we use $\sum_{j\in\J(k_n)}\sigma_j^s\lesssim 1$.
\qed
~\\

\paragraph{Remark} The following result is a variant of Lemma A.7 in the paper~\citep{Zhilova:2015}.
\begin{lemma}\label{lem:hellinger}
Let $A$ and $B$ be positive semi-definite matrices in \smash{$\R^{d\times d}\setminus\{0\}$} whose columns span the same subspace of $\R^d$.
Define two multivariate normal random vectors $\xi\sim N(0, A)$ and $\zeta\sim N(0,B)$. 
Let $r\leq d$ be the dimension of the subspace spanned by the columns of $A$ and $B$, and let $Q\in\R^{d\times r}$ have columns that are an orthonormal basis for this subspace. Define the $r\times r$ positive definite matrices $\tilde A=Q\ttop A Q$ and $\tilde B=Q\ttop B Q$, and let $H$ be any square matrix satisfying $H\ttop H=\tilde A$. Finally, let $\e>0$ be a number such that $\|(H^{-1})\ttop \tilde B (H^{-1}) -\mathbf{I}_r\|\op \leq \e$. Then, there is an absolute constant $c>0$ such that
\begin{equation}\label{hellenger}
 \sup_{t\in\R}\bigg|\P\Big(\max_{1\leq j\leq d} \xi_j\leq t\Big)-\P\Big(\max_{1\leq j\leq d} \zeta_j\leq t\Big)\bigg| \ \leq  c\sqrt{r}\, \e.
\end{equation}
\end{lemma}

\proof 
We may assume that $\sqrt{r} \,\e\leq 1/2$, for otherwise the claim trivially holds with $c=2$. Define the $r$-dimensional random vectors $\tilde \xi=Q\ttop  \xi$ and $\tilde \zeta=Q\ttop \zeta$. As a consequence of the assumptions, the random vector $\xi$ lies in the column-span of $Q$ almost surely, which gives $Q\tilde \xi=\xi$ almost surely. It follows that for any $t\in\R$, the event $\{\max_{1\leq j\leq d} \xi_j \leq t\}$ can be expressed as $\{\tilde \xi\in\mathcal{A}_t\}$ for some Borel set $\mathcal{A}_t\subset \R^r$. Likewise, we also have  $\{\max_{1\leq j\leq d} \zeta_j \leq t\} = \{\tilde \zeta\in\mathcal{A}_t\}$. Hence, the left hand side of~\eqref{hellenger} is upper-bounded by the total variation distance between $\L(\tilde \xi)$ and $\L(\tilde \zeta)$, and in turn, Pinsker's inequality implies this is upper-bounded by $c\sqrt{d_{\textup{KL}}(\L(\tilde \zeta),\L(\tilde \xi))}$, where $c>0$ is an absolute constant, and $d_{\textup{KL}}$ denotes the KL divergence. Since the random vectors $\tilde \xi$ and $\tilde \zeta$ are Gaussian, the following exact formula is available if we let $\tilde C=(H\ttop)^{-1}\tilde B(H^{-1})-\mathbf{I}_r$, 
\begin{equation}
\begin{split}
 d_{\textup{KL}}(\L(\tilde \zeta),\L(\tilde \xi)) & \ = \ \ts\frac{1}{2}\Big(\tr(\tilde C)-\log\det(\tilde C+\mathbf{I}_r)\Big)\\[0.2cm]
 & \ = \ \ts\frac{1}{2}\sum_{j=1}^r \lambda_j(\tilde C)-\log(\lambda_j(\tilde C)+1).
 \end{split}
\end{equation}
Using the basic inequality  $|x-\log(x+1)|\leq x^2/(1+x)$ that holds for any $x\in(-1,\infty)$, as well as the condition $|\lambda_j(\tilde C)|\leq \e\leq 1/2$, we have
\begin{equation}
\begin{split}
 d_{\textup{KL}}(\L(\tilde \zeta),\L(\tilde \xi)) & \ \leq \ c \, r\, \|\tilde C\|\op^2\\[0.2cm]
 & \ \leq c \, r \, \e^2,
 \end{split}
\end{equation}
for some absolute constant $c>0$. \qed

\begin{lemma}\label{lem:Snjnorm}
Suppose the conditions of Theorem~\ref{THM:G} hold, and let $q=\max\{\frac{2}{\beta_n}, \log(n),3\}$.
 Then, there is a constant $c>0$ not depending on $n$ such that for any $j\in\{1,\dots,p\}$, we have
\begin{equation}\label{eqn:firstsnjnorm}
\|\ts\frac{1}{\sigma_j}S_{n,j}\|_q \leq c\, q,
\end{equation}
and the following event holds with probability 1,
\begin{equation}
\Big(\E\big[|\ts\frac{1}{\hat \sigma_j}S_{n,j}^{\star}|^q|X\big] \Big)^{1/q} \ \leq c\, q.
\end{equation}
\end{lemma}
\proof We only prove the first bound, since the second one can be obtained by repeating the same argument, conditionally on $X$.
Since $q>2$, Lemma~\ref{lem:rosenthal} gives
\begin{equation}\label{eqn:lemrosenthalfirstD}
 \|\ts\frac{1}{\sigma_j}S_{n,j}\|_q \ \lesssim \ q\cdot \max\Big\{  \|\ts\frac{1}{\sigma_j}S_{n,j}\|_2 \, , \,  n^{-1/2+1/q} \|\ts\frac{1}{\sigma_j}(X_{1,j}-\mu_j)\|_q\Big\}.
\end{equation}
Clearly,
\begin{equation}
 \|\ts\frac{1}{\sigma_j}S_{n,j}\|_2^2 \ = \  \var(\ts\frac{1}{\sigma_j}S_{n,j}) \ = \ 1.
\end{equation}
Furthermore, if we define the vector $u:=\ts\frac{1}{\sigma_j}\Sigma^{1/2}e_j$ in $\R^p$, which satisfies $\|u\|_2=1$, then 
\begin{equation}
\begin{split}
\big\|\ts\frac{1}{\sigma_j}(X_{1,j}-\mu_j)\big\|_q
&=\big\| Z_1\ttop u\big\|_q \ \lesssim  \ q \\[0.2cm]
\end{split}
\end{equation}
where the last step  follows from Assumption~\ref{A:model}. Applying the work above to the bound~\eqref{eqn:lemrosenthalfirstD} gives
\begin{equation}
 \|\ts\frac{1}{\sigma_j}S_{n,j}\|_q \ \lesssim \ q\cdot \max\Big\{ 1, n^{-1/2+1/q}\cdot q\Big\}.
\end{equation}
Finally, the stated choice of $q$ implies that the right side in the last display is of order $q$.\qed

\begin{lemma}\label{lem:white}
 Let the random vectors $Z_1,\dots,Z_n\in\R^p$ be as in Assumption~\ref{A:model}, and let $Q\in \R^{p\times r}$ be a fixed matrix having orthonormal columns with $r\leq k_n$. Lastly, let
 \begin{equation}\label{eqn:Wndef}
 W_n=\frac 1n \sum_{i=1}^n  (Z_i-\bar Z)(Z_i-\bar Z)\ttop,
 \end{equation}
 where $\bar Z=\frac 1n \sum_{i=1}^n Z_i$. Then, there is a constant $c>0$ not depending on $n$, such that the event
  \begin{equation}\label{eqn:opnormevent}
 \big\|Q\ttop W_nQ-\mathbf{I}_{r}\big\|\op \ \leq \ \ts\frac{c\log(n)k_n}{\sqrt n},
 \end{equation}
 holds with probability at least $1-\frac cn$. 
\end{lemma}
\proof 
Let $\e\in (0,1/2)$, and let $\mathcal{N}$ be an $\e$-net (with respect to the $\ell_2$-norm) for the unit $\ell_2$-sphere in $\R^{r}$. It is well known that $\mathcal{N}$ can be chosen so that $\text{card}(\mathcal{N})\leq (3/\e)^{r}$, and the inequality
$$\big\|Q\ttop W_n Q- \mathbf{I}_{r}\big\|\op \leq \ts\frac{1}{1-2\e}\cdot \displaystyle \max_{u\in\mathcal{N}} \Big|u\ttop\Big(Q\ttop W_n Q-\mathbf{I}_{r}\Big)u\Big|,$$
holds with probability 1~\citep[Lemmas 5.2 and 5.4]{Vershynin:2012}. For a fixed $u\in\mathcal{N}$, put $\xi_{u,i}:=Z_i\ttop Qu$,
and consider the simple algebraic relation
\begin{equation}
u\ttop\Big(Q\ttop W_n Q-\mathbf{I}_{r}\Big)u = \underbrace{\Big(\ts\frac{1}{n}\sum_{i=1}^n \xi_{i,u}^2-1\Big)}_{=: \Delta(u)}- \underbrace{\Big(\ts\frac{1}{n}\sum_{i=1}^n \xi_{i,u}\Big)^2}_{=:\Delta'(u)}.
\end{equation}
We will show that both terms on the right side are small with high probability, and then take a union bound over $u\in\mathcal{N}$. The high-probability bounds will be obtained by using Lemma~\ref{lem:rosenthal} to control $\|\Delta(u)\|_q$ and $\|\Delta'(u)\|_q$ when $q$ is sufficiently large.\\

To apply Lemma~\ref{lem:rosenthal}, first observe the following bounds, which are consequences of Assumption~\ref{A:model},
\begin{equation}
\|\xi_{i,u}\|_q \ \lesssim \  q, 
\end{equation}
and
\begin{equation}
\begin{split}
\|\xi_{i,u}^2-1\|_q \ \lesssim \ q^2.
\end{split} 
\end{equation}
Therefore, when $q>2$, Lemma~\ref{lem:rosenthal} gives
\begin{equation}\label{eqn:justbefore}
\begin{split}
\|\Delta(u)\|_q  & \ \lesssim \ q \max\Big\{ \|\Delta(u)\|_2 \, , \, \ts\frac{1}{n}\big(\tsum_{i=1}^n \|\xi_{i,u}^2-1\|_q^q\big)^{1/q}\Big\}\\[0.3cm]
& \lesssim \ q \max\Big\{ \ts\frac{1}{\sqrt n} \, , \, n^{-1+1/q} \cdot q^2\Big\}.
\end{split}
\end{equation} 
Due to Chebyshev's inequality,
$$\P\Big( |\Delta(u)| \ \geq e \|\Delta(u)\|_q\Big) \ \leq \ e^{-q},$$
and so if we take $q= \max\{C\log(n)k_n,3\}$ for some constant $C>0$ to be tuned below, then~\eqref{eqn:justbefore} gives $\|\Delta(u)\|_q\lesssim q/\sqrt{n}$, and the following inequality holds for all $u\in\mathcal{N}$,
$$\P\Big( |\Delta(u)| \ \geq \ts\frac{c\log(n) k_n}{\sqrt n}\Big) \ \leq \ \exp\Big\{-C\log(n) k_n\Big\}.$$
The random variable $\Delta'(u)$ can be analyzed with a similar set of steps, which leads to the following inequality for all $u\in\mathcal{N}$,
\begin{equation}
\P\bigg( |\Delta'(u)| \ \geq c\Big(\ts\frac{\log(n) k_n}{\sqrt n}\Big)^2\bigg) \ \leq \ \exp\Big\{-C\log(n) k_n\Big\}.
\end{equation}
Combining the previous work with a union bound, if we consider the choice $\e=\min\{c\log(n)k_n/\sqrt{n},\ts\frac{1}{4}\}$, then 
$$ \P\Big(\big\|Q\ttop W_n Q- \mathbf{I}_{r}\big\|\op \geq \e\Big) \ \leq \ 2\exp\Big\{ - C\cdot k_n\cdot\log(n) \ + \ r\cdot \log(3/\e)\Big\}.$$
Finally, choosing $C$ sufficiently large implies the stated result.\qed
 
\paragraph{Remark} For the next results, define the correlation 
$$\rho_{j,j'}=\ts\frac{\Sigma_{j,j'}}{\sigma_j\sigma_{j'}},$$ 
and its sample version 
$$\hat{\rho}_{j,j'}=\ts\frac{\hat\Sigma_{j,j'}}{\hat\sigma_j\hat\sigma_{j'}},$$
 for any $j,j'\in\{1,\dots,p\}$.
\begin{lemma}\label{lem:cor}
Suppose the conditions of Theorem~\ref{THM:G} hold. Then, there is a constant $c>0$ not depending on $n$ such that the three events
\begin{equation}\label{eqn:firstcor}
 \max_{j\in\J(k_n)}\Big| \ts\frac{\hat\sigma_j}{\sigma_j}-1\Big| \leq \ts \frac{c\log(n)}{\sqrt n},
\end{equation}
\begin{equation}
 \min_{j\in\J(k_n)}\hat\sigma_j^{1-\tau_n} \ \geq \ \Big(\min_{j\in\J(k_n)} \sigma_j^{1-\tau_n}\Big)\cdot \Big(1-\ts\frac{c\log(n)}{\sqrt n}\Big),
\end{equation}
and
\begin{equation}\label{eqn:firstcor}
 \max_{j,j'\in\J(k_n)}\big| \hat \rho_{jj'}-\rho_{jj'}\big| \leq \ts \frac{c\log(n)}{\sqrt n}
\end{equation}
each
hold with probability at least $1-\ts\frac{c}{n}$.
\end{lemma}
\proof The result is a direct consequence of Lemma~\ref{lem:corbasic} below. The details are essentially algebraic manipulations, and so are omitted.\qed

\begin{lemma}\label{lem:corbasic}
Suppose the conditions of Theorem~\ref{THM:G} hold, and fix any two (possibly equal) indices $j,j'\in\{1,\dots,p\}$.
 Then, for any number $\kappa\geq 1$, there are positive constants $c$ and $c_1(\kappa)$ not depending on $n$ such that the event
\begin{equation}\label{eqn:corevent}
 \Big| \ts\frac{\hat \Sigma_{j,j'}}{\sigma_j\sigma_{j'}}-\rho_{j,j'}\Big| \leq \frac{c_1(\kappa)\log(n)}{\sqrt n}
\end{equation}
holds with probability at least $1-cn^{-\kappa}$. 
\end{lemma}
\paragraph{Remark}The event in the lemma has been formulated to hold with probability at least $1-cn^{-\kappa}$, rather than $1-\frac{c}{n}$, in order to accommodate a union bound for proving Lemma~\ref{lem:cor}.

\proof
 Consider the $\ell_2$-unit vectors $u=\Sigma^{1/2}e_j/\sigma_j$ and $v=\Sigma^{1/2}e_{j'}/\sigma_{j'}$  in $\R^p$.
Letting $W_n$ be as defined in~\eqref{eqn:Wndef}, observe that
\begin{equation}\label{eqn:decomp}
\begin{split}
\frac{\hat \Sigma_{j,j'}}{\sigma_j\sigma_{j'}} -\rho_{j,j'}
 & \ = u\ttop (W_n-\mathbf{I}_p)v.
\end{split}
\end{equation}
For each $1\leq i\leq n$, define the random variables $\zeta_{i,u}=Z_i\ttop u$ and $\zeta_{i,v}=Z_i\ttop v$. In this notation, the relation~\eqref{eqn:decomp} becomes
\begin{equation*}
\frac{\hat \Sigma_{j,j'}}{\sigma_j\sigma_{j'}} -\rho_{j,j'} =\underbrace{\Big(\ts\frac 1n \sum_{i=1}^n \zeta_{i,u}\zeta_{i,v}- u\ttop v\Big)}_{=:\Delta(u,v)} \ - \ \underbrace{\Big(\ts\frac 1n \sum_{i=1}^n \zeta_{i,u}\Big)\Big(\ts\frac 1n \sum_{i=1}^n \zeta_{i,v}\Big)}_{=:\Delta'(u,v)}.
\end{equation*}
Note that $\E[\zeta_{i,u}\zeta_{i,v}]=u\ttop v$.  Also, if we let $q=\max\{\kappa\log(n),3\}$, then
$$\|\zeta_{i,u}\zeta_{i,v}-u\ttop v\|_q \ \lesssim \ q^2,$$
which follows from Assumption~\ref{A:model}. Therefore, Lemma~\eqref{lem:rosenthal} gives the following bound for $q>2$,
\begin{equation}
\begin{split}
\|\Delta(u,v)\|_q & \ \lesssim  \ q \max\Big\{ \|\Delta(u,v)\|_2 \, , \, \ts\frac{1}{n}\big(\sum_{i=1}^n \|\zeta_{i,u}\zeta_{i,v}-u\ttop v\|_q^q\big)^{1/q}\Big\}\\[0.3cm]
& \ \lesssim \ q \max\Big\{\ts\frac{1}{\sqrt n} \, , \,  n^{-1+1/q}\cdot q^2\Big\}\\[0.3cm]
& \ \lesssim \ \ts\frac{\log(n)}{\sqrt n}.
\end{split}
\end{equation}
Using the Chebyshev inequality 
$$\P(|\Delta(u,v)|\geq e \|\Delta(u,v)\|_q) \ \leq e^{-q},$$
 we have
 $$\P\Big(|\Delta(u,v)|\geq \ts\frac{c\kappa\log(n)}{\sqrt n}\Big) \ \leq \ \ts\frac{1}{n^{\kappa}}.$$
Similar reasoning leads to the following tail bound for $\Delta'(u,v)$, 
  $$\P\bigg(|\Delta'(u,v)|\geq \big(\ts\frac{c\kappa\log(n)}{\sqrt n}\big)^2\bigg) \ \leq \ \ts\frac{1}{n^{\kappa}},$$
and combining with the previous tail bound gives the stated result.\qed

\begin{lemma}\label{lem:coupling}
Fix any number $\delta\in(0,1/2)$, and suppose the conditions of Theorem~\ref{THM:G} hold. Then, there is a constant $c>0$ not depending on $n$ such that the event~\eqref{eqn:maincouple} holds with probability at least $1-\frac cn$.
\end{lemma}

\proof Let $(a_1,\dots,a_{k_n})$ and $(b_1,\dots,b_{k_n})$ be real vectors, and note the basic fact
$$\bigg|\max_{1\leq j\leq k_n} a_j -\max_{1\leq j\leq k_n}b_j\bigg| \ \leq \ \max_{1\leq j\leq k_n} |a_j-b_j|.$$
From this, it is simple to derive the inequality
\begin{equation}
\big| \breve M_{k_n}^{\star}-M_{k_n}^{\star}\big| \leq \max_{j\in\J(k_n)}\Big|\big(\ts\frac{\hat\sigma_j}{\sigma_j}\big)^{\tau_n}-1\Big|\cdot \displaystyle \max_{j\in\J(k_n)}\Big|S_j^{\star}/\hat\sigma_j^{\tau_n}\Big|.
\end{equation}
To handle the first factor on the right side, it follows from Lemma~\ref{lem:corbasic} that the event
\begin{equation}
\max_{j\in\J(k_n)}\Big|\big(\ts\frac{\hat\sigma_j}{\sigma_j}\big)^{\tau_n}-1\Big| \ \leq c\cdot n^{-1/2}\cdot \log(n)
\end{equation}
holds with probability at least $1-\frac cn$. Next,  consider the random variable
\begin{equation}\label{eqn:Ustar}
U^{\star}:=\max_{j\in\J(k_n)} |S_{n,j}^{\star}/\hat\sigma_j^{\tau_n}|.
\end{equation}
It suffices to show there is possibly larger constant $c>0$, such that the event
\begin{equation}\label{eqn:Ustarbound}
\P\Big( U^{\star} \geq  c \log(n)^{3/2} \Big| X\Big) \leq \ts\frac 1n 
\end{equation}
holds with probability at least $1-\frac cn$.
Using Chebyshev's inequality  with $q\geq \log(n)$ gives
 $$\P\Big(U^{\star} \geq e\, \big(\E[|U^{\star}|^q|X]\big)^{1/q}\Big| X\Big)\leq e^{-q}.$$
Likewise, if the event
 \begin{equation}\label{eqn:UstarLq}
 (\E[|U^{\star}|^q|X])^{1/q} \ \leq c\log(n)^{3/2}
 \end{equation}
holds for some constant $c>0$, then the event~\eqref{eqn:Ustarbound} also holds. For this purpose, the argument in the proof of Lemma~\ref{lem:Bboundboot}\eqref{bbound} can be essentially repeated with $q=\max\{\frac{2}{\beta_n},\log(n),3\big\}$ to show that the event~\eqref{eqn:UstarLq} holds with probability at least $1-\frac cn$. The main detail to notice when repeating the argument is that $U^{\star}$ involves a maximum over $\J(k_n)$, whereas the argument for Lemma~\ref{lem:Bboundboot}\eqref{bbound} involves a maximum over $\J(k_n)^c$. This distinction can be handled by using the bound~\eqref{eqn:boundmathfrakt} in Lemma~\ref{lem:mathfrak}.\qed

\paragraph{The case when \lowercase{ $p\leq k_n$}}\label{lowdimcase}
The previous proofs relied on the condition $p>k_n$ only insofar as this implies $k_n\geq n^{\frac{1}{\log(n)^a}}$ and $\ell_n\geq \log(n)^3$. (These conditions are used in the analyses of $\I_n$ and $\III_n$, as well as $\I_n'$ and $\III_n'(X)$.)  However, if $p\leq k_n$, then the definition of $k_n$ implies that $p=k_n$, which causes the quantities $\I_n$, $\III_n$, $\I_n'(X)$ and $\III_n'(X)$ to become exactly 0. In this case, the proofs of Theorems~\ref{THM:G} and~\ref{THM:BOOT} reduce to bounding $\II_n$ and $\II_n'(X)$, and these arguments can be repeated as before.

\section{Proof of Theorem~\ref{THM:M}}\label{app:mult}
\paragraph{Notation and remarks}  
An important piece of notation for this appendix (and the next one) is the integer $d_n$, which we define to be the largest index $d_n\in\{1,\dots,p\}$ such that $\sigma_{(d_n)}^2\geq \e_0^2/\sqrt{n}$, with $\e_0$ as in Assumption~\ref{A:mult}. (Such an index must exist under Assumption~\ref{A:mult}.)
Also, it is simple to check that $k_n\leq d_n$ holds for all large $n$ under Assumption~\ref{A:mult}. Furthermore, as in the previous appendices, we will assume $p>k_n$, since the case $p\leq k_n$ can be handled using similar reasoning to that explained in the previous paragraph. Lastly, it will be helpful to note that $\sigma_{(j)}^2=\pi_{(j)}(1-\pi_{(j)})$, since it can be checked that $\pi_i\leq \pi_j$ implies $\sigma_i\leq \sigma_j$.

\paragraph{Outline of proof} Since the number $d_n$ will often play the role that $p$ did in previous proofs, we will use a slightly different notation for the analogues of the earlier quantities $\I_n$, $\II_n$, $\II_n'(X)$, and $\III_n'(X)$. This will also serve as a reminder that new details are involved in the context of the multinomial model. The new quantities are:
\begin{align*}
\bI_n &= d_{\text{K}}\big(\mathcal{L}(M_{d_n}) \, , \, \mathcal{L}(M_{k_n})\big)\\[0.2cm]
\bII_n &=d_{\text{K}}\big( \mathcal{L}(M_{k_n})\, , \, \mathcal{L}(\tilde M_{k_n})\big)\\[0.2cm]
\bII'_n(X) &\ = \ d_{\text{K}}\big(\mathcal{L}(\tilde M_{k_n}) \, , \, \mathcal{L}( M^{\star}_{k_n}|X\big)\big)\\[0.2cm]
\bIII'_n(X) & \ = \ d_{\text{K}}\big(\mathcal{L}(M^{\star}_{k_n}| X\big) \, , \, \mathcal{L}( M_{d_n}^{\star} |X\big)\big).
\end{align*}
The overall structure of the proof is based on the simple bound
\begin{equation*}\label{eqn:threepartsM}
\begin{split}
   d_{\textup{K}}\big(\L(\M)\,,\,\L(\M^{\star}|X)\big) \ \leq & \  \ d_{\textup{K}}\big(\L(\M) \, ,\,\L(M_{d_n})\big)\\[0.2cm]
   & \ + \bI_n +\bII_n+\bII_n'(X)+\bIII_n'(X)\\[0.2cm]
   & \ + d_{\textup{K}}\big(\L(M_{d_n}^{\star}|X) \, , \, \L(\M^{\star}|X)\big).
   \end{split}
\end{equation*}
Most of the proof will be completed through four separate lemmas, showing that each of the quantities in Roman numerals is at most of order $n^{-1/2+\delta}$. These lemmas are labeled as~\ref{lem:IIM} (for $\bII_n$),~\ref{lem:IM} (for $\bI_n$),~\ref{lem:IIprimeM} (for $\bII'_n(X)$) and~\ref{lem:IIIprimeM} (for $\bIII_n'(X)$). Finally, the other two quantities in the first and third lines %
will be shown to be at most of order $n^{-1/2+\delta}$ in Lemma~\ref{lem:otherM}. \qed
~\\

\begin{lemma}\label{lem:IIM}Fix any number $\delta\in(0,1/2)$, and suppose the conditions of Theorem~\ref{THM:M} hold. Then,
\begin{equation}
\bII_n \ \lesssim \  n^{-1/2+\delta}.
\end{equation}
\end{lemma}
\proof The argument is rougly similar to the proof of Proposition~\ref{prop:bentkus}, and we retain some of the notation used there. For the new proof, define the $k_n\times k_n$ matrix 
\begin{equation}\label{eqn:sigmakdef}
\Sigma(k_n):=\Pi_{k_n}\Sigma\Pi_{k_n}\ttop,
\end{equation}
where $\Pi_{k_n}\in\R^{k_n\times p}$ is the projection onto the coordinates indexed by $\J(k_n)$, as explained on page~\pageref{eqn:IInbound}. Let $r_n$ denote the rank of $\Sigma(k_n)$, and let 
$$\Sigma(k_n) = Q\Lambda_{r_n} Q\ttop$$
be a spectral decomposition of $\Sigma(k_n)$, where $Q\in\R^{k_n\times r_n}$ has orthonormal columns, and $\Lambda\in\R^{r_n\times r_n}$ is diagonal and invertible. In addition, define
$$\mathfrak{C}=\Lambda_{r_n}^{1/2}Q\ttop D_{k_n}^{-\tau_n},$$
as well as the $r_n$-dimensional random vector
$$Z_i'=\Lambda_{r_n}^{-1/2}Q\ttop \Pi_{k_n}(X_i-\boldsymbol\pi),$$
and the sample average
$$\breve Z_n'= \ts\frac{1}{\sqrt n}\sum_{i=1}^n Z_i'.$$
Since the random vector $\Pi_{k_n}(X_i-\boldsymbol \pi)$ lies in the column span of $\Sigma(k_n)$ almost surely, it follows that the relation $QQ\ttop \Pi_{k_n}(X_i-\boldsymbol\pi)=\Pi_{k_n}(X_i-\boldsymbol\pi)$ holds almost surely. In turn, this gives
$$Q\Lambda_{r_n}^{1/2} Z_i' = \Pi_{k_n}(X_i-\boldsymbol \pi),$$
and hence
 $$D_{k_n}^{-\tau_n}\Pi_{k_n}S_n = \mathfrak{C}\ttop\breve Z_n'.$$
 Consequently, for any $t\in\R$, there is a Borel convex set $\mathcal{A}_t\subset\R^{r_n}$ such that 
$$\P\big(M_{k_n}\leq t) \ = \  \P\Big(\max_{j\in\J(k_n)}D_{k_n}^{-\tau_n}\Pi_{k_n}S_{n,j}  \leq t\Big) \ = \  \P\big(\breve Z_n'\in\mathcal{A}_t\big).$$
Similar reasoning also can be applied to $\tilde M_{k_n}$ in order to obtain the expression \smash{$\P(\tilde M_{k_n}\leq t)= \gamma_{r_n}(\mathcal{A}_t)$}, where $\gamma_{r_n}$ is the standard Gaussian distribution on $\R^{r_n}$. 
Therefore, we have the bound
\begin{equation}
\bII_n \ \leq  \ \sup_{\mathcal{A}\in\mathscr{A}}\Big| \P(\breve Z' \in\mathcal{A}) - \gamma_{r_n}(\mathcal{A})\Big|,
\end{equation}
with $\mathscr{A}$ being the collection of Borel convex subsets of $\R^{r_n}$.

Since the vectors $Z_1',\dots,Z_n'$ are i.i.d., with mean 0 and identity covariance matrix, we may apply the Berry-Esseen bound of Bentkus~(Lemma~\ref{lem:bentkus}). The only remaining detail is to bound $\E[\|Z_1'\|_2^3]$, and show that it is at most a fixed power of $k_n$. To do this, Lemma~\ref{lem:lambdamin} implies there is a constant $c>0$ not depending on $n$ such that the following inequalities hold with probability 1,
\begin{equation*}
\begin{split}
\|Z_1'\|_2^2  & \ \leq \ \ts\frac{1}{\lambda_{r_n}(\Sigma(k_n))}\big\|\Pi_{k_n}(X_i-\boldsymbol\pi)\big\|_2^2\\[0.2cm]
& \ \leq \ ck_n^c \cdot k_n,
\end{split}
\end{equation*}
where we have used the facts that $\|Q\|\op\leq 1$, and the entries of $\Pi_{k_n}(X_i-\boldsymbol\pi)$ are bounded in magnitude by 1. Thus, we have $\E[\|Z_1'\|_2^3]\lesssim k_n^{(3/2)(c+1)}$, which completes the proof.\qed

\begin{lemma}\label{lem:IM}
Fix any number $\delta\in(0,1/2)$, and suppose the conditions of Theorem~\ref{THM:M} hold. Then,
\begin{equation}
\bI_n \ \lesssim \ n^{-1/2+\delta}.
\end{equation}
\end{lemma}
\proof By repeating the proof of Proposition~\ref{prop:IandIII}, it follows that if $\texttt{t}_{1,n}$ and $\texttt{t}_{2,n}$ are any two real numbers with $\texttt{t}_{1,n}\leq \texttt{t}_{2,n}$, then 
\begin{equation}
\bI_n \ \leq \ \P(\texttt{A}(\texttt{t}_{2,n})) \ + \ \P(\texttt{B}(\texttt{t}_{1,n})),
\end{equation}
where we write
\small
\begin{equation*}
\texttt{A}(t)=\Big\{\max_{j\in\J(k_n)} S_{n,j}/\sigma_j^{\tau_n} \leq t\Big\} \ \ \ \text{ and } \ \  \ \ \texttt{B}(t)=\Big\{\max_{j\in\J(d_n)\setminus\J(k_n)} S_{n,j}/\sigma_j^{\tau_n}> t\Big\},
\end{equation*}
\normalsize
for any $t\in \R$. (Note that $\texttt{B}(t)$ differs from $B(t)$ only insofar as $\J(d_n)\setminus\J(k_n)$ replaces $\J(k_n)^c$.)

To handle the probability $\P(\texttt{A}(\texttt{t}_{2,n}))$, we mimic the definition~\eqref{eqn:t2} and let
\begin{equation}
 \texttt{t}_{2,n}=\e_0\cdot \ell_n^{-\beta_n}\cdot \sqrt{\log(\ell_n)},  
\end{equation}
with $\e_0\in(0,1)$ as in Assumption~\ref{A:mult}. Having made this choice, the proof of Lemma~\ref{lem:bounds}\eqref{aboundlem} may be repeated essentially verbatim to show that $\P(\texttt{A}(\texttt{t}_{2,n}))\lesssim n^{-1/2+\delta}$. In particular, it is important to note that the correlation matrix $R(\ell_n)$ in the multinomial case satisfies the conditions needed for that argument to work, because $R^+(\ell_n)=\mathbf{I}_{\ell_n}$. Also, this argument relies on Lemma~\ref{lem:IIM} in the same way that the proof of Lemma~\ref{lem:bounds}\eqref{aboundlem} relies on Proposition~\ref{prop:bentkus}.

To handle $\P(\texttt{B}(\texttt{t}_{1,n}))$, the proof of Lemma~\ref{lem:bounds}\eqref{bboundlem} can be mostly repeated to show that this probability is of order $1/n$. However, there are a few differences. First, we may regard the $\alpha$ from the context of Lemma~\ref{lem:bounds}\eqref{bboundlem} as being equal to $1/2$, due to the basic fact that $\sigma_{(j)}\leq  j^{-1/2}$ always holds in the multinomial model. Likewise, we define
 \begin{equation}
 \texttt{t}_{1,n}=c\cdot k_n^{-(1-\tau_n)/2}\cdot \log(n)  
\end{equation}
as the analogue of $t_{1,n}$ in~\eqref{eqn:t1}. The only other issues to notice are that the set $\J(k_n)^c$ is replaced by $\J(d_n)\setminus\J(k_n)$, and that we must verify the following condition. Namely, there is a constant $c>0$ not depending on $n$ such that the inequality
 $$\max_{j\in\J(d_n)} \big\|\ts\frac{1}{\sigma_j}S_{n,j}\big\|_q \ \leq c q$$
 holds when $q=\max\{\frac{2}{(1/2)(1-\tau_n)},\log(n),3\}$. This will be verified later in Lemma~\ref{lem:SnjnormM}. Finally, it is simple to check that  $\texttt{t}_{1,n}  \leq \texttt{t}_{2,n}$ holds for all large $n$. \qed
~\\

\begin{lemma}\label{lem:IIprimeM}
Fix any number $\delta\in(0,1/2)$, and suppose the conditions of Theorem~\ref{THM:M} hold. Then, there is a constant $c>0$ not depending on $n$ such that the event
\begin{equation}
\bII_n'(X) \ \leq \ c\, n^{-1/2+\delta}
\end{equation}
holds with probability at least $1-\frac cn$.
\end{lemma}

\proof Let the random variable $\breve M_{k_n}^{\star}$ be as defined in the proof of Proposition~\ref{prop:IIprime}, and consider the triangle inequality
 \begin{equation}\label{eqn:IIprimeM}
 \small
\bII_n'(X) \ \leq \ \dK\Big(\L(\tilde M_{k_n})\, , \L(\breve M_{k_n}^{\star}|X)\Big) \ + \ \dK \Big(\L(\breve M_{k_n}^{\star}|X) \, , \, \L(M_{k_n}^{\star}|X)\Big).
\end{equation}
Regarding the second term on the right, the proof of Proposition~\ref{prop:IIprime} shows that this term can be controlled with a coupling inequality for $\breve{M}_{k_n}^{\star}$ and $M_{k_n}^{\star}$, as well as an anti-concentration inequality for $M_{k_n}^{\star}$. In the multinomial context, both of these inequalities can be established using the same overall approach as before. The main items that need to be updated are to replace $\J(k_n)^c$ with $\J(d_n)\setminus \J(k_n)$, and to control quantities involving \smash{$\{\hat\sigma_j | j\in\J(d_n)\}$} by using Lemmas~(\ref{lem:corM},~\ref{lem:sighatrnormM}) instead of Lemmas~(\ref{lem:cor},~\ref{lem:sighatrnorm}).

Next, to control the first term in the bound~\eqref{eqn:IIprimeM}, we will use an argument based on  Lemma~\ref{lem:hellinger}, which requires quite a few pieces of notation. First, let
$\Sigma(k_n)=\Pi_{k_n}\Sigma\Pi_{k_n}\ttop \in\R^{k_n\times k_n}$, and let the rank of this matrix be denoted by $r_n$. Next, write the spectral decomposition of $\Sigma(k_n)$ as
 $$\Sigma(k_n)=Q\Lambda_{r_n}Q\ttop,$$
 where $Q\in\R^{k_n\times r_n}$ has orthonormal columns and $\Lambda_{r_n}\in\R^{r_n\times r_n}$ is diagonal. In addition define 
\begin{equation}\label{eqn:whitemultdef}
\begin{split}
\hat\Sigma(k_n)&=\Pi_{k_n}\hat\Sigma \Pi_{k_n}\ttop\\[0.2cm]
W_n&= \Lambda_{r_n}^{-1/2}Q\ttop \hat\Sigma(k_n) Q\Lambda_{r_n}^{-1/2}.
\end{split}
\end{equation}
Since each vector $\Pi_{k_n}(X_i-\bar X)$ lies in the column span of $\Sigma(k_n)$ almost surely, it follows that the relation $QQ\ttop\hat\Sigma(k_n)QQ\ttop=\hat\Sigma(k_n)$ holds almost surely, which is equivalent to
$$\hat\Sigma(k_n)=Q\Lambda_{r_n}^{1/2}W_n\Lambda_{r_n}^{1/2}Q\ttop.$$
With this in mind,  let $\mathfrak{C}=\Lambda_{r_n}^{1/2}Q\ttop D_{k_n}^{-\tau_n}$, and also define
\begin{equation}
\begin{split}
\mathfrak{S}&=\mathfrak{C}\ttop \mathfrak{C}\\[0.2cm]
\breve{\mathfrak{S}} &= \mathfrak{C}\ttop W_n\mathfrak{C}.
\end{split}
\end{equation}
It is straightforward to check that $\tilde M_{k_n}$ is the coordinate-wise maximum of a Gaussian vector drawn from $N(0,\mathfrak{S})$, and similarly, $\breve M_{k_n}^{\star}$ is the coordinate-wise maximum of a Gaussian vector drawn from $N(0,\breve{\mathfrak{S}})$.

We will now compare $\tilde M_{k_n}$ and $\breve M_{k_n}^{\star}$ by applying Lemma~\ref{lem:hellinger}. For this purpose, let
\begin{equation}\label{eqn:Udef}
\mathfrak{C}=ULV\ttop
\end{equation}
be an s.v.d.~for $\mathfrak{C}$, where the matrix $L\in\R^{r_n\times r_n}$ is diagonal and invertible, and the matrices $U\in\R^{r_n\times r_n}$ and $V\in\R^{k_n\times r_n}$ each have orthonormal columns. From these definitions, it follows that
\begin{equation}
\Big(V\ttop \mathfrak{S} V\Big)^{-1/2}\Big( V\ttop \breve{\mathfrak{S}}V\Big) \Big(V\ttop\mathfrak{S}V\Big)^{-1/2} \ = \ U\ttop W_n U.
\end{equation}
Thus, the matrix $(V\ttop \mathfrak{S} V)^{1/2}$ will play the role of $H$ in the statement of Lemma~\ref{lem:hellinger}.
Also, in order to apply that lemma, we need that the columns of $\mathfrak{S}$ and $\breve{\mathfrak{S}}$ span the same subspace of $\R^{k_n}$ (with high probability). Noting that $\mathfrak{S}=VL^2 V\ttop$ and $\breve{\mathfrak{S}} =VL (U\ttop W_n U)L V\ttop$, it follows that $\mathfrak{S}$ and $\breve{\mathfrak{S}}$ have the same column span whenever $U\ttop W_n U$ is invertible, and by Lemma~\ref{lem:whitemult}, this holds with probability at least $1-c/n$. Therefore, Lemma~\ref{lem:hellinger} ensures that if the event
 \begin{equation}\label{eqn:tempevent}
 \|U\ttop W_n U-\mathbf{I}_{r_n}\|_{\text{op}} \ \leq \ \e,
 \end{equation}
 holds for some number $\e>0$, then the event
\begin{equation}
 \dK\Big(\L(\tilde M_{k_n})\, , \L(\breve M_{k_n}^{\star}|X)\Big) \leq c\cdot  k_n^{1/2} \cdot \e
\end{equation}
also holds. Finally, Lemma~\ref{lem:whitemult} shows that if we take $\e= c\log(n)k_n^c/\sqrt{n}$, then the event~\eqref{eqn:tempevent} holds with probability at least $1-c/n$, which completes the proof.\qed

\begin{lemma}\label{lem:IIIprimeM}
Fix any number $\delta\in(0,1/2)$, and suppose the conditions of Theorem~\ref{THM:M} hold. Then, there is a constant $c>0$ not depending on $n$ such that the event
\begin{equation}
\bIII_n'(X) \ \leq \ c\, n^{-1/2+\delta}
\end{equation}
holds with probability at least $1-\frac cn$.
\end{lemma}

\proof The proof follows the argument outlined in Section~\ref{sec:III'}. The main details to be updated for the multinomial context arise in controlling the quantities $\{\hat\sigma_j | j\in\J(d_n)\}$. Specifically, the index set $\J(k_n)^c$ must be replaced with  $\J(d_n)\setminus \J(k_n)$, and Lemmas~(\ref{lem:corM},~\ref{lem:sighatrnormM}) must be used in place of Lemmas~(\ref{lem:cor},~\ref{lem:sighatrnorm}).\qed

\begin{lemma}\label{lem:otherM}
Suppose the conditions of Theorem~\ref{THM:M} hold. Then, there is a constant $c>0$ not depending on $n$, such that
\begin{equation}\label{eqn:reductionfirst}
 d_{\textup{K}}\big(\L(\M), \L(M_{d_n})\big) \ \leq \ \bI_n+\ts\frac{c}{n},
\end{equation}
and the event
\begin{equation}\label{eqn:reductionsecond}
  d_{\textup{K}}\big(\L(\M^{\star}|X), \L(M_{d_n}^{\star}|X)\big) \ \leq \   \bIII_n'(X) 
\end{equation}
occurs with probability at least $1-\frac cn$.
\end{lemma}
\proof Fix any $t\in\R$. By intersecting the event $\{\max_{j\in\hat\J_n}S_{n,j}/\sigma_j^{\tau_n}\leq t\}$ with the events $\{\J(k_n)\subset \hat\J_n\}$ and $\{\J(k_n)\not \subset\hat\J_n\}$, and noting that the maximum can only become smaller on a subset, we have
\begin{equation*}
\footnotesize
\begin{split}
\P\Big(\max_{j\in\hat\J_n}S_{n,j}/\sigma_j^{\tau_n}\leq t\Big) \ \leq \
 \P\Big(\max_{j\in\J(k_n)}S_{n,j}/\sigma_j^{\tau_n}\leq t\Big) +\P\big(\J_{k_n}\not\subset\hat\J_n).
\end{split}
\end{equation*}
Therefore, by subtracting from both sides the probability involving the maximum over $\J(d_n)$, we have
\footnotesize
\begin{equation}\label{eqn:firstsign}
 \P\Big(\max_{j\in\hat\J_n}S_{n,j}/\sigma_j^{\tau_n}\leq t\Big) - \P\Big(\max_{j\in\J(d_n)}S_{n,j}/\sigma_j^{\tau_n}\leq t\Big) \ \leq \ \bI_n+\P\big(\J_{k_n}\not\subset\hat\J_n).
\end{equation}
\normalsize
Similarly, by intersecting with the events $\{\hat\J\subset \J(d_n)\}$ and $\{\hat\J_n\not \subset\J(d_n)\}$, we have
\footnotesize
\begin{equation*}
\P\Big(\max_{j\in\J(d_n)}S_{n,j}/\sigma_j^{\tau_n}\leq t\Big) 
 \ \leq \  \P\Big(\max_{j\in\hat\J_n}S_{n,j}/\sigma_j^{\tau_n}\leq t\Big)+ \P(\hat\J_n\not \subset\J(d_n)).
\end{equation*}
\normalsize
 Next, subtracting the probability involving the maximum over $\hat\J_n$ gives
\footnotesize
\begin{equation}\label{eqn:secondsign}
\P\Big(\max_{j\in\J(d_n)}S_{n,j}/\sigma_j^{\tau_n}\leq t\Big) 
 \   - \ \P\Big(\max_{j\in\hat\J_n}S_{n,j}/\sigma_j^{\tau_n}\leq t\Big) \ \leq \ 0 + \P(\hat\J_n\not \subset\J(d_n)).
\end{equation}
\normalsize
Combining~\eqref{eqn:firstsign} and~\eqref{eqn:secondsign} implies
\begin{equation*}
d_{\textup{K}}\big(\L(\M)\, ,\, \L(M_{d_n})\big) \ \leq  \ \I_n+\P\big(\J_{k_n}\not\subset\hat\J_n)+\P(\hat\J_n\not \subset\J(d_n)),
\end{equation*}
and furthermore, Lemma~\ref{lem:Jhat} shows that the last two probabilities on the right are at most $c/n$. This proves ~\eqref{eqn:reductionfirst}. The inequality~\eqref{eqn:reductionsecond} follows by similar reasoning, and is actually easier, because conditioning on $X$ allows us to work under the assumption that the events $\{\J(k_n)\subset \hat\J_n\}$ and $\{\hat\J_n\subset\J(d_n)\}$ hold, since they occur with probability at least $1-\frac cn$. \qed

\section{Technical Lemmas for Theorem~\ref{THM:M}}\label{app:lemmamult}

\begin{lemma}\label{lem:Jhat}
Suppose the conditions of Theorem~\ref{THM:M} hold. Then, with probability at least $1-\frac cn$, the following two events hold simultaneously,
\begin{equation}\label{eqn:firstJhat}
\J(k_n) \ \subset \ \hat{\mathcal{J}}_n, 
\end{equation}
and
\begin{equation}\label{eqn:secondJhat}
\hat{\mathcal{J}}_n \ \subset \ \J(d_n).
\end{equation}
\end{lemma}
\proof We first address the event~\eqref{eqn:firstJhat}. By a union bound, the following inequalities hold for all large $n$,
\begin{equation}
\begin{split}
\P\Big(\J(k_n)\not\subset \hat{\mathcal{J}}_n \Big) 
& \ \leq \  \sum_{j\in\J(k_n)} \P\Big(\hat\pi_j <\sqrt{\log(n)/n}\Big)\\[0.3cm]
&  \leq \  k_n\cdot \max_{j\in\J(k_n)}\P\Big(\sqrt{n}|\hat\pi_j-\pi_j| >n^{1/4}\Big)
\end{split}
\end{equation}
where the last step follows from the fact that if $j\in\J(k_n)$, then the crude (but adequate) inequality $\sqrt{n}\pi_j- \sqrt{\log(n)}\geq n^{1/4}$ holds for all large $n$. In turn, Hoeffding's inequality~\cite[p.460]{vanderVaart:Wellner:2000} implies
$$ k_n\cdot\max_{j\in\J(k_n)}\P\Big(\sqrt{n}|\hat\pi_j-\pi_j| >n^{1/4}\Big) \ \leq \ 2k_n e^{-2\sqrt{n}} \ \lesssim \ \ts\frac{1}{n},$$
which establishes~\eqref{eqn:firstJhat}.

We now turn to~\eqref{eqn:secondJhat}. Under Assumption~\ref{A:mult}, it is simple to check that $\J(k_n)\subset\J(d_n)$ holds for all large $n$. Consequently, in the low-dimensional situation where $k_n=p$, we must have $\J(d_n)=\J(p)$ for all large $n$, and then \eqref{eqn:secondJhat} is clearly true. Hence, we may work in the situation where $k_n<p$. To  begin, observe that Assumption~\ref{A:mult} gives $\pi_{(1)}(1-\pi_{(1)})=\sigma_{(1)}^2\geq \e_0^2$, which implies $1-\pi_{(1)}\geq \e_0^2$, and hence $1-\pi_j\geq \e_0^2$ for all $j\in\{1,\dots,p\}$. Based on this observation, if we consider the set 
$$\J_n':=\Big\{j\in\{1,\dots,p\}\, \Big| \, \pi_j\geq 1/\sqrt{n}\Big\},$$
then $j\in\J_n'$ implies $\sigma_j^2=\pi_j(1-\pi_j)\geq \e_0^2/\sqrt{n}$. Now, recall that $\J(d_n)$ is defined so that $j\in\J(d_n) \Longleftrightarrow \sigma_j^2\geq \e_0^2/\sqrt{n}$. As a result of this definition, we have $\J_n'\subset\J(d_n)$. Therefore, in order to show that the event  $\{\hat\J_n\subset\J(d_n)\}$ holds with probability at least $1-c/n$, it suffices to show that the event $\{\hat\J_n\subset\J_n'\}$ holds with at least the same probability.

To proceed,
observe that the following inclusion always holds,
\begin{equation}
\begin{split}
\hat{\mathcal{J}}_n
& \ \subset \ \Big\{ j\in\{1,\dots,p\} \ \Big| \  \sqrt{n}\pi_j \ \geq \sqrt{\log(n)} - \max_{1\leq j\leq p}\sqrt{n}(\hat\pi_j-\pi_j)\Big\}.
\end{split}
\end{equation}
Therefore, if we can show that the event 
$$\mathcal{E}:=\Big\{ \sqrt{\log(n)} - \max_{1\leq j\leq p}\sqrt{n}(\hat\pi_j-\pi_j) \geq 1\Big\}$$ 
occurs with probability at least $1-\frac cn$, then the event $\{\hat\J_n\subset \J_n'\}$ will also occur with probability at least $1-\frac cn$. Now, consider the union bound
\begin{equation}
\P(\mathcal{E}^c) \ \leq \  \sum_{j=1}^p \P\Big(\sqrt{n}(\hat\pi_j-\pi_j) > \sqrt{\log(n)}-1\Big).
\end{equation}
We will bound this sum by considering two different sets of indices. For the indices $j\in\J(k_n)^c$, the values $\pi_j$ are mostly small. This motivates the use of Kiefer's inequality (Lemma~\ref{lem:Kiefer}), which implies there is some $c>0$ not depending on $n$, such that the following bound holds for all large $n$, and $j\in\J(k_n)^c$,
\begin{equation}
\begin{split}
\P\Big(\sqrt{n}(\hat\pi_j-\pi_j) > \sqrt{\log(n)}-1\Big) & \ \leq \ 2\exp\Big\{-c\log(n)\log(\ts\frac{1}{\pi_j})\Big\}\\[0.2cm]
& \ = \ 2 \pi_j^{c\log(n)}.
%
\end{split}
\end{equation}
On the other hand, if $j\in\J(k_n)$, then $\pi_j$ is of moderate size, and Hoeffding's inequality~\cite[p.460]{vanderVaart:Wellner:2000} implies the following inequality for all large $n$,
\begin{equation}
\begin{split}
\P\Big(\sqrt{n}(\hat\pi_j-\pi_j) > \sqrt{\log(n)}-1\Big) & \ \leq \ 2\exp\big\{-(3/2)\log(n)\big\}\\[0.2cm]
& \ = \ 2 n^{-3/2}.
\end{split}
\end{equation}
(The constant $3/2$ in the exponent has been chosen for simplicity, and is not of special importance.) Combining the two different types of bounds, and using the fact that any probability vector satisfies $\pi_{(j)}\leq j^{-1}$ for all $j\in\{1,\dots,p\}$, we obtain
\begin{equation}
\begin{split}
\P(\mathcal{E}^c) 
& \ \lesssim \ k_n n^{-3/2} \ + \ \sum_{j=k_n+1}^p \pi_{(j)}^{c\log(n)}\\[0.3cm]
& \ \lesssim \ \ts\frac{1}{n} \ +  \displaystyle\int_{k_n}^p x^{-c\log(n)}dx\\[0.3cm]
& \ \lesssim \ \ts\frac{1}{n} \ + k_n^{-c\log(n)+1}\\[0.3cm]
& \ \lesssim \ \ts\frac{1}{n},
\end{split}
\end{equation}
as needed.\qed

\begin{lemma}\label{lem:sighatrnormM}
Suppose the conditions of Theorem~\ref{THM:M} hold, and let 
$q=\max\{\ts\frac{2}{(1/2)(1-\tau_n)}, \log(n),3\}$.
Then, there is a constant $c>0$ not depending on $n$, such that for any $j\in\mathcal{J}(d_n)$, we have
\begin{equation}
\|\hat\sigma_{j}\|_q \ \leq \  c \cdot \sigma_j \cdot \sqrt{q}.
\end{equation}
\end{lemma}
\proof By direct calculation
\begin{equation*}
\begin{split}
\ts\frac{1}{\sigma_j}\|\hat\sigma_{j}\|_q
&= \ \ts\frac{1}{\sigma_j} \Bigg\|\ts\frac 1n \sum_{i=1}^n (X_{i,j}^2-\pi_j)+\big(\pi_j-\pi_j^2\big)+\big(\pi_j^2-\bar{X}_j^2\big)\Bigg\|_{q/2}^{1/2}\\[0.2cm]
&\leq \ \ts\frac{1}{\sigma_j} \Bigg(\big\|\ts\frac 1n \sum_{i=1}^n (X_{i,j}^2-\pi_j)\big\|_{q/2}^{1/2}+(\pi_j-\pi_j^2)^{1/2}+\big\|\pi_j^2-\bar{X}_j^2\big\|_{q/2}^{1/2}\Bigg)\\[0.3cm]
&\leq \ \ts\frac{1}{\sigma_jn^{1/4}}\big\| S_{n,j}\big\|_{q/2}^{1/2} \ + \ 1 \ + \ \ts\frac{\sqrt{2}}{\sigma_jn^{1/4}}\big\| S_{n,j}\big\|_{q/2}^{1/2}.
\end{split}
\end{equation*}
Since $1/(\sigma_j n^{1/4})\leq 1/\e_0$ when $j\in\J(d_n)$, it follows from Lemma~\ref{lem:rosenthal} that the first and third terms are at most of order $\sqrt{q}$. \qed

\paragraph{Remark} For the next lemma, put $\boldsymbol \pi_{k_n}:=(\pi_{(1)},\dots,\pi_{(k_n)})$, and recall the definition $\Sigma(k_n)=\Pi_{k_n}\Sigma\Pi_{k_n}\ttop$ from line~\eqref{eqn:sigmakdef}.

\begin{lemma}\label{lem:lambdamin}
If $r_n$ denotes the rank of $\Sigma(k_n)$, and the conditions of Theorem~\ref{THM:M} hold, then there is a constant $c>0$ not depending on $n$ such that
\begin{equation}
\lambda_{r_n}(\Sigma(k_n)) \ \gtrsim \ k_n^{-c}.
\end{equation}
\end{lemma}
\proof  We first consider the case $k_n<p$, and then handle the case $k_n=p$ separately at the end of the proof. 
 Under the multinomial model, we have for all $i,j\in\{1,\dots,k_n\}$,
$$\Sigma_{i,j}(k_n) = \begin{cases} \  \pi_{(i)}(1-\pi_{(i)}) & \text{ if } i=j, \\ \ -\pi_{(i)}\pi_{(j)} & \text{ if } i\neq j. \end{cases}$$
For each $i\in\{1,\dots,k_n\}$, define the ``deleted row sum'', 
$$\varrho_i:=\sum_{j\neq i}|\Sigma_{i,j}(k_n)|.$$
By the Ger\v{s}gorin disc theorem~\citep[Sec. 6.1]{Horn:Johnson:1990}, 
\begin{equation}
\begin{split}
\lambda_{r_n}(\Sigma(k_n)) & \ \geq \ \lambda_{\min}(\Sigma(k_n))\\[0.2cm]
 & \ \geq \  \min_{1\leq i\leq k_n}\Big\{\Sigma_{i,i}(k_n)  - \varrho_i\Big\}\\[0.2cm]
& \ = \ \min_{1\leq i\leq k_n} \Big\{(\pi_{(i)} -\pi_{(i)}\tsum_{j=1}^{k_n}\pi_{(j)}\Big\},\\[0.2cm]
& \ = \ \pi_{(k_n)}\Big(1-\tsum_{j=1}^{k_n} \pi_{(j)}\Big).
\end{split}
\end{equation}
When $k_n<p$, it follows from Assumption~\ref{A:mult} that
\begin{equation}
\begin{split}
1-\tsum_{j=1}^{k_n}\pi_{(j)}&\geq \pi_{(k_n+1)}\\[0.2cm]
&\geq \sigma_{(k+1)}^2\\[0.2cm]
&\geq \e_0^2 (k_n+1)^{-2\alpha}\\[0.2cm]
&\gtrsim k_n^{-2\alpha}.
\end{split}
\end{equation}
Hence, the previous steps show that $\lambda_{r_n}(\Sigma(k_n))$ is at least of order $k_n^{-4\alpha}$.

Finally, consider the case when $k_n=p$. In this case, it is a basic fact that the rank of $\Sigma=\Sigma(k_n)$ satisfies $r_n=p-1$ (where we note that  Assumption~\ref{A:mult} ensures $\pi_{(p)}>0$). Also, it is known from matrix analysis~\citep[Theorem 1]{Benasseni:2012} that 
$$\lambda_{p-1}(\Sigma)\geq \pi_{(p)},$$
and therefore Assumption~\ref{A:mult} leads to
\begin{equation}
\begin{split}
\lambda_{r_n}(\Sigma(k_n))&\geq \pi_{(k_n)}\\[0.2cm]
&\geq\sigma_{(k_n)}^2\\[0.2cm]
&\geq \e_0^2 k_n^{-2\alpha},
\end{split}
\end{equation}
as needed.\qed

\begin{lemma}\label{lem:whitemult} 
Let the deterministic matrix $U\in\R^{r_n\times r_n}$ and the random matrix $W_n\in\R^{r_n\times r_n}$ be as defined in~\eqref{eqn:Udef} and~\eqref{eqn:whitemultdef} respectively. Also, suppose that the conditions of Theorem~\ref{THM:M} hold. Then, there is a constant $c>0$ not depending on $n$ such that the event
\begin{equation}
\big\|U\ttop W_n U-\mathbf{I}_{r_n}\big\|_{\textup{op}}\ \leq \ \ts\frac{c\, k_n^c \log(n)}{\sqrt n}
\end{equation}
holds with probability at least $1-\frac cn$.
\end{lemma}
\proof Let the notation from the proof of Lemma~\ref{lem:IIprimeM} be in force, and observe that
\footnotesize
\begin{equation*}
\begin{split}
\big\|U\ttop W_n U-\mathbf{I}_{r_n}\big\|_{\textup{op}} & \ = \ \Big\|U\ttop\Big(\Lambda_{r_n}^{-1/2}Q\ttop \hat\Sigma(k_n) Q\Lambda_{r_n}^{-1/2} - \mathbf{I}_{r_n}\Big)U\Big\|\op\\[0.2cm]
& \ \leq \ \Big\|\Lambda_{r_n}^{-1/2}Q\ttop \Big(\hat\Sigma(k_n) -\Sigma(k_n)\Big)Q\Lambda_{r_n}^{-1/2}\Big\|\op\\[0.2cm]
& \ \leq \ \frac{1}{\lambda_{r_n}(\Sigma(k_n))} \, \big\|\hat\Sigma(k_n)-\Sigma(k_n)\big\|\op.
\end{split}
\end{equation*}
\normalsize
With regard to the first factor in the previous line, Lemma~\ref{lem:lambdamin} implies
$$\ts\frac{1}{\lambda_{r_n}(\Sigma(k_n))} \lesssim k_n^c.$$
To complete the proof, let $\boldsymbol\pi_{k_n}=\Pi_{k_n}\boldsymbol \pi$, and let $u\in\R^{k_n}$ be a generic unit vector. Consider the decomposition
\begin{equation*}
\footnotesize
\begin{split}
|u\ttop\big(\hat\Sigma(k_n)-\Sigma(k_n)\big)u| & \leq  \bigg|\frac1n\sum_{i=1}^n \Big((X_i\ttop \Pi_{k_n}\ttop u)^2-u\ttop\text{diag}(\boldsymbol\pi_{k_n})u\Big)\bigg| + \Big| (\bar X\ttop \Pi_{k_n}\ttop u)^2-(\boldsymbol\pi_{k_n}\ttop u)^2\Big|\\[0.3cm]
& \ =: \ \Delta_n(u)+\Delta_n'(u).
\end{split}
\end{equation*}
\normalsize
In order to control these terms, note that each random variable $(X_i\ttop \Pi_{k_n}\ttop u)^2$ is bounded in magnitude by 1, and has expectation equal to $u\ttop \text{diag}(\boldsymbol \pi_{k_n})u$. In addition, we have
$$|\Delta_n'(u)| \ \leq \ 2|\bar X\ttop \Pi_{k_n}\ttop u-\boldsymbol \pi_{k_n}\ttop u|.$$
Based on these observations, the proof of Lemma~\ref{lem:white} can be essentially repeated to show that there is a constant $c>0$ not depending on $n$ such that the event
$$ \big\|\hat\Sigma(k_n)-\Sigma(k_n)\big\|_{\text{op}} \ \leq \ \ts\frac{c\log(n) k_n}{\sqrt{n}}$$
holds with probability at least $1-\frac cn$. This completes the proof.\qed\\

\begin{lemma}\label{lem:SnjnormM} Let $q=\max\{\ts\frac{2}{(1/2)(1-\tau_n)}, \log(n),3\}$, and suppose that the conditions of Theorem~\ref{THM:M} hold.  Then, there is a constant $c>0$ not depending on $n$, such that 
\begin{equation}
\max_{j\in\mathcal{J}(d_n)}\|\ts\frac{1}{\sigma_j}S_{n,j}\|_q \leq cq.
\end{equation}
In addition, the following event holds with probability 1,
\begin{equation}
\max_{j\in\mathcal{J}(d_n)} \Big(\E\big[|\ts\frac{1}{\hat \sigma_j}S_{n,j}^{\star}|^q|X\big] \Big)^{1/q} \ \leq c\, q.
\end{equation}
\end{lemma}
\proof The second inequality can be obtained by repeating the proof of Lemma~\ref{lem:Snjnorm}, since $S_n^{\star}$ is still Gaussian under the setup of Assumption~\ref{A:mult}. To prove the first inequality, note that since $q>2$, Lemma~\ref{lem:rosenthal} gives
\begin{equation}\label{eqn:lemrosenthalfirst}
 \|\ts\frac{1}{\sigma_j}S_{n,j}\|_q \ \lesssim \ q\cdot \max\Big\{  \|\ts\frac{1}{\sigma_j}S_{n,j}\|_2 \, , \,  n^{-1/2+1/q} \|\ts\frac{1}{\sigma_j}(X_{1,j}-\pi_j)\|_q\Big\}.
\end{equation}
Clearly,
\begin{equation*}
 \|\ts\frac{1}{\sigma_j}S_{n,j}\|_2^2 \ = \  \var(\ts\frac{1}{\sigma_j}S_{n,j}) \ = \ 1.
\end{equation*}
For the stated choice of $q$, the second term inside the maximum satisfies
$$n^{-1/2+1/q} \|\ts\frac{1}{\sigma_j}(X_{1,j}-\pi_j)\|_q \ \lesssim \frac{1}{\sqrt n \sigma_j},$$
and also, since $j\in\mathcal{J}(d_n)$, we have $\frac{1}{\sqrt{n}\sigma_j}\lesssim \frac{1}{n^{1/4}}$, which leads to the stated claim. \qed

\begin{lemma}\label{lem:corM}
Suppose the conditions of Theorem~\ref{THM:M} hold. Then, there is a constant $c>0$ not depending on $n$ such that the events
\begin{equation}\label{eqn:firstcor}
 \max_{j\in\J(k_n)}\Big| \ts\frac{\hat\sigma_j}{\sigma_j}-1\Big| \leq \ts \frac{c\cdot k_n^c\cdot\sqrt{\log(n)}}{n^{1/2}},
\end{equation}
and
\begin{equation}\label{eqn:minvarhattemp}
 \min_{j\in\J(k_n)}\hat\sigma_j^{1-\tau_n} \ \geq \ \Big(\min_{j\in\J(k_n)} \sigma_j^{1-\tau_n}\Big)\cdot \Big(1-\ts\frac{c\cdot k_n^c\cdot \sqrt{\log(n)}}{n^{1/2}}\Big)
\end{equation}
each
hold with probability at least $1-\ts\frac{c}{n}$.
\end{lemma}

\proof  Note that if~\eqref{eqn:firstcor} occurs, then~\eqref{eqn:minvarhattemp} also occurs, and so we only deal with the former event. Fix any number $\kappa\geq 2$. By a union bound, it suffices to show there are positive constants $c$ and $c_1(\kappa)$ not depending on $n$, such that for any $j\in\J(k_n)$, the event
\begin{equation}\label{eqn:corevent}
 \Big| \ts\frac{\hat{\sigma}_j^2}{\sigma_j^2}-1\Big| \ \leq \ \frac{c_1(\kappa)\cdot k_n^c\cdot \sqrt{\log(n)}}{n^{1/2}}
\end{equation}
holds with probability at least $1-\frac{c}{n^{\kappa}}$. 
To this end, observe that
\begin{equation*}
\begin{split}
\Big| \ts\frac{\hat{\sigma}_j^2}{\sigma_j^2}-1\Big| & \ \leq \ \Big|\ts\frac{1}{\sigma_j^2n} \sum_{i=1}^n (X_{i,j}^2-\pi_j) \Big| \ + \ \ts\frac{1}{\sigma_j^2}\big|\pi_j^2-\bar{X}_j^2\big|\\[0.3cm]
  & \ \leq \  \ts\frac{3}{\sigma_j^2 \sqrt{n}}|S_{n,j}|.
 \end{split}
\end{equation*}
Due to Assumption~\ref{A:mult}, and the fact that $j\in\J(k_n)$, we have 
$$\ts\frac{1}{\sigma_j^2\sqrt n} \ \lesssim  \ \ts\frac{k_n^{2\alpha}}{\sqrt n}.$$
In addition, Hoeffding's inequality ensures there is a constant $c_1(\kappa)$ such that the event
\begin{equation*}
|S_{n,j}| \ \leq \ c_1(\kappa)\sqrt{\log(n)}
\end{equation*}
holds with probability at least $1-\frac{c}{n^{\kappa}}$.\qed


\section{Background results} \label{app:background}
The following result is a multivariate version of the Berry-Esseen theorem due to~\cite{Bentkus:2003}. 
\begin{lemma}[Bentkus' multivariate Berry-Esseen theorem]\label{lem:bentkus}
\hspace{-0.4cm} Let\ $V_1,\dots,V_n$ be i.i.d.~random vectors $\R^d$, with zero mean, and identity covariance matrix. Furthermore, let $\gamma_d$ denote the standard Gaussian distribution on $\R^d$,  and let $\mathscr{A}$ denote the collection of all Borel convex subsets of $\R^d$. Then, there is an absolute constant $c>0$ such that
\begin{equation}
 \sup_{\mathcal{A}\in\mathscr{A}} \Big| \P\big(\ts\frac{1}{\sqrt{n}}(V_1+\dots+V_n)\in \mathcal{A}\big) -\gamma_d(\mathcal{A})\Big| \ \leq  \ \displaystyle \frac{c\cdot  d^{1/4}\cdot  \E\big[\|V_1\|_2^3\big]}{n^{1/2}}.
\end{equation}
\end{lemma}

\noindent The following is a version of Nazarov's inequality~\citep{Nazarov:2003,Klivans:2008}, as formulated in~\cite[Lemma 4.3]{CCK:SPA}.
\begin{lemma}[Nazarov's inequality]\label{lem:Nazarov}
Let $(\xi_1,\dots,\xi_m)$ be a multivariate normal random vector, and suppose the parameter $\underline{\sigma}^2:=\min_{1\leq j\leq m} \var(\xi_j)$ is positive. Then, for any $r>0$,
\begin{equation}
\sup_{t\in\R}\, \P\Big(\Big|\max_{1\leq j\leq m} \xi_j -t\Big| \ \leq r \Big)  \ \leq \  \frac{2r}{\underline{\sigma}}\cdot (\sqrt{2\log(m)}+2).
\end{equation}
\end{lemma}

\noindent The result below is a version of Slepian's lemma, which is adapted from~\cite[Theorem 2.2]{Li:Shao:2002}. (See references therein for earlier versions of this result.)
\begin{lemma}[Slepian's lemma]\label{lem:slepian} 
Let $m\geq 3$, and let $\mathsf{R}\in\R^{m\times m}$ be a correlation matrix with $\max_{i\neq j}\mathsf{R}_{i,j}<1$. Also, let $\mathsf{R}^+$ be the matrix with $(i,j)$ entry given by $\max\{\mathsf{R}_{i,j},0\}$, and suppose $\mathsf{R}^+$ is positive semi-definite. Furthermore, let $\zeta\sim N(0,\mathsf{R})$ and $\xi\sim N(0,\mathsf{R}^+)$. Then, the following inequalities hold for any $t\geq 0$,
\begin{equation}
 \P\Big(\max_{1\leq j\leq m} \zeta_j  \leq t\Big) \ \leq \ \P\Big(\max_{1\leq j\leq m} \xi_j\leq t\Big) \ \leq \  K_m(t) \cdot \Phi^m(t),
\end{equation}
where
\begin{equation}
K_m(t)=\exp\Bigg\{\sum_{1\leq i<j\leq m} \log\Big(\ts\frac{1}{1-\ts\frac{2}{\pi}\arcsin(\mathsf{R}_{i,j}^+)}\Big)\exp\Big(-\ts\frac{t^2}{1+\mathsf{R}_{i,j}^+}\Big)\Bigg\}.
\end{equation}

\end{lemma}

\noindent The following inequalities are due to~\cite{Zinn:1985}. 
\begin{lemma}[Rosenthal's inequality with best constants]\label{lem:rosenthal}
Fix $r\geq 1$ and put $\textup{Log}(r):=\max\{\log(r),1\}$. Let $\xi_1,\dots,\xi_m$ be independent random variables satisfying $\E[|\xi_j|^r]<\infty$ for all $1\leq j\leq m$. Then, there is an absolute constant $c>0$ such that the following two statements are true.
\begin{itemize}
\item[(i).] If $\xi_1,\dots,\xi_m$ are non-negative random variables, then 
	\begin{equation}
	\big\|\tsum_{j=1}^m \xi_i\big\|_r\leq c\cdot \ts\frac{r}{\textup{Log}(r)} \cdot \max\bigg\{ \big\|\tsum_{j=1}^m \xi_j\big\|_1 \, , \, \big(\tsum_{j=1}^m \big\|\xi_i\|_r^r\big)^{1/r}\bigg\}.
	\end{equation}
\item[(ii).] If $r>2$, and the random variables $\xi_1,\dots,\xi_m$ all have mean 0, then
	\begin{equation}
	\big\|\tsum_{j=1}^m \xi_i\big\|_r\leq c\cdot \ts\frac{r}{\textup{Log}(r)} \cdot \max\bigg\{ \big\|\tsum_{j=1}^m \xi_j\big\|_2 \, , \, \big(\tsum_{j=1}^m \big\|\xi_i\|_r^r\big)^{1/r}\bigg\}.
	\end{equation}
	\end{itemize}
\end{lemma}

\paragraph{Remark} The non-negative case is handled in~\cite[Theorem 2.5]{Zinn:1985}.  With regard to the mean 0 case, the statement above differs slightly from \cite[Theorem 4.1]{Zinn:1985}, which requires symmetric random variables, but the 
 remark on page 247 of that paper explains why the variables $\xi_1,\dots,\xi_m$ need not be symmetric as long as they have mean 0.\\

The result below is a sharpened version of Hoeffding's inequality for handling the binomial distribution when the success probability is small~\citep[Corollary A.6.3]{vanderVaart:Wellner:2000}.
\begin{lemma}[Kiefer's inequality]\label{lem:Kiefer}
Let $\xi_1,\dots,\xi_m$ be independent Bernoulli random variables with success probability $\pi_0\in(0,1/e)$, and let $\bar\xi=\frac{1}{m}\sum_{i=1}^m\xi_i$. Then, the following inequality holds for any $m\geq 1$ and $t>0$,
\begin{equation}
\P\Big(\sqrt{m}|\bar\xi-\pi_0|\geq t\Big) \ \leq \ 2\exp\Big\{-t^2\big[\log\big(\ts\frac{1}{\pi_0}\big)-1\big]\Big\}.
\end{equation}
\end{lemma}

\section{Related work on Gaussian approximation}\label{app:Gaussian}

Although our focus is primarily on rates of bootstrap approximation, this topic is closely related to rates of Gaussian approximation in the central limit theorem --- for which there is a long line of work in finite-dimensional and infinite-dimensional settings. We refer to the chapter \citep{Bentkus:2000} for a general survey, as well as Appendix L of the paper~\citep{CCK:2013} for a discussion that is oriented more towards high-dimensional statistics.

 To describe how our work fits into this literature, we fix some notation. Let $\mathbb{B}$ denote a Banach space, with $\mathscr{A}$ being a collection of its subsets, and $\varphi:\mathbb{B}\to\R$ being a function. For the moment, we will regard the summands of $S_n$ as centered~i.i.d.~random elements of $\mathbb{B}$, and let $\tilde S_n$ denote a centered Gaussian element of $\mathbb{B}$ with the same covariance as $S_n$. Broadly speaking, the literature on rates of Gaussian approximation is concerned with bounding the  quantities 
 $\rho_n(\mathscr{A})=\sup_{A\in\mathscr{A}}\rho_n(A)$, or $\|\Delta_n\|_{\infty}=\sup_{t\in\R}|\Delta_n(t)|$, where
 \begin{align}
 \rho_n(A) & \ = \  \Big|\P(S_n\in A) - \P(\tilde S_n\in A)\Big|,\label{eqn:rhodef}\\[0.2cm]
 \Delta_n(t) & \ = \ \Big|\P\big(\varphi(S_n)\leq t\big) - \P\big(\varphi(\tilde S_n)\leq t\big)\Big|.\label{eqn:deltadef}
 \end{align}
  In particular, note that if $\mathbb{B}=\R^p$ and $\tilde T=\max_{1\leq j\leq p} \tilde S_{n,j}$, then the distance $d_{\text{K}}(\mathcal{L}(T),\mathcal{L}(\tilde T))$ can be represented as either $\rho_n(\mathscr{A})$ or $\|\Delta_n\|_{\infty}$, by taking $\mathscr{A}$ to be a class of rectangles, or by taking $\varphi$ to be the coordinate-wise maximum function.

   Typically, the rates at which $\rho_n(\mathscr{A})$ and $\|\Delta_n\|_{\infty}$ decrease with $n$ is dependent on the distribution of $S_n$, the dimension of $\mathbb{B}$, as well as the smoothness of $\mathscr{A}$ and $\varphi$, among other factors. Although the study of $\rho_n(\mathscr{A})$ and $\|\Delta_n\|_{\infty}$ is highly multifaceted, it is a general principle that the smoothness of $\mathscr{A}$ and $\varphi$ tends to be much more influential when $\mathbb{B}$ is infinite-dimensional, as compared to the finite-dimensional case. (A discussion may be found in~\citep{Bentkus:1993}.) Indeed, this is worth emphasizing in relation to our work, since the choices of $\mathscr{A}$ and $\varphi$ corresponding to $d_{\text{K}}(\mathcal{L}(T),\mathcal{L}(\tilde T))$ are notable for their lack of smoothness.
   
   To illustrate this point, consider the following two facts that are known to hold when $\mathbb{B}$ is separable and infinite-dimensional~\citep{Bentkus:1993}: (1) If $\E[\|X_1\|_{\mathbb{B}}^3]<\infty$, and $\varphi$ satisfies smoothness conditions stronger than having three Fr\'echet derivatives,  then $\|\Delta_n\|_{\infty}$ is of order $\mathcal{O}(n^{-1/2})$. (2) If the previous conditions hold, except that $\varphi$ is only assumed to have one Fr\'echet derivative, then an example can be constructed so that $\|\Delta_n\|_{\infty}$ is bounded below by a sequence that converges to 0 arbitrarily slowly.  (Other examples of lower bounds may be found in the papers~\citep{Bentkus:1986,Talagrand:1984}, among others.) By contrast, in the finite-dimensional case where $\mathbb{B}=\R^p$ with $p$ held fixed as $n\to\infty$ and $\E[\|X_1\|_2^3]<\infty$, it is known that  non-smooth choices of $\varphi$ and $\mathscr{A}$ can lead to a $n^{-1/2}$ rate. Namely, it is known that $\|\Delta_n\|_{\infty}$ and $\rho_n(\mathscr{A})$ are of order $\mathcal{O}(n^{-1/2})$ when $\varphi$ is convex, or when $\mathscr{A}$ is the class of Borel convex sets. Beyond the case where $p$ is held fixed, many other works allow $n$ and $p$ to diverge together. When $p$ grows relatively slowly compared to $n$, the leading rate for the choices of $\varphi$ and $\mathscr{A}$ just mentioned is essentially $\mathcal{O}(p^{7/4}/n^{1/2})$~\citep{Bentkus:2003,Bentkus:2005}. See also~\citep{Sazonov:1968,Nagaev:1976,Senatov:1981,Portnoy:1986,Gotze:1991,Chen:2011,Zhai:2018} for additional background. Meanwhile, when $p\gg n$, it has recently been established that rates of the form $\mathcal{O}(\log(p)^b/n^{1/6})$ can be achieved if $\varphi$ is the coordinate-wise maximum function, or if $\mathscr{A}$ is a class of ``sparsely convex sets'', such as rectangles~\citep{CCK:AOP}.

   In light of the previous paragraph, it is natural that results in the infinite-dimensional setting have focused predominantly on smooth choices of $\mathscr{A}$ and $\varphi$. Nevertheless, there are some special cases where notable results have been obtained for non-smooth choices of $\varphi$ and $\mathscr{A}$ that correspond to max statistics. One such result is established in the paper~\citep{Asriev:1986}, which deals with bounds on $\rho_n(A(\mathbf{r}))$ for rectangular sets such as $A(\mathbf{r}):= \prod_{j=1}^{\infty}(-\infty,r_j]$,  with $\mathbf{r}=(r_1,r_2,\dots)$, in the infinite-dimensional Euclidean space $\R^{\infty}$. More specifically, if we continue to let $\sigma_j^2=\var(X_{1,j})$ for $j=1,2,\dots$, and define the ``effective dimension'' parameter
   $ {\tt{d}}(\mathbf{r})=\sum_{j=1}^{\infty} \ts\frac{\sigma_j}{\sigma_j+r_j},$
   then the following bound holds under certain conditions on the distribution of $S_n$,
   \begin{equation}\label{eqn:asriev}
    \rho_n(A(\mathbf{r})) \ \lesssim \  \ts\frac{\log(n)^{3/2}}{n^{1/2}} \cdot {\tt{d}}(\mathbf{r})\cdot(1+\tt{d}(\mathbf{r})^2).
   \end{equation}
To comment on how the bound~\eqref{eqn:asriev} relates with our Gaussian approximation result in Theorem~3.1, note that they both involve near-parametric rates, and are governed by the parameters $(\sigma_1,\sigma_2,\dots)$. However, there are also some crucial differences. First, the bound~\eqref{eqn:asriev} is \emph{non-uniform} with respect to the set $A(\mathbf{r})$, whereas Theorem~3.1 is a uniform result. In particular, this difference becomes apparent by setting all $r_j$ equal to $\sigma_{(n)}$, which implies ${\tt{d}}(\mathbf{r})\geq n/2$, and causes the bound~\eqref{eqn:asriev} diverge as $n\to\infty$.  Second, the bound~\eqref{eqn:asriev} relies on the assumption that $S_n$ has a \emph{diagonal covariance} matrix, whereas Theorem~3.1 allows for much more general covariance structures.

Some further examples related to max statistics arise in the context of empirical process theory. Let $\G_n(f)=\sqrt{n}\sum_{i=1}^n (f(X_i)-\E[f(X_i)])$ denote an empirical process that is generated from i.i.d.~observations $X_1,\dots,X_n$, and indexed by a class of functions $f\in\mathscr{F}$. Also let the Gaussian counterpart of $\G_n$ be denoted by $\tilde{\G}$ (i.e.~a Brownian bridge)~\citep[][]{vanderVaart:Wellner:2000}. In this setting, the paper~\citep{Norvaivsa:1991} studies the quantity
\begin{equation}\label{eqn:Deltanew}
\Delta_n(t)=\Big|\P\big(\|\G_n\|_{\mathscr{F}}\leq t\big) - \P\big(\|\tilde\G\|_{\mathscr{F}}\leq t\big)\Big|,
\end{equation}
which can be understood in terms of the earlier definition~\eqref{eqn:deltadef} by letting $\varphi(\G_n)=\|\G_n\|_{\mathscr{F}}=\sup_{f\in\mathscr{F}}|\G_n(f)|$. Under the assumption that $\mathscr{F}$ is a VC subgraph class of uniformly bounded functions, it is shown in the paper~\citep{Norvaivsa:1991} that the bound
\begin{equation}\label{eqn:empbound1}
\Delta_n(t) \ \lesssim \ \ts\frac{1}{(1+t)^3}\ts\frac{\log(n)^2}{n^{1/6}}
\end{equation}
holds for all $t\geq 0$. In relation to the current work, the condition that the functions in $\mathscr{F}$ are uniformly bounded is an important distinction, because the max statistic $T=\max_{1\leq j\leq p} \G_n(f_j)$ arises from the functions $f_j(x)=x_j$, which are unbounded. 

 More recently, the papers~\citep{CCK:suprema,CCK:SPA} have derived bounds on the quantity~$\|\Delta_n\|_{\infty}$ (as well as coupling probabilities) under weaker assumptions on the class $\mathscr{F}$. For instance, these works allow for classes of functions that are non-Donsker or unbounded, provided that a suitable envelope function is available. Nevertheless, the resulting bounds on $\|\Delta_n\|_{\infty}$ involve restrictions on how quickly the parameter $\inf_{f\in\mathscr{F}}\var(\G_n(f))$ can decrease with $n$ --- whereas our results do not involve such restrictions. Also, the rates developed in these works are broadly similar to~\eqref{eqn:empbound1}, and it is unclear if a modification of the techniques would lead to near-parametric rates in our setting.
Lastly, a number of earlier results that are related to bounding $\|\Delta_n\|_{\infty}$, such as invariance principles and couplings, can be found in the papers~\citep[][and references therein]{Csorgo:1986,Massart:1986,Massart:1989,Paulauskas:Stieve:1990,Bloznelis:1997}. However, these works are tailored to fairly specialized processes (e.g., dealing with c\`adl\`ag functions on the unit interval, rectangles in $\R^p$ when $p$ is fixed, or processes generated by Uniform[0,1] random variables).

\end{document}